\documentclass[12pt]{amsart}
\usepackage{amsmath,amsfonts,amssymb,amscd,amsthm,epsf}
\usepackage{comment}
\newtheorem{theorem}{Theorem}
\newtheorem{lemma}[theorem]{Lemma}
\newtheorem{prop}[theorem]{Proposition}
\newtheorem{de}[theorem]{Definition}
\theoremstyle{remark}
\newtheorem{remark}[theorem]{Remark}
\newtheorem{algorithm}[theorem]{Algorithm}

\def\real{{\mathbb R}}
\def\torus{{\mathbb T}}
\def\integer{{\mathbb Z}}

\def\A{{\mathcal{A}}}
\def\B{{\mathcal{B}}}
\def\C{{\mathcal{C}}}

\def\J{{\mathcal{J}}}
\def\K{{\mathcal{K}}}
\def\L{{\mathcal{L}}}

\def\R{{\mathcal{R}}}
\def\S{{\mathcal{S}}}

\def\j{m}
\def\a{a(\hat J)}

\def\ep{\varepsilon}
\def\Id{ {\rm Id}}
\def\NHIM{normally hyperbolic invariant manifolds}

\def\Y{ {\mathcal Y}}
\def\M{ {\mathcal M}}

\def\dist{\operatorname{dist}}

\newcommand{\abs}[1]{\left | #1 \right |}
\newcommand{\norm}[1]{\left \| #1 \right \|}

\newcommand{\Or}{\mathop{\rm O}\nolimits}

\renewcommand{\d}{\,{\rm d}}
\newcommand{\cte}{{\rm cte.}\,}

\newcommand{\st}{{\rm s}}           
\newcommand{\un}{{\rm u}}           
\def\cal{\mathcal}
\def\ep{\varepsilon}

\def\itr{\pitchfork}
\def\Tau{ {\mathcal T} }

\def\A{ {\cal A}}
\def\B{ {\cal B}}
\def\D{ {\cal D}}

\def\I{ {\cal I}}
\def\K{ {\cal K}}
\def\LL{ {\cal L}}
\def\N{ {\cal N}}

\def\T{ {\cal T}}
\def\U{ {\cal U}}

\def\C{{\cal C}}

\def\tLambda{\tilde \Lambda}
\def\vp{\varphi}

\def\II{\cal I}
\def\JJ{\cal J}

\date{June 11, 2013}
\begin{document}
\title{Instability of high dimensional Hamiltonian
Systems: Multiple resonances do not impede diffusion. }
\author{Amadeu Delshams}
\email{Amadeu.Delshams@upc.es}
\address{Departament de Matem\`{a}tica Aplicada I\\
           Universitat Polit\`{e}cnica de Catalunya\\
           Diagonal 647, 08028 Barcelona, Spain
}
\author{Rafael de la Llave}
\address{School of Mathematics\\
Georgia Inst. of Technology\\
686 Cherry St. \\
Atlanta, GA 30332} \email{rafael.delallave@math.gatech.edu}
\author{ Tere M.~Seara}
\address{Departament de Matem\`{a}tica Aplicada I\\
           Universitat Polit\`{e}cnica de Catalunya\\
           Diagonal 647, 08028 Barcelona, Spain
} \email{Tere.M-Seara@upc.es}

\begin{abstract}
We consider models given by Hamiltonians of the form
\[
H(I, \vp, p, q, t; \ep)
= h(I) + \sum_{j = 1}^n \pm\left( \frac{1}{2} p_j^2 + V_j(q_j) \right)
+ \ep Q(I, \vp, p, q, t; \ep)
\]
where $I \in \II \subset \real ^{d}, \vp \in \torus^d$, $p, q \in \real^n$,
$t \in \torus^1$. These are higher dimensional analogues, both in the center
and hyperbolic directions, of the models studied in
\cite{DelshamsLS03a,DelshamsLS06a,GideaL06a, GideaL06b}. All these models
present the \emph{large gap problem}.

We show that, for $0 < \ep \ll 1$, under regularity and explicit
non-degeneracy conditions on the model,  there are orbits whose action
variables $I$ perform rather arbitrary excursions in a domain of size $O(1)$.
This domain  includes resonance lines and, hence, large gaps among
$d$-dimensional  KAM tori.

The method of proof follows closely the strategy of
\cite{DelshamsLS03a,DelshamsLS06a}. The main new phenomenon that appears when
the dimension $d$ of the center directions is larger than one, is the
existence of multiple resonances. We show that, since these multiple
resonances happen in sets of codimension greater  than one in the space of
actions $I$, they can be contoured. This corresponds to the mechanism called
\emph{diffusion across resonances} in the Physics literature.

The present paper, however, differs substantially from
\cite{DelshamsLS03a,DelshamsLS06a}. On the technical details of the proofs,
we  have taken advantage of the theory of the  scattering map
\cite{DelshamsLS08}, not available when the above papers were written. We
have analyzed the  conditions imposed on the resonances in more detail.

More precisely, we have found that there is a simple condition on the
Melnikov potential which allows us to conclude that the resonances are
crossed. In particular, this condition does not depend on the resonances. So
that the results are new even when applied to the models in
\cite{DelshamsLS03a,DelshamsLS06a}.
\end{abstract}

\maketitle \tableofcontents

\section{Introduction}
The goal of this paper is to present some explicit sufficient conditions for
instability in models of the form
\begin{equation} \label{modelsconsidered}
H(I, \vp, p, q, t; \ep) = h(I) + P(p,q) + \ep Q(I, \vp, p, q, t; \ep),
\end{equation}
where
\[
P(p,q)=\sum_{j = 1}^n P_j(p_j, q_j), \quad P_j(p_j, q_j) = \pm
\left(\frac{1}{2} p_j^2 + V_j(q_j) \right).
\]

We will assume that  $I \in \II \subset \real ^{d}, \vp \in \torus^d$, $\II$
an open set, $p, q \in \real^n$. The symplectic form of the phase space is
$\Omega = \sum_i d I_i \wedge d\vp _i + \sum_j \d p_j \wedge  d q_j$. We will
assume that the dependence of $Q$ on $t$ is $1$-periodic, so that $t \in
\torus^1$. Moreover, for simplicity, we will assume that $Q$ is a
trigonometric polynomial in the variables $(\vp, t)$. This is not a crucial
assumption and can be eliminated. See remark \ref{trigonometricpol} for
details about this assumption.

We will show (see Theorem~\ref{main} for precise statements) that, under
suitable regularity and non-degeneracy assumptions---that can be checked by
studying the $3$ jet in the $\ep$ variable of the perturbation $\ep
Q$---there are orbits of the system in which the $I$ variables can perform
largely arbitrary excursions in a set $\II^* \subset \II$ of size of order
$1$ (that is, independent of $\ep$ as $\ep \to 0$). An explicit example is
shown in Section~\ref{sec:example}.

The main point is that the set $\II^*$ can include simple resonant surfaces,
that is, resonances of multiplicity one. This makes the models considered
here present the \emph{large gap problem}. This problem, which will be
discussed in more detail in Section~\ref{averaged}, consists in the fact that
the customary perturbation theory does not produce chains of whiskered KAM
tori with transverse heteroclinic intersections. The reason is that a
perturbation of size $\ep$ causes gaps in the set of whiskered KAM tori of
order $\ep^{1/2}$ near the (first order in $\ep$) resonance surfaces of order
$1$. On the other hand, the effect of the perturbation on the stable and
unstable manifolds of these whiskered tori is only of $O(\ep)$. Hence, a
naive perturbation theory cannot establish the  existence of chains of
transition tori traversing the resonance surface. The main goal of this paper
is to describe a geometric mechanism, based on a the computation of a
\emph{scattering map} defined on the normally hyperbolic invariant manifold
of the system, which is used to establish the existence of chains of
transition tori  that can traverse these resonance regions devoid of primary
KAM tori, that is, tori that can be continued from invariant tori of the
integrable system. The main idea of the mechanism proposed is to use
secondary tori, that is, tori generated by the resonances that can not be
continued from  invariant tori of the integrable system. We will show that
these secondary tori generated by the resonances fill the gaps created in the
set of KAM primary tori.

The result for $d =1$, $n = 1$ was already established in
\cite{DelshamsLS03a,DelshamsLS06a}. The problem was reexamined in
\cite{DelshamsH09}, where the hypothesis about $Q$ being a trigonometric
polynomial in $(\varphi,t)$ was eliminated. The works
\cite{GideaL06a,GideaL06b} supplemented the methods of
\cite{DelshamsLS03a,DelshamsLS06a} with the use of the method of
correctly-aligned windows. This method allowed to simplify the proof and to
obtain explicit estimates on the time spent by the diffusion trajectories. In
\cite{GideaL06b}, which considers $n$ arbitrary, a particularly careful
version of the method of correctly aligned windows produces essentially
optimal estimates on time and allows to weaken the non-degeneracy assumptions
in \cite{DelshamsLS03a,DelshamsLS06a}.

The case $d \ge 2$ presents a difficulty that was not present in the case $d
= 1$, namely, that there are points where the resonances have higher
multiplicity (the multiplicity of a resonance is the dimension of its
resonance module, see~\eqref{multiplicity}).

The technique used in \cite{DelshamsLS03a,DelshamsLS06a} was to take
advantage of the fact that, in the neighborhood of simple resonances, that
is, resonances of multiplicity $1$, it is possible to introduce a normal form
which is integrable, and  can be analyzed with great accuracy. Unfortunately,
it is well known that multiple resonances, that is,  resonances of
multiplicity greater or equal than $2$, lead to normal forms that  are not
integrable and need other techniques to be analyzed (see
\cite{ArnoldKN88,Haller99}). Recent progress in the analysis of double
resonances can be found in \cite{Marco13,KaloshinZ12,Cheng12}. We also note
that \cite{Treschev12} establishes diffusion far from \emph{strong
resonances} for the case $n=1$, $d\ge 2$, using the method of the separatrix
map.

In this paper we adapt the methods of the previous papers
\cite{DelshamsLS03a,DelshamsLS06a} to show instability under explicit
conditions. The basic observation is that multiple resonances happen in
subsets of codimension greater than $1$ in $\II$. We will adapt the methods
of \cite{DelshamsLS03a,DelshamsLS06a} to analyze the behavior of the system
in regions of simple resonances  and show that the diffusing trajectories can
contour the multiple resonances.

As we will see in Theorem~\ref{main}, we can choose largely arbitrary paths
in the action space $\II$ (we just require that they do not pass through some
higher codimension subsets of multiple resonances) and then, show that there
are orbits whose $I$-projection follows these paths up to an  error, which
becomes arbitrarily small with $\ep$. The sets that can be reached are of
size $O(1)$ and they include simple resonances.

Similar definitions of diffusion along paths were also used in
\cite{ChierchiaG94,ChierchiaG98}, but the methods of these and related papers
\cite{BertiB02,BertiBB03} only established the existence of diffusion in sets
completely devoid of resonances (the so-called ``gap bridging mechanism'').
The orbits that we produce cross the codimension $1$ resonant surfaces of
multiplicity one. Similar phenomena have been observed in the heuristic
literature \cite{Chirikov79,LiebermanT83}. In \cite{Tennyson82}, a similar
mechanism of diffusion is called \emph{diffusion across resonances}.

Of course, we are far from believing that the mechanism discussed in this
paper  is the only one to produce changes of order one in the actions. In
particular, combinations of variational and geometric methods have been
recently applied in Hamiltonians with $2+1/2$ degrees of freedom
\cite{Mather04,Cheng12,Zhang11,KaloshinZ12,Bernard10,Marco13,BessiCV01} which
require the Hamiltonian to be positive definite, which is a non-generically
(albeit open) property.

Also, the paper of \cite{Tennyson82} suggests several other mechanisms that
should be at play. It seems to be a  very challenging problem to make
rigorous the heuristic discussions on statistical and quantitative properties
of different instability mechanisms in the heuristic literature
\cite{Chirikov79,LiebermanT83,Tennyson82}. Of course, the heuristic
literature is convinced that double resonances help diffusion because they
are one of the ingredients in some of the heuristic mechanisms (but not in
others!). It is somewhat paradoxical that the rigorous mathematical theory
has difficulty precisely at the places which heuristics considers favorable.

In this paper, we show rigorously that double resonances can be contoured.
This can enhance the believe that there are several mechanisms.

Note however that if we consider \eqref{modelsconsidered} as a model of what
happens in a resonance (the pendulum being the resonant variable), then the
multiplicity of the resonances in the real system is one more than the order
that appears in the model. Hence, what we call simple resonances in our
Hamiltonian (modeling a resonance in a real model) would be double resonances
in the real Hamiltonian.

Since the proof presented here is quite modular and has well defined
milestones, we think that it is  almost certain  that other methods can be
applied to improve some of our arguments. In particular, we expect that the
method of correctly aligned windows can also give alternative proofs or to
improve several steps of the proof. The field of instability has experienced
a great deal of activity in recent years and there is a large variety of
results that have been obtained or announced. For a more detailed survey of
recent results, we refer to
\cite{DelshamsGLS08,PiftankinT07,Cheng08,Cheng10,Bernard10}.

\begin{remark}\label{twoparameter}
It is customary in some literature to refer to models of the form
\eqref{modelsconsidered} as a-priori-unstable models. We note, however, that
this distinction only  makes sense in considering analytic models depending
only on one small parameter. The results we  present here apply just as well
when the potentials $V_i$ in \eqref{modelsconsidered} are arbitrarily close
to $0$. In such a case, we just need to choose $\ep$ very small (even
exponentially small) relatively  to the hyperbolicity properties of the
$V_i$. Several papers in the literature, notably those dealing with generic
results, which occur in typical or cusp-residual Hamiltonian, call these
systems ``a priori stable". In particular, one can use this method to produce
systems that present instability but which are as close to integrable as
desired.  This procedure was pioneered in \cite{Arnold64}.
\end{remark}

\begin{remark} Hamiltonian  \eqref{modelsconsidered}   can
be considered as a simplified model of what happens in a neighborhood of a
resonance of multiplicity $n$ in a near integrable Hamiltonian. The averaging
method \cite{LochakM88, Haller97, Haller99, ArnoldKN88} shows that near a
resonance of multiplicity $j$, one can reduce a near integrable Hamiltonian
to a Hamiltonian of the form
\begin{equation}\label{eq:model}
h(I)+ \sum _{i=1}^{n}\frac{p_{i}^{2}}{2}+\ep V(q_{1}, \dots ,q_{n},I)
+ O(\ep ^{2}).
\end{equation}

The assumption that the averaged system  is given by uncoupled pendula is
made often \cite{HolmesM82,Haller97}. It is a {\emph {generic}} assumption
for $n =1 $. Hence we expect that the mechanism presented here is typical in
a neighborhood of a resonance. Of course, the hyperbolicity will be weak in
systems close to integrable, but in families with two parameters, it would
suffice to exclude wedges. See Remark~\ref{twoparameter}.

For $n \ge 2$, the above model \eqref{eq:model} is, in general, not
integrable whereas the pendulum part of \eqref{modelsconsidered} is .
Nevertheless, we point out that  the only think we need for our analysis is
that $\sum _{i=1}^{n}\frac{p_{i}^{2}}{2}+\ep V(q_{1}, \dots ,q_{n},I) $
admits transversal homoclinic orbits to a hyperbolic equilibrium point.
Systems of the form \eqref{modelsconsidered} appear naturally in several
physical models. A motivation to include this generality is that there is
very little difference dealing with any $n$ and it allows to emphasize that
the geometric methods allow to deal with systems that are not positive
definite.

\end{remark}

\section{Notation, assumptions and results}
\label{sec:notation}

In this section, we will present an overview of the argument and formulate
precisely most of the non-degeneracy assumptions we will assume. We will
postpone the precise formulation of the most technical ones till we have
developed the notation for them and motivated their explicit expressions.

The proof is divided in well defined steps and each of them can be
accomplished using  standard tools. We hope that the experts in these
techniques can fill in the arguments better than the authors, so that for
many possible readers, the heuristic discussion will be enough.

We have found it convenient to present the argument in an order slightly
different from the one followed in \cite{DelshamsLS03a, DelshamsLS06a} so
that, even if our hypotheses correspond closely to the assumptions of these
papers, the numbers do not correspond. In some cases, we have chosen to
present the result under slightly different assumptions than in
\cite{DelshamsLS03a, DelshamsLS06a} to simplify the exposition.  For the same
reason---simplifying and shortening the exposition---some of the objects
whose expansions were explicitly computed  in \cite{DelshamsLS06a} will now
be given through existence theorems, so that the final conditions will become
slightly less explicit. Nevertheless, since the procedures here are rather
constructive, explicit formulas can be given through more detailed work. On
the other hand, we note that the tool of the scattering map and its
symplectic properties \cite{DelshamsLS08},  a tool which was not available
when \cite{DelshamsLS06a} was written, simplifies significantly the
computations and thus, the conditions we obtain in this paper are simpler to
verify and more generally applicable than those in \cite{DelshamsLS06a}.

The precise statement of the main result (Theorem \ref{main}) requires the
definition of the resonance web, which depends on assumption {\bf H3}. The
statement of the last non-degeneracy conditions,  {\bf H6}, {\bf H7}, {\bf
H8}, can only be made after the system has been analyzed near resonances. We
note that these conditions are verifiable in concrete models with a finite
calculation, as it is performed in the example \eqref{eq:example}.

\subsection{Some  elementary notation: the extended flow, the
time-one map}

We will always  consider the extended flow $\tilde \Phi_{\ep,t}(\tilde x)$
which is obtained by supplementing the standard Hamilton equations with the
equation $\dot s  = 1$:
\begin{eqnarray}\label{eq:hamiltonianequations}
\dot I &=&-  \ep  \frac{\partial Q}{\partial \vp}(I,
\vp,p,q,s;\ep)\nonumber \\
\dot \vp &=& \frac{\partial h}{\partial
I}(I)+ \ep \frac{\partial
Q}{\partial I}(I, \vp,p,q,s;\ep)\nonumber \\
\dot p &=& -\frac{\partial P}{\partial q}(p,q)- \ep  \frac{\partial
Q}{\partial q}(I, \vp,p,q,s;\ep)\\
\dot q &=& \frac{\partial P}{\partial p}(p,q)+ \ep  \frac{\partial
Q}{\partial p}(I, \vp,p,q,s;\ep)\nonumber \\
\dot s &=&1 \nonumber
\end{eqnarray}

To the extended differential equations \eqref{eq:hamiltonianequations}
 corresponds the extended
phase space  $\tilde M := \II \times \torus^d \times \real^n \times \real^n
\times \torus $ associated to the variables $\tilde x=(I, \vp, p, q, s)$
respectively.

Some of our calculations are made easier by considering the time-$1$ map of
the flow. We will use the notation $f_\ep$ to denote the time one map
starting at the initial condition $t = 0$.

\subsection{The first elementary assumptions:
regularity,  hyperbolicity of the pendula and  non-degeneracy of the
integrable part}

We will be making the following assumptions:
\begin{itemize}
\item{{\bf H1}} We will assume that  the functions $h, V_j, Q$ are $C^r$
    in their corresponding domains with $r \ge r_0$ sufficiently large.
\end{itemize}
\begin{itemize}
\item{{\bf H2}} We will assume that the potentials $V_j$ have
    non-degenerate local maxima each of which gives rise, at least,  to a
    homoclinic orbit of the pendulum $P_{j}$.
\end{itemize}

Without loss of generality and to simplify the notation, we will assume that
the maxima of the potentials $V_j$ happen at $q_j = 0 $. That is, we will
assume that $V'_j (0) = 0$, $V''_j(0) = -\alpha_j^2$ with $\alpha _j \ge
\alpha > 0$, $j=1,\dots, n$.

We will denote by $(p^*_j(t), q^*_j(t))$ a parameterization by the natural
time of the homoclinic orbit we have chosen. That is,
\begin{equation}\label{eq:homoclinics}
\begin{split}
& \frac{d}{dt} q^*_j (t) = p^*_j(t); \\
&\frac{d}{dt} p^*_j (t) = -V'_j(q^*_j(t));  \\
&\lim_{t \to \pm \infty } (p^*_j (t),q^*_j (t))  = (0,0).
\end{split}
\end{equation}

When the variables $q_j$ have the physical interpretation of angles it is
natural to assume that the $V_j$ are periodic. In such a case, the limit in
(\ref{eq:homoclinics}) is understood modulus the period of the potential. The
method of proof only requires the existence of the homoclinic orbits to
hyperbolic saddles. Hence it applies to the coupling of an integrable system
in the $(I,\vp)$ variables to a chaotic system in the $(p,q)$ variables.

We note that the choice of a parameterization of the homoclinic orbit in the
full space involves the choice of $n$ origins of time in each of the
homoclinic orbits. Subsequent hypotheses will be independent of these
choices. The possibility of choosing the origin of the parameterizations of
each of the homoclinic orbits independently will play an important role in
our discussion of the Poincar\'{e} function in \eqref{Melnikovpot}.

Once we have chosen a homoclinic orbit to the origin in any pendulum, we
obtain a homoclinic connection in the space of the pendula. We will denote by
$\U \subset \real ^n \times \real ^n$ a neighborhood of the homoclinic
connection chosen in the $p,q$ space.

In this paper, we will assume that the equilibrium  points are  hyperbolic.
It would be interesting to extend the result to degenerate maxima (leading to
weakly hyperbolic points). This case has been proposed in the literature
\cite{DuToitMM09}.

{From} now on during the paper, we will consider $\II^* \subset \II$, and we
consider the compact set
\begin{equation}\label{difusionset}
\D:= \II^*\times \torus ^d  \times \U \times \torus ^1
\end{equation}
to be the domain of our problem. So, all the hypotheses refer to this domain.

\begin{itemize}
\item{ {\bf H3} } The mapping $I \rightarrow
    \omega(I):=\frac{\partial}{\partial I} h(I)$ is a diffeomorphism from
    $\II^*$ to its image.
\end{itemize}

\subsection{Assumption on the  structure of the perturbation}
We will furthermore assume:

\begin{itemize}
\item{ {\bf H4} } The function $Q$ in (\ref{modelsconsidered}) is   a
    trigonometric polynomial on $(\vp, t)$.
\end{itemize}
That is, we can write
\begin{equation}\label{expansionformula}
Q(I, \vp, p, q, t; \ep) = \sum_{(k, l) \in \N_Q} Q_{k,l}(I,p, q;
\ep) e^{2 \pi i (k \cdot \vp + l t)}
\end{equation}
with $\N_Q \subset \integer^d \times \integer$ a {\emph{finite set}}, with
$Q_{k,l}\not\equiv 0$ in $\II ^*\times \U$, if $(k,l)\in \N_Q$.

Hypothesis {\bf H4} clearly does not belong in the problem and we hope to
eliminate it in future treatments. Since the main goal of this paper is to
deal with the issue of multiple resonances, we have thought it convenient to
make the result as easy to read as possible, even if we do not achieve the
largest possible level of generality.

\begin{remark}\label{trigonometricpol}
Hypothesis {\bf H4}  appeared in \cite{DelshamsLS06a} for the case $d = 1$,
$n = 1$. In that case, the trigonometric polynomial hypothesis  has been
eliminated in \cite{DelshamsH09, GideaL06b} under generic assumptions. The
paper \cite{GideaL06b} eliminated the trigonometric polynomial hypothesis for
$d = 1$, $n$ arbitrary (one can argue that, possibly, the orbits produced are
not the same as the orbits in the previous papers). Similar improvements  are
clearly possible in the higher dimensional case $d>1$.
\end{remark}

\begin{remark}\label{rem:weakerhypotheses}
The methods we will use here can reach the same conclusions under slightly
weaker hypotheses.

We only need that, for some big enough but finite $m\le r$ the sets of
integer indexes
\begin{equation}\label{hypothesisgeneral}
\{ (k,l) \in \N _{Q}; \  \ \frac{ \partial^i}
{\partial \ep^{i}} Q_{k,l}(I, p = 0, q = 0 ; \ep = 0 ) \ne 0 \}
\end{equation}
for $i\le m$ are finite. This  happens, in fact, for some models in celestial
mechanics  \cite{FejozGKR11}. Nevertheless, since the writing in this case
becomes more cumbersome, we will only claim the weaker result.
\end{remark}

\subsection{Non degeneracy assumptions}
The first non-degeneracy assumption concerns the averaged Hamiltonian near
simple resonances, and are stated in Section \ref{nondegeneracy7}.

\begin{itemize}
\item {{\bf H5}} Consider the set of integer indexes $\N^{[\le 2]}= \N_1
    \cup \N_2 \subset \integer^{d+1}$ where $\N_1$ is the support of the
    Fourier series of the perturbation $Q(I,\varphi,p,q,t;0)$, $\N_2=
    (\N_1+\N_1)\cup \bar \N$, where $\bar \N$ is the support of the
    Fourier series of $\frac{\partial Q}{\partial
    \ep}(I,\varphi,p,q,t;0)$.

Then we assume that, for any $(k, l) \ne (0,0) \in \N ^{[\le 2]} $, the
set
\begin{equation}\label{degenerateresonances}
\{ I\in \II ^*, D h(I) k+l =0, \quad k^{\top} D^{2}h(I)k=0\}
\end{equation}
is empty or a manifold of codimension at least two in $\II^*$ .

In the case the map $\tilde h(I_0,I)= I_0 +h(I)$ is a quasi convex
function the set \eqref{degenerateresonances} is an empty  set for any
$(k, l) \ne (0,0) \in \integer^d \times \integer$ and a fortiori for any
$(k, l) \ne (0,0) \in \N ^{[\le 2]}$. Therefore Hypothesis {\bf H5} is
true for any perturbation $Q$ in this case.
\item{ {\bf H6}} Assume that the perturbation $Q$ satisfies some
    non-degene\-racy conditions stated in Section \ref{nondegeneracy7} in
    the connected domain $\II^* \times \torus ^{d+1}$ related to the
    averaged Hamiltonian.
\end{itemize}

The following non-degeneracy assumptions concern the so called Poincar\'e
function (or Melnikov potential) associated to the homoclinic connection
$(p^*, q^*)$ chosen in assumption {\bf H2}:
\begin{equation} \label{Melnikovpot}
\begin{split}
L(\tau, I, \vp, s)  & =- \int_{-\infty}^{\infty}  \big[ Q(I, \vp + \omega(I)
\sigma, p^*( \tau+\sigma ),q^*( \tau+ \sigma ), s+\sigma ;0 )  \\
&   -
Q(I, \vp + \omega(I)\sigma,0, 0, s+\sigma  ;0)   \big] \, d\sigma \\
\end{split}
\end{equation}
where
\begin{eqnarray*}
\tau  &=&  (\tau_1,\dots,\tau_n) \\
 p^*(\tau+ \sigma )&=& \left(p^*_1( \tau_1+ \sigma ),\ldots, p^*_n(\tau_n + \sigma )\right)\\
 q^*(\tau+ \sigma )&=& \left(q^*_1( \tau_1+ \sigma ),\ldots, q^*_n(\tau_n + \sigma )\right)
\end{eqnarray*}

\begin{itemize}
\item{ {\bf H7}} Assume that, for any value of $I\in \II^*$, there exists
    a non-empty set $\J _I \subset \torus ^{d+1}$, with the property that
    when $(I,\varphi,s)\in H_-$, where
\begin{equation}\label{domainscattering}
H_-=\bigcup _{I\in \I^*}\ \{I\}\times \J_I\subset \II^*\times \torus ^{d+1},
\end{equation}
the system of equations
\begin{equation}
\frac{\partial}{\partial \tau} L(\tau, I, \vp, s) = 0
\end{equation}
admits a non degenerate solution $\tau = \tau^*(I, \vp, s)$ with $\tau
^*$ a smooth function.

\item{ {\bf H8}} Define the auxiliary functions (related to the
    scattering map that will be introduced in Section
    \ref{sec:scattering0})
\begin{equation} \label{generator}
\L(I, \vp, s) =  L( \tau^*(I,\vp,s), I, \vp, s), \quad
\L^{*}(I,\theta)= \L(I, \theta, 0)
\end{equation}

Assume that the {\emph{reduced Poincar\'{e} function}} $\L^{*}(I,
\vp-\omega(I)s)$ satisfies some non-degeneracy conditions stated in
Section \ref{sec:scattering} in the domain $H_-$ (see
\eqref{eq:nondegeneratenores}, \eqref{eq:nondegenerateresonant}).
Nevertheless we anticipate that an informal description of the hypothesis
will be discussed after the main theorem \ref{main}.
\end{itemize}

We also note that the hypothesis {\bf H8} is simplified by the very simple
hypothesis:
\begin{itemize}
\item {\bf {H8'}}: $\forall I \in \II^*$, the reduced Poincar\'{e}
    function $\L^{*}(I,\theta)$ defined in \eqref{generator} has non
    degenerate critical points.
\end{itemize}

\subsection{Statement of the main result}
The main result of this paper is the following:
\begin{theorem}\label{main}
Let $H$ be a Hamiltonian of the form \eqref{modelsconsidered} satisfying  the
elementary assumptions {\bf H1, H2},  the regularity assumption {\bf H3}, the
simplifying assumption  {\bf H4} and the non-degeneracy assumptions {\bf H5,
H6, H7, H8}.

Then, for every $\delta>0$, there exists $\ep_0>0$, such that for every $0 <
|\ep| < \ep_0$, given  $I_{\pm}\in \II^*$, there exists an orbit $\tilde
x(t)$ of \eqref{modelsconsidered} and $T>0$, such that:
\begin{equation} \label{conclusion}
\begin{array}{lcr}
| I(\tilde x(0))  - I_- | \le C \delta \\
| I(\tilde x(T))  - I_+ | \le C \delta .
\end{array}
\end{equation}

\end{theorem}

Actually, we will show that given a largely arbitrary path $\gamma (s)
\subset \II^*$, we can find orbits $\tilde x(t)$ such that $I(\tilde x(t))$
is $\delta$-close to $\gamma (\Psi(t))$ for some reparameterization $\Psi$.
We postpone the precise statement till we have developed the notation. See
Theorem \ref{thm:path}.

The set $\II^*$ will be described precisely in the course of the proof. The
set is determined by the non-degeneracy assumptions {\bf H5}, {\bf H6}, {\bf
H7} and {\bf H8}. Given any concrete system, the assumptions can be verified
from the finite jet in $\varepsilon$ of $H$. Therefore, these conditions hold
in regions of order $1$ of $\I$.

The main restriction to obtain the set $\II^*$ from the domain of definition
$\II$ is given by assumption ${\bf H7}$, which guarantees the transversality
of the intersection of the  stable and unstable manifolds of the perturbed
normally hyperbolic manifold $\Lambda_\ep$ obtained in section
\ref{normalhyp}. Once ${\bf H7}$ is imposed, one needs to eliminate some sets
of codimension two from it to obtain $\II^*$:

\begin{itemize}
\item {\bf H8} Eliminates the values of $I$ for which the scattering map
    is not transversal to the inner map. More precisely, the invariant
    KAM tori of the inner map are transverse to their images under the
    scattering map. See Section~\ref{sec:scattering}.

In fact, the different conditions given in  Section~\ref{sec:scattering}
which constitute hypothesis {\bf H8} can be replaced by the sufficient
condition {\bf H8'}.

\item {\bf H5} and {\bf H6} Eliminate the region in the resonances where
    the leading term of the averaged system is degenerate. That is, given
    the codimension $1$ resonant surfaces, we have to eliminate the place
    when some function vanishes. See Section~\ref{nondegeneracy7}.
\end{itemize}

Since $\II ^*$ can contain multiple resonances appearing up to finite order
averaging  theory, and our mechanism is based on avoiding these multiple
resonances, we choose $\delta >0$ and contour these codimension two sets up
to a distance of order $O(\delta)$, obtaining a reduced domain $\II_\delta
\subset \II^*$. The precise definition of the sets to eliminate is deferred
to Section~\ref{averaging}.  Roughly, we eliminate the double resonances in
which one of the resonances are or order 1 or 2 (double resonances in which
both resonances are of order higher than $2$ are allowed in $\II_\delta$).

Note that, since  by ${\bf H4}$, the perturbation is a  trigonometric
polynomial and we assume hypothesis {\bf H3}, we only need to eliminate the
intersection of finite number of codimension one manifolds.

We note that all the conditions {\bf H5-H8} are generic: $C^2$ open in the
space of Hamiltonians and  hold except in sets of infinite codimension. Note
that all these hypothesis are transversality conditions among objects that
are independent. That is, we require transversality conditions among objects
that depend on the perturbation restricted to different places. See the
details later.

The only hypothesis that is not generic in the set-up is  assumption {\bf
H4}. It seems clear that this assumption {\bf H4}  can be eliminated using
the techniques developed in \cite{DelshamsH09}. However, we have preferred to
maintain it to simplify the exposition.  Roughly,  the idea is that, for
every $\ep >0 $ one can approximate the perturbation by a trigonometric
polynomial. If the trigonometric polynomial verifies the hypotheses of
Theorem~\ref{main}, one can obtain the existence of wandering paths and
information about their robustness. As it turns out, if the perturbation is
smooth enough, one can show that truncation error is much smaller than the
robustness allowed by the mechanism. Of course, there are quite a number of
details to be verified and we will not endeavor to do them now.

Once we have defined the set $\II_\delta \subset \II^*$, we will show that
KAM tori--either primary or secondary--are closely spaced on it. We will show
that,  given any KAM torus with $I$ coordinates in $\II_\delta$, it has
transversal heteroclinic connections with all the other KAM tori in a small
neighborhood. Applying the shadowing lemma, we can find orbits whose actions
follow almost  arbitrary paths inside $\II_\delta$. This is, of course,
slightly stronger than the conclusion \eqref{conclusion} of
Theorem~\ref{main}.

It is important to note that  codimension $2$ objects do not separate the
regions and can be contoured so that they do not obstruct the change along
the paths.

\section{Proof of Theorem~\ref{main}}

\subsection{First step: The use of normal hyperbolicity}\label{normalhyp}
We note that for $\ep = 0$, the manifold
$$
\Lambda_0 = \{p = 0, q = 0, \  I \in \II^*, \vp \in \torus^d \}
$$
is locally invariant under $f_0$, the time-$1$ map of the flow.

In this paper, we will not only work with the time one map but also  with the
flow $\tilde \Phi _{0,t}$ of system  \eqref{eq:hamiltonianequations}, so we
will work in the extended phase space $\tilde \M = \real ^{n}\times \real
^{n}\times \II \times\torus ^{d} \times \torus $.

In the extended phase space, we consider the invariant manifold
$$
\tLambda_0 = \{ p = 0, q = 0,\  I \in \II^*, \vp \in \torus^d,  s \in
\torus^1\}
$$
which is a normally hyperbolic invariant manifold under the flow $\tilde \Phi
_{0,t}$ in the sense of \cite{Fenichel71, Fenichel74, HirschPS77, Pesin04}.
That is, for every $\tilde x \in \tLambda_0$, there is a decomposition
\begin{equation}\label{decomposition}
T_{\tilde x}  \tilde \M = E^s_{\tilde x} \oplus E^u_{\tilde x} \oplus T_{\tilde x} \tLambda_0
\end{equation}
and numbers $C > 0$, $0 < \beta < \alpha $,  such that the decomposition
\eqref{decomposition} is characterized by:
\begin{equation}\label{characterization}
\begin{split}
& v \in E^s_{\tilde x} \iff |D\tilde \Phi_{0,t}(\tilde x) v | \le C e^{- \alpha t}  |v| \quad \forall t \ge 0\\
& v \in E^u_{\tilde x} \iff |D\tilde \Phi_{0,t}(\tilde x) v | \le C e^{- \alpha |t|}  |v| \quad \forall t \le 0\\
& v \in T_{\tilde x} \tLambda_0 \iff |D\tilde \Phi_{0,t}(\tilde x) v | \le C e^{ \beta |t|}  |v| \quad \forall t \in \real \\
\end{split}
\end{equation}

It is clear that the stable and unstable  spaces $E^{s,u}_{\tilde x}$  are
the direct sum of the stable and unstable spaces at the critical point of
each of the pendula $P_{j}$,  and $\alpha$ is given in assumption {\bf H2}.
Furthermore, we can take for $\beta $ any number satisfying
$0<\beta<<\alpha$.

The standard theory of  persistence of \NHIM \cite{Fenichel71, Fenichel74,
HirschPS77, Pesin04} implies that, for $|\ep| < \ep_0$ there is a locally
invariant normally hyperbolic manifold $\tilde \Lambda_\ep$ verifying
\eqref{decomposition} and \eqref{characterization} for the perturbed flow
$\tilde \Phi_{\ep, t}$ (with $\alpha_{\ep}, \beta_{\ep}, C_{\ep}$ close to
$\alpha,\beta,C$ respectively). The theory in \cite{Fenichel71, Fenichel74,
HirschPS77, Pesin04} guarantees that $\tilde \Lambda_\ep$ is a somewhat
smooth family of manifolds but the degree of smoothness can be limited by
ratios of normal and tangential exponents $\alpha$ and $\beta$. In our case,
since the motion on the manifold for $\ep = 0$ is integrable and therefore
the Liapunov exponents are zero, we can guarantee that for $|\ep|$ small
enough, the family $\tilde \Lambda_\ep$ will be a $C^{r-1}$ family if $\tilde
\Phi_{\ep,t}$ is a $C^r$ family.

Moreover, as it was shown in  \cite[Theorem 24]{DelshamsLS08} there is a
naturally defined symplectic parametrization $k_{\ep}$, such that the
perturbed manifold can be written as $\tilde \Lambda_{\ep}= k_{\ep}(\tilde
\Lambda_0)$, using as the reference manifold the unperturbed manifold $\tilde
\Lambda_{0}$, and choosing $k_{0}=\mathrm{Id}$.

Using this symplectic parameterization, one can show that the reduced flow
$\tilde \phi _{\ep, t}$ on $\tilde \Lambda_0$, characterized by $k_{\ep}\circ
\tilde \phi _{\ep,t} = \tilde \Phi _{\ep,t}\circ k_{\ep}$ is a Hamiltonian
flow. The following proposition makes explicit its Hamiltonian.
\begin{prop}\label{prop:restrictedflow}
The reduced  flow $\tilde \phi _{\ep,t} $ on $\tilde \Lambda_{0}$ defined
through $k_{\ep}\circ \tilde \phi _{\ep,t} = \tilde \Phi _{\ep,t}\circ
k_{\ep}$ is generated by a $\C^{r-1}$ time dependent Hamiltonian vector field
with Hamiltonian of the form
\begin{equation}\label{eq:hamiltonianrestricted}
K_{\ep}(I,\vp, s) = h(I) + \sum_{i = 1}^N \ep^i K^i(I, \vp, s) +
O_{C^{r -N -2}}(\ep^{N+1}),
\end{equation}
where each of the terms $K^i$ is a trigonometric polynomial in the $\vp, s$
variables.

Moreover, $K^i$ is an algebraic expression in terms of $\nabla^\ell Q(I, \vp,
p=0, q=0,s;\ep=0)$, for $\ell=0,\dots,i-1$. In particular, $K^1(I, \vp, s) =
Q(I, \vp, 0, 0, s;0)$.
\end{prop}

\subsection{Analyzing the flow restricted to the invariant manifold}
\label{sec:analyzing}
\def\tLambda{{\tilde \Lambda}}

The goal of next sections is to study the objects in $\tilde \Lambda_{\ep}$
invariant by the flow $\tilde \Phi_{\ep,t}$. Using Proposition
\ref{prop:restrictedflow}, this is equivalent to studying the objects in
$\tLambda_0$ invariant under the Hamiltonian flow $\tilde \phi_{\ep,t}$ of
Hamiltonian $K_\ep (I, \vp,s)$ given in \eqref{eq:hamiltonianrestricted}.

The main tool used to obtain invariant objects in $\tilde \Lambda_{\ep}$ will
be averaging theory \cite{DelshamsLS06a} and KAM theorem applied to  the
Hamiltonian  $K_{\ep}(I,\vp,s)$. We see that, after we add some extra
variable $I_0$ conjugated to the  variable $s\in\torus$, to make it
symplectic and autonomous, we are lead to considering a Hamiltonian of the
form:
\begin{equation}\label{tildeRep}
\tilde K_\ep = I_{0} + K_{\ep}= \tilde h(\tilde I) +
 \sum_{i=1}^N \ep^i K^i(I,\tilde \vp) + O_{C^{r-2-N}}(\ep^{N+1})
\end{equation}
were we have introduced the notation $\tilde \varphi =(\varphi,s)$, $\tilde
I=(I,I_0)$, and $\tilde h(\tilde I) = I_0 + h(I)$. We recall that Proposition
\ref{prop:restrictedflow}  tells us that, by choosing $|\ep|$ small enough
and assuming  the regularity $r$ of the Hamiltonian $H (p,q,I,\vp,s;\ep)$  in
\eqref{modelsconsidered} large enough, we can take $N$ as large as we want
and the regularity of the remainder in \eqref{tildeRep} is as large as we
want.

Furthermore, in our case, using the assumption {\bf H4} it follows that  the
$K^i$ are trigonometric polynomials in the angle variables $\tilde \vp$ with
Fourier coefficients that depend on $I$, but not on $I_0$:

\begin{equation}\label{trigonometric}
 K^{i}(I,\tilde \vp)= \sum _{(k,l) \in \N_i }K^{i}_{k,l} (I) e^{2\pi i (k\vp+ls)}
\end{equation}
$\N_{i}$ being  finite sets. Very explicit formulas for the coefficients
$K_{k,l}^{i}(I)$ are given in \cite{DelshamsLS06a}.

The fact that the perturbation terms do not depend  on the $I_0$ variable is
a reflection of the fact that $I_0$ is just a variable introduced to keep the
time rotating at unit speed. This is independent of the perturbations.

\subsection{The averaging method}\label{averaging}
In this section we recall the basis of the averaging method for time periodic
perturbations. The averaging method is a rather standard tool in Hamiltonian
dynamical systems and has an extensive literature. Modern surveys are
\cite{LochakM88,ArnoldKN88}.

The basic idea of the averaging method is to make symplectic changes of
variables carefully chosen so that the resulting Hamiltonian presents a
particularly simple form. There are many averaging theories depending on what
is the simple form to be achieved and what is the method used to keep track
of the simplifying transformations. In this paper, we will follow
\cite{DelshamsLS06a} and use the method of Lie transforms. The averaged
Hamiltonians we will use here are different from those used in
\cite{DelshamsLS06a} to accommodate the fact that resonances  for systems
with two or more degrees of freedom are manifolds in the action space whereas
for one degree of freedom, resonances are just points. In this paper, we also
consider more general unperturbed Hamiltonians $h(I)$---in
\cite{DelshamsLS06a}, the unperturbed Hamiltonian was just quadratic---but,
under hypothesis  \textbf{H3}, this  makes little difference.

\subsubsection{Some generalities on the averaging method}

We will follow the method of Lie transforms \cite{Cary81, Meyer91},
considering transformations obtained as the time-$1$ map of a Hamiltonian
$G(I,\tilde \vp)$ in $\real ^{d}\times \torus ^{d+1}$.

Given a Hamiltonian function  $G(I,\tilde \varphi)$ in the extended phase
space, we denote  by $\exp( G)$ the time-$1$ map of the Hamiltonian flow
generated by $G$.

The main technical result we will use about the time-$1$ map is a direct
consequence of Taylor's expansions and the regularity of the solutions of an
ordinary differential equation as well as the expression of the derivatives
of Hamiltonian functions in terms of Poisson brackets \cite{Thirring97}. So
that the following Lemma \ref{differentiability} is just a Taylor expansion
along trajectories.

\begin{lemma}\label{differentiability}
Let $A \subset B \subset \real^{d+1} $ be compact sets and $\ell \ge 2$.
There exists a constant $C=C(l,k)$, such that, given $G \in C^\ell( B \times
\torus^{d+1})$ satisfying
\[
|| G ||_{C^1( B \times \torus^{d+1} )}  <   d( A, \real^{d+1}
\setminus B)
\]
so that the Hamiltonian flow associated to $G$ starting in $A \times
\torus^{d+1}$ stays in the interior  of $B \times \torus^{d+1}$, we have:
\begin{itemize}
\item[a)] $\exp{G} \in C^{\ell -1}(A \times\torus^{d+1} ) $
\item[b)] $||\exp( G) - \Id||_{C^{\ell -1}(A \times\torus^{d+1} ) } \le C
    || G||_{C^{\ell}(B \times\torus^{d+1} )}$
\item[c)] given $H \in C^\ell( B \times \torus^{d+1})$ then:
\begin{multline*}
|| H\circ \exp( G) - H - \{H, G\}||_{ {C^{\ell - 2}(A
\times\torus^{d+1})}  }\\
\le C || H ||_{{C^\ell}(B \times \torus^{d+1} )} || G||_{{C^\ell}(B
\times \torus^{d+1} )}^2
\end{multline*}
where $\{\cdot, \cdot\}$ denotes the Poisson bracket in $\real
^{d+1}\times \torus^{d+1}$.
\item[d)] More generally, if $\ell > k+1$, there is an asymptotic
    expansion

\[
\begin{split}
|| H\circ \exp( G) - H& - \{H, G\} - \frac{1}{2} \{ \{H, G\},G \} \\
& - \ldots - \frac{1}{k!} \{\{\dots \{\{H, G\},G\},G\}\dots \}
||_{C^{\ell - k -1 }(A\times\torus^{d+1} )  }\\
&\le C || H || _{C^\ell(B \times\torus^{d+1})} || (G||_{C^\ell(B
\times\torus^{d+1} )})^{k+1}
\end{split}
\]
\end{itemize}
\end{lemma}

As a consequence  of Lemma~\ref{differentiability}, we obtain the following
algorithm, which is the main formal  step of the general averaging method and
which allows computations to high order.

\begin{algorithm}\label{recursive}
Given a sufficiently smooth  Hamiltonian averaged up to order $N-1\ge 0$

\begin{equation}\label{eq:restrictedN}
\tilde K_\ep ^{N-1}=  \tilde h + \sum_{i = 1}^{N-1} \ep^i {\bar K}^i
+ \ep^{N} K^{N} +  O(\ep^{N+1})
\end{equation}

Assume that we can find sufficiently smooth  ${\bar K}^{N}$, $G ^N$ solving:
\begin{equation}\label{averagingstep}
{\bar K}^{N} = K^{N} + \{\tilde h, G^N \}
\end{equation}

Then
\[
\begin{split}
\tilde K_\ep ^{N}&=K_\ep ^{N-1}\circ \exp( \ep^{N} G^{N}) =  \tilde h
+ \sum_{i = 1}^{N-1} \ep^i {\bar K}^i +
\ep^{N} {\bar K}^{N} +  O(\ep^{N+1}) \\
& =  \tilde h + \sum_{i = 1}^N \ep^i {\bar K}^i   + \ep^{N+1} {K}^{N+1}+ O(\ep^{N+2})
\end{split}
\]
\end{algorithm}

In general there are many choices for ${\bar K^{N}}$ and $G^N$.
In the following subsections, we will specify the choices that we make for
our case and establish estimates for the transformation $G^{N}$ and the new
Hamiltonian  $ \tilde K_{\ep}^{N}$. In particular, we will have estimates for
the averaged Hamiltonian $\bar K^{N}$.

\subsubsection{One step of averaging: The infinitesimal equations.
Resonances in one averaging step}

Our Hamiltonian $\tilde K_\ep$ given in  \eqref{tildeRep} is of the form
$$
\tilde K_\ep = I_0 + K_\ep
$$
where $K_\ep$ is given in \eqref{eq:hamiltonianrestricted} and only depends
on $(I,\tilde \varphi)$. Then, in our case, we will take $\tilde K_\ep ^N$ in
\eqref{eq:restrictedN}, with $\tilde h = I_0+h(I)$ and $\tilde K_\ep
^N-\tilde h $ only depends on $(I,\tilde \varphi)$. So that the function $G$
will depend only on $(I,\tilde \varphi)$.

The fact that $\tilde h $ in \eqref{eq:restrictedN} is given by $\tilde h=
I_0 + h(I)$, will allow us  to treat the averaging equation
\eqref{averagingstep} using Fourier series. We  find it convenient to divide
the phase space into different regions and perform different averaging
procedures in each region.

At every step of the iteration algorithm \ref{recursive}, given $K(I,\tilde
\varphi)$, we have to solve equation \eqref{averagingstep} for the unknowns
$\bar K(I,\tilde \varphi), G(I,\tilde \varphi)$:
\begin{equation}
\label{eq:generalaveragingstep} \bar K = K + \{\tilde h,G\}
\end{equation}
with $\tilde h (\tilde I)= I_0+h(I)$. Writing
\[
K=\sum_{\tilde k  \in \N \subset \integer^{d+1}}  K_{\tilde k} (I)
e^{2\pi i \tilde k \cdot \tilde \vp} =\sum_{(k, l)  \in \N \subset
\integer^{d+1}}  K_{k,l} (I) e^{2\pi i( k \cdot \vp+ls)},
\]
it is clear, because $\tilde h$ depends only on $\tilde I$, that the Poisson
bracket with $\tilde h$ is diagonal in Fourier series. Hence, it is natural
to search for $\bar K, G$ such that their Fourier series are supported  in
$\N$, the support  of the Fourier series of $K$. Hence we look for:
\[
\begin{split}
&G = \sum_{(k,l)  \in \N \subset \integer^{d+1}} G_{k,l} (I)
e^{2\pi i ( k \cdot \vp+ ls)}\\
& \bar K  = \sum_{(k,l)  \in \N \subset \integer^{d+1}} {\bar
K}_{k,l} (I) e ^{2\pi i ( k \cdot \vp+ ls)}.
\end{split}
\]

Using that
\[
\{\tilde h,G\} = - 2\pi i\sum_{(k,l)  \in \N \subset \integer^{d+1}}
(\omega(I) \cdot k +l)\, G_{k,l} (I) e^{2\pi i  (k \cdot \vp+l s)}
\]
where $\omega(I)=\frac{ \partial h}{\partial I}(I)= \nabla h(I)$, equation
\eqref{eq:generalaveragingstep} becomes a set of equations for the Fourier
coefficients:
\begin{equation}\label{averagingstepFourier}
 {\bar K}_{k,l}(I) - K_{k,l}(I)  = - 2\pi
i(\omega(I) \cdot k +l) G_{k,l}(I)
\end{equation}
The solution of \eqref{averagingstepFourier} is obtained by choosing $\bar
K_{k,l}$ and then, setting
\begin{equation}\label{Fouriersolution}
G_{k,l} (I) = \frac{K_{k,l}(I) - {\bar K}_{k,l}(I)}{2\pi i(
\omega(I) \cdot k +l)}.
\end{equation}

It is clear that the solution \eqref{Fouriersolution} requires special
treatment when
\begin{equation}\label{eq:res}
\omega(I) \cdot k+l=0, \ (k,l)\in \N.
\end{equation}

This motivates next definition.
\begin{de} \label{resonantset}
Given a Hamiltonian $h(I)$ we define a resonance as the set
\begin{equation}\label{defresonanceset}
\R_{k,l}=\{ I\in\II, \ \omega(I)  \cdot k+l=0 \}=
\omega^{-1} \left(\{ \Omega \in \real ^{d}, \ \Omega \cdot k +l=0 \}\right).
\end{equation}
where $\omega(I)= \frac{\partial h}{\partial I}(I)=\nabla h(I)$.
\end{de}

Let us observe that, by Hypothesis {\textbf{H3}}, for $I\in \II^*$, the
resonances are smooth  surfaces for $k\ne 0$,  as smooth as the map
$\omega=\frac{\partial h}{\partial I}$.

Moreover, $\R_{m\, k,m\, l}=\R_{k,l}$ for any $m \in \integer$, therefore any
two of these sets $\R_{k,l}$ and $\R_{\bar k,\bar l}$ either:
\begin{itemize}
\item Are identical, if and only if $(k,l) = (m \bar k, m\bar l )$ for
    some $m \in \integer$
\item Do not intersect
\item Intersect transversally in a manifold of codimension two without
    boundary
\end{itemize}

More generally, if we consider different resonances, the intersection
$$
\R_{k_1,l_1}\cap \dots \cap \R_{k_j,l_j} =
\omega^{-1} \left(\{ \Omega \in \real ^{d}, \ k_i \cdot\Omega+l_i=0 , i=1,\dots , j\}\right)
$$
will be a manifold of dimension $m$, where $m$ is the dimension of the
$\integer$-module $\M [(k_1,l_1),\dots ,(k_j,l_j)]$ generated by
$(k_1,l_1),\dots ,(k_j,l_j)$.

Given $\N \subset \integer ^{d+1}$, the support of the Fourier series of $K$,
we introduce the notation
\begin{equation}\label{tilden}
\tilde \N = \{ (k,l)\in \N \ | \ \nexists \
m \ne 1\in \integer , (\tilde k,\tilde l)\in \N, (k,l) = m (\tilde k,\tilde l)\}
\end{equation}
Notice that, with this notation  if $ \R_{k_1,l_1} = \R_{k_2,l _2}, \quad
\mathrm{with } \quad (k_i,l_i) \in\tilde \N, \quad i=1,2 $ one has that
$(k_1,l_1) = (k_2,l_2)$ and then
$$
\cup _{(k,l)\in \N }\R_{k,l} = \cup _{(k,l)\in \tilde \N }\R_{k,l}.
$$

Since resonances are sets of codimension $1$, it is natural to give in them a
system of $d-1$ coordinates. In the next lemma we will consider the function
$$
\Gamma_{k,l} :\II ^* \to \R_{k,l}
$$
which is a projection onto the resonance $\R_{k,l}$ along a transversal
bundle to it. The standard choice is the {\emph {orthogonal}} projection, in
such a way that the projection $\Gamma _{k,l} (I)$ is the closest point to
$I$ in $\R_{k,l}$. The orthogonal projection is well defined in a tubular
neighborhood with respect to the normal bundle of $\R_{k,l}$.

There is another simpler choice, which appears naturally when one realizes
that the dynamics close to resonances is generically ``pendulum-like" (see
equation~\eqref{averagedJ}). It is just  projecting along the bundle
$\R_{k,l}+\langle k\rangle$, that is, defining $\Gamma_{k,l}(I)$ as the
intersection of the straight line $\{I+tk, t\in\real\}$ with the resonance
$\R_{k,l}$. This projection will  be not well  defined close to points $I\in
\R_{k,l} \subset \II^*$ such that the direction given by the vector $k$ is
tangent to $\R_{k,l}$, that is, close to points $I$ satisfying:
\begin{equation}\label{projectionnodefined}
\begin{array}{rcl}
Dh(I) \cdot k + l&=&0\\
k^{\top} D^2h(I)k& = & 0
\end{array}
\end{equation}
Nevertheless, under hypothesis \textbf{H5}, for {\emph{secular resonances}},
that is, for $\R_{k.l}$ such that $(k,l) \in \N^{[\le 2]}$, these points are
a codimension two set (as one would expect from naive parameter counting
because \eqref{projectionnodefined} are two conditions). In the averaging
procedure used in this paper, we just need to deal carefully with secular
resonances to obtain a ``pendulum like" Hamiltonian near them. So, our
strategy will be to use this projection near the secular resonances and the
standard orthogonal one near the rest of resonances we encounter in the
higher averaging steps.

We define $\JJ _{\delta} \subset \II ^*$ as  the set of points  of $\II^*$
taking out a neighborhood of size $\delta$ of the points verifying
\eqref{projectionnodefined} for some $(k,l) \in \N^{[\le 2]}$. In the
resulting set $\JJ _{\delta} $ we will have a well defined projection in the
$k$ direction  for any secular resonance:
\begin{equation}\label{jotadelta}
\mbox{If} \quad I^* \in \JJ _{\delta} \cap  \R_{k_0,l_0}
\quad \mbox{then} \quad k_0^{\top} D^2 h(I^*)k_0 \ne 0 \ \mbox{for} \ (k_0,l_0)\in \N^{[\le 2]}
\end{equation}
Of course, if our Hamiltonian $h$ is quasiconvex, since the set satisfying
\eqref{projectionnodefined} is empty, we obtain that $\JJ _{\delta}=\II^*$.

For future reference, let us write the characterization of the projection
$\Gamma_{k,l}(I)$ along the $k$-direction:

\begin{equation}\label{projection}
I^*=\Gamma_{k,l}(I) \Longleftrightarrow I-I^*\in\langle k\rangle \mbox{ and }
\frac{\partial h}{\partial I} (I^*)\cdot k+l=0.
\end{equation}

Let us observe that for $I\in \JJ_{\delta}$ in a neighborhood of a secular
resonance  $\R_{k,l}$, $(k,l) \in \N^{[\le 2]}$, there exists a constant
$C\ge 1$ such that
\begin{equation}\label{eq:goodprojection}
\dist (I,\R_{k,l}) \leq \dist (I,\Gamma_{k,l}(I))\leq C\dist (I,\R_{k,l}).
\end{equation}
Indeed, $C$ can be chosen as any constant satisfying
$$
C> \frac{\norm{k}\norm{D^2h(I)k}}{\norm{k^{\top}D^2h(I)k}}, \quad \forall I
\in \JJ_{\delta}, \ (l,k)\in \N^{[\le 2]}
$$

In the case of the orthogonal projection the value of the constant $C$
\eqref{eq:goodprojection} is $1$. We emphasize that, in the averaging
procedure, the resonant sets are determined by the integrable Hamiltonian
$h(I)$ and not by the perturbation. Of course, given a concrete system, many
resonances do not play any role and only the resonances excited by the
perturbation play a role, as we will see in next lemma.

As we indicate before, the solution of the homological equation involves
choices of which terms are eliminated and which terms are kept in the
averaged Hamiltonian. The following lemma indicates the choice we will
follow. We will use in it a  general projection $\Gamma _{k,l}$, just
assuming that it verifies \eqref{eq:goodprojection}.

\begin{lemma}\label{lem:homological}
Let
\[
K(I,\vp ,s)=\sum _{(k,l)\in \N}K_{k,l}(I)e^{2\pi i(k\vp+ls)}
\]
be a Hamiltonian, with $\N = \N (K)\subset \integer ^{d+1}$ a finite set.
Assume that  $K$ is of class $\C ^{\ell}$ with respect to $I \in \JJ
_{\delta}\subset \II ^*\subset \real ^{d}$ and consider the resonant set
$\R_{\N}= \{I\in \JJ _{\delta}, \ \omega(I)\cdot k+l=0, \ (k,l) \in \N\}=
\cup _{(k,l) \in \N} \R_{k,l} \subset \JJ _{\delta}$.

Choose  $0<L<1$ small enough such that for any $(k,l), (\bar k, \bar l) \in
\N$, either $\R_{k,l} = \R_{\bar k, \bar l}$ or the tubular neighborhood of
$\R_{k,l}$ of radius $L$ does not contain $\R_{\bar k, \bar l}$.

Assume that we have a projection $\Gamma _{k,l} :\JJ _{\delta} \to \R_{k,l}$
such that it verifies \eqref{eq:goodprojection}.

Then, there exist $G(I,\vp,s)$ of class $\C^{\ell-1}$ with respect to $I$,
and $\bar K$ of class $\C ^{\ell}$ with respect to $I$, with $\N(G)$ and $ \N
(\bar K)$  finite sets. Moreover  $\N(G), \  \N (\bar K)\subset \N(K)$.

The functions $G(I,\vp,s)$ and $\bar K$ solve the homological
equation~\eqref{eq:generalaveragingstep} in $\JJ _{\delta}$ and satisfy:
\begin{itemize}
\item [a)] If $d(I,\R _{\N})\ge 2L$, then
\[
\bar K(I,\vp,s)= K_{0,0}(I).
\]
\item [b)] If $d(I, \R_{k_{i},l_{i}}) \le L$ for $i=1, \dots j$, then
\begin{eqnarray*}
\bar K(I,\vp,s) &=&
K_{0,0}(I) +\sum_{i=1}^{j} \left( \sum_{\nu=-N_i}^{N_i} K_{\nu k_{i}, \nu l_{i}}
(\Gamma_{k_{i},l_{i}}(I)) e^{2\pi i \nu (k_{i} \vp +l_{i}s)} \right) \\
&=& K_{0,0}(I)+U_{k_{1},l_{1}, \dots ,k_{j},l_{j}}(I,k_{1} \vp +l_{1}s, \dots , k_{j} \vp +l_{j}s),
\end{eqnarray*}
where $0<N_i<\infty $ are such that if $(\nu k_{i},\nu l_{i})\in \N$,
then $\abs{\nu}\le N_i $.
\item [c)] The function $\bar K$ verifies: ${\| \bar K\|}_{\C ^{\ell}}
    \le (1+\frac{C}{L^{\ell+1}}){\|K\|}_{\C ^{\ell}}$, where $C$ is a
    constant independent of $L$.
\item [d)] The function $G$ verifies ${\|G\|}_{\C ^{\ell-1}} \le
    \dfrac{C}{L^{\ell +1 }} {\|K\|}_{\C ^{\ell}}$.
\end{itemize}
\end{lemma}

\begin{remark}\label{redundantangles}
We observe that in the formula of $\bar K$ the angles can be redundant
because some of the angles included in the sum can be combination of others.
One can be more precise by  considering  the module generated by
$(k_{1},l_{1}), \dots ,(k_{j},l_{j})$ and the dimension of this module  gives
us the number of independent angles among $k_{1} \vp +l_{1}s, \dots , k _{j}
\vp +l_{j}s$. However, in this paper, this will not be needed. Our strategy
later will be to reduce the domain in such a way that we will not need to
deal with multiple resonances.
\end{remark}

\begin{proof}
If we write the homological equation~\eqref{eq:generalaveragingstep} in
Fourier coefficients, we obtain equation \eqref{Fouriersolution}. Our first
choice is  $\bar K_{k,l}(I) =G_{k,l}(I)=0$, if $(k,l) \notin \N$.
For $(k,l) \in \N$, we solve  equation~\eqref{Fouriersolution} choosing:
\begin{enumerate}
\item If $(0,0)\in\N$ we take $\bar K_{0,0}(I) =K_{0,0}(I)$.
\item If $(0,l) \in \N$, $l\ne 0$, $\bar K_{0,l}(I) =0$.
\item If $(k,l)\in \N$, $k\ne 0$, we choose $\bar K_{k,l}(I)$ as:
\[
\bar K_{k,l}(I)= K_{k,l}(\Gamma_{k,l}(I))\psi\left(\frac{1}{L}(d(I,\R_{k,l}))\right)
\]
where $\psi(t)$ is a fixed $\C^{\infty}$ function such that: $\psi(t)=1$,
 if $t\in [-1,1]$, and $\psi(t)=0$, if $t\notin
[-2,2]$. With this choice we have that $\bar K_{k,l}$ verifies:
\begin{enumerate}
\item If $d(I,\R_{k,l})\le L$ then  $\bar K_{k,l}(I) =
    K_{k,l}(\Gamma_{k,l}(I))$,
 \item if $d(I,\R_{k,l})\ge 2L$ then  $\bar K_{k,l}(I) = 0$.
\end{enumerate}
\end{enumerate}
Once we have defined $\bar K=\sum _{(k,l)\in \N}\bar
K_{k,l}e^{i(k\varphi+ls)}$, it is clear that it has the form announced in a)
and b), is a $\C^{\ell}$ function with respect to $I$, and that it verifies
the bounds c), where the constant $C$ only depends on the cut-off
$\C^{\infty}$ function $\psi$, the  functions $\Gamma_{k,l}$, the degree of
the Fourier polynomials and $\ell$.

Now, we choose $G$  that satisfies equation~\eqref{averagingstepFourier}:
\begin{enumerate}
\item $G_{0,0}(I)=0$,
\item For $(0,l)\in \N$, $l\ne 0$, $G_{0,l}(I)=\dfrac{K_{0,l}(I)}{2\pi i
    l}$,
\item If $ (k,l)\in \N$, $k\ne 0$,  we choose $G_{k,l}(I)$ as:
\begin{enumerate}
 \item If $\omega(I)\cdot k +l \ne 0$ then $G_{k,l}(I)=\dfrac{
     K_{k,l}(I)-\bar K_{k,l}(I)}{2\pi i(\omega(I)\cdot k +l)}$,
 \item If $\omega(I)\cdot k +l = 0$ and we are using the standard
     orthogonal projection,  then $G_{k,l}(I)=\dfrac{\nabla
     K_{k,l}(\Gamma _{k,l}(I))\cdot D^2 h(I) k} {2\pi
 i||D^2 h(I) k||^{2}}$.
\item If $\omega(I)\cdot k +l = 0$ and we are using the
    $k$-projection, then $G_{k,l}(I)=\dfrac{DK_{k,l}(I)k}    {2\pi
    i\norm{k^{\top}D^2h(I)k}}$.
\end{enumerate}
\end{enumerate}

To bound  the function $G$ we first bound its Fourier coefficients
$G_{k,l}(I)$:
\begin{enumerate}
\item For $(0,l)\in\N$, $l\ne 0$, ${\|G_{0,l}\|}_{\C ^{\ell -1}} \le
    C{\|K_{0,l}\|}_{\C ^{\ell -1}}$.
\item Given  $(k_{0},l_{0})\in \N$, $k_0 \ne 0$, by the definition of
    $\bar K$ and $G$, we have:
\begin{enumerate}
\item On $\{ I \in \II, \ d(I,\R_{k_{0},l_{0}})\le L\}$, we have
    ${\|G_{k_{0},l_{0}}\|}_{\C ^{\ell -1}} \le
    C\dfrac{{\|K_{k_{0},l_{0}}\|}_{\C ^{\ell }}}{\abs{k_{0}}}.$

\item On $\{ I\in \II, \ d(I,\R_{k_0,l_0})\ge 2L\}$, we have
    ${\|G_{k_{0},l_{0}}\|}_{\C ^{\ell -1}} \le
    C\dfrac{{\|K_{k_{0},l_{0}}\|}_{\C
    ^{\ell}}}{(\abs{k_{0}}L)^{\ell}}.$

\item On $\{I\in \II, \  L\le d(I,\R_{k_0,l_0})\le 2L\}$, we have
 \[
{\|G_{k_{0},l_{0}}\|}_{\C ^{\ell -1}}
 \le C\frac{{\|K_{k_{0},l_{0}}\|}_{\C ^{\ell -1}}}{(\abs{k_{0}}L)^{\ell-1}}.
 \]
 \end{enumerate}
\end{enumerate}
Therefore, $G(I,\vp,s)$ is a trigonometric polynomial in $(\vp, s)$, and of
class $\C^{\ell-1}$ with respect to $I$ and satisfies the bounds in d).
\end{proof}

\begin{remark}
We note that the above estimates use the fact that the function $K$ is a
trigonometric polynomial in $(\vp, s)$ and the constant $C$ in the bounds of
$G$ depend on the degree of the polynomial.

One can follow the same procedure when the function is not a polynomial by
estimating the Fourier coefficients using Cauchy bounds:
$$
\| K_{k,l}\| \le C \|K\|_{\C ^r} (|k|+|l|)^{-r}
$$
and then using the estimates we presented. The only difference  is that for
general functions we obtain estimates:
\[
{\|G\|}_{\C ^{\ell-d-1}} \le \dfrac{C}{L^{\ell+1 }} {\|K\|}_{\C ^{\ell}}.
\]
where $C$ depends only on the dimension of the space.
\end{remark}

As we will use  Lemma \ref{lem:homological} in the averaging algorithm
\ref{recursive}, in the  next definition we want to  emphasize that
resonances only play a role when they are present in the Fourier transform of
the Hamiltonian at the step $N$, that is, when the numerator in the the
equation \eqref{Fouriersolution} for $G_{k,l}(I)$ is not zero.

Note that the denominator in the expression \eqref{Fouriersolution} depends
on the unperturbed system, but the numerator depends on the perturbation. The
places where the denominator  vanishes are the resonances. Clearly the
resonances do not matter unless the numerator is not zero.

\begin{de}\label{def:activated}
Given a resonance $\R_{k,l}$ as defined in definition \ref{resonantset}, we
say that it is activated at order $N$ if $N$ is the smallest value such that
$(k,l)\in \N_{N}$, were   $\N_{N}$ is  the support of the Fourier transform
of the term of order
 $\ep^N$,  after applying $N-1$ steps of the averaging algorithm
\ref{recursive}. That is, $\N_{N}$ is the support of the Fourier transform of
$K^N$ in \eqref{eq:restrictedN}.

We denote the set of resonances activated at order $N$ by
\begin{equation}\label{activatedn}
\R^{N}= \cup
_{(k,l)\in \N_{N}}\R_{k,l}
\end{equation}
and we introduce the resonances open up to order $N$ as the set
$$
\R^{[\le N]}=\cup _{i=1}^{N}\R^{i} = \cup
_{(k,l)\in \N^{[\le N]}}\R_{k,l}
$$
where
$$
\N^{[\le N]}=\cup _{i=1}^{N}\N_{i}
$$

The resonances open up to order $2$, $\R^{[\le 2]}$, are called secular
resonances.
\end{de}

We note that if a resonance has been activated at order $q$, we have that, in
a neighborhood of that resonance, the system is reduced to integrable up to
order $\ep ^{q}$. If the term of order $\ep ^{q}$ in the averaged system does
not vanish, averaging to higher order does not change the leading order term
anymore. We will only need the cases $q=1,2$.
\begin{de}
Given a frequency $\omega$ and an order of averaging $N$, we define its
active resonances up to order $N$:
\begin{equation}\label{activeesonance}
\A(\omega,N)= \{ (k,l) \in \N^{[\le N]}, \omega \cdot k +l=0\}
\end{equation}
and $\mathbf{m}(\omega,N)$, the multiplicity of $\omega$ up to  order  $N$,
as the dimension of the $\integer$-module generated by $\A(\omega,N)$:
\begin{equation}\label{multiplicity}
\mathbf{m}(\omega,N)= \mathrm{dim}( \M[\A(\omega,N)]).
\end{equation}
\end{de}
Of course $\mathbf{m}(\omega,N+1)\ge \mathbf{m}(\omega,N)$ and inequality can
be strict.

The relevance of the concept of multiplicity comes because, as we emphasize
in Remark \ref{redundantangles}, the main result of the averaging method is
that in a neighborhood of a point $I\in \J_\delta$ such that the frequency
$\omega (I)$ is of multiplicity $\mathbf{m}(\omega, N)$ up to order $N$,
there exists a change of variables that reduces the Hamiltonian $K$ in Lemma
\ref{lem:homological} to a function of $I$ and $\mathbf{m}(\omega,N)$ angles
up to an error of order $O(\ep ^{N+1})$.

Using the method of Lemma \ref{lem:homological} in the Algorithm
\ref{recursive} we obtain straightforwardly Lemma~\ref{inductivelemma}, that
is the main iterative step in the averaging procedure.

The hypothesis of next Lemma~\ref{inductivelemma} are  that we have a
Hamiltonian averaged up to order $q$. The conclusions  are that we can
produce another Hamiltonian which is averaged up to a higher order $q+1$ in
$\ep$.

\begin{lemma}\label{inductivelemma}
Consider a Hamiltonian of the form:
\begin{equation} \label{startinghamiltonian}
K_{q}(I,\vp,s;\ep)=K_{q}^{0}(I,\vp,s;\ep)+\ep ^{q+1}K_{q}^{1}(I,\vp,s;\ep).
\end{equation}
Assume that  $K_{q}$ is of class $\C ^{\ell}$ with respect to $I \in \JJ
_{\delta}\subset \II ^*\subset \real ^{d}$.

Consider the  finite collection of  sets $\R ^{s} \subset \J_\delta $, called
resonances activated at  order $s$, $s=1, \dots ,q$, and $\R^{[\le q]}=\cup
_{1\le s\le q}\R ^{s}$ the set of resonances open up to order $q$.

Consider a number $L<1$ small enough such that:
\begin{itemize}
\item [L1]\label{L} for any $\R_{k,l}, \R_{\tilde k, \tilde l} \in
    \R^{[\le q]}$, either $\R_{k,l} = \R_{\tilde k, \tilde l}$ or the
    tubular neighborhood of $\R_{k,l}$ of radius $L$ does not contain
    $\R_{\tilde k, \tilde l}$.
\end{itemize}
And assume that:
\begin{itemize}
\item
$K_{q}^{0}(I,\vp,s;\ep)$ satisfies:

\begin{itemize}
\item 
If $q=0$, $K_{0}^{0}(I,\vp,s;\ep)= I_{0}+h(I)$.
\item
If $q\ge 1$, $K_{q}^{0}(I,\vp,s;\ep)$ is a  $\C ^{n+2-2q}$ function that
verifies:
\item [1.1.]\label{noresonant} If $d(I,\R  ^{[ \le q]})\ge 2L$, then
\[
K_{q}^{0}(I,\vp,s;\ep) = I_{0}+h(I)+ \ep K_{q}^{0,0}(I;\ep),
\]
where $\ep  K_{q}^{0,0}(I;\ep)$ is a polynomial of degree $q$ in $\ep$.

\item [1.2.]\label{resonantq} If $d(I,\R  ^{[ \le q]})\le L$, then we can
    find at least one (may be more)  $0\le j \le q$ such that $d(I,\R
 ^{j})\le L$, and therefore at least one $k_i^j,l_i^j$ such that
$$
\R_{k_i^j,l_i^j} \in \R^s, \quad i=1,\dots ,n_j,
$$
and  that $d(I,\R _{k_i^j,l_i^j} ) \le L$.

Then,
\begin{equation} \label{avhamres2}
\begin{array}{rcl}
\bar K^0_q(I,\vp,s;\ep) &=& h(I)+ \ep \bar K^{0,0}(I;\ep) \\
&+& \sum _{0\le j \le q} \ep ^j U^{j,q}(I, k_{1}^j\varphi + l_{1}^j s,
\dots , k_{n_j}^j\varphi +l_{n_j}^j s;\ep)\\
\end{array}
\end{equation}
where the functions $U^{j,q}(I, \theta_1^j, \dots ,\theta _{n_j}^j;\ep)$
are polynomials in $\ep$ and trigonometric polynomials in the angle
variables $\theta_i^j$, $i=1,\dots, n_j$, with support of the Fourier
transform with respect to the $(\varphi,s)$ contained in $\N_1\cup \dots
\cup \N_q$.

Moreover, for $j=1$ the function $U^{1,q}(I, \theta_1^1, \dots ,\theta
_{n_1}^1;\ep)$ is given by:
\[
\begin{split}
U^{1,q}& =  \sum _{i=1}^{n_1} \left(\sum _{p=-N_i}^{N_i}
K^{1}_{pk_{i}^1,pl_{i}^1} (\Gamma_{k_{i}^1,l_{i}^1}(I) )e^{2\pi i p(k_{i}^1 \vp+l_{i}^1s)}\right)\\
&  + \Or (\ep)
\end{split}
\]
where  $ K^{1}_{k,l}(I)$ are the Fourier coefficients of
$K^{1}_{0}(I,\vp,s;0 )$ with respect to the angle variables $(\vp,s)$.
\end{itemize}
\item [2.] $ K_{q}^{1}(I,\vp,s;\ep) $ is a  $\C ^{n-2q}$ function whose
    Taylor series coefficients with respect to $\ep$ are trigonometric
    polynomials in $(\vp,s)$.
\end{itemize}
Denote by
$$
K = K^1_q(I, \vp, s; 0 )= \sum _{(k,l) \in \N_{q+1}}
K_{q}^{k,l}(I) e^{2\pi i (k\vp+l s)},
$$
where  $\N_{q+1}$ is assumed to be a finite set. $ K $  is the term of the
perturbation of order exactly $q+1$. Introduce also the set of resonances
activated at order $q+1$:
\begin{equation}
\label{resonancesetq}
 \R^{q+1} = \cup _{(k,l)\in \N^{q+1}}\R_{k,l} \setminus \R ^{[\le q]}
\end{equation}
Choose $0<\tilde L <L$ such that $L)$ holds for $\R^{[\le q+1]}$.

Let  $G(I,\vp,s)$ be the $\C ^{n-2q-1}$ function whose Fourier coefficients
$G_{k,l}$ verify equation~\eqref{averagingstepFourier}, for $(k,l)\in
\N_{q+1}$,  with $K=K_{q}^{1}(I,\vp,s;0)$.

Then,   the   $\C ^{n-2q-2}$  change of variables
\[
(I,\vp,s) = g(\B,\alpha,s),
\]
given by the time one flow of the Hamiltonian $\ep ^{q+1}G(\B,\alpha,s)$,
transforms the Hamiltonian $K_{q}(I,\vp,s;\ep)$ into a Hamiltonian
\[
K_{q+1}(\B,\alpha,s;\ep)=K_{q+1}^{0}(\B,\alpha,s;\ep)+\ep ^{q+2}K_{q+1}^{1}(\B,\alpha,s;\ep),
\]
with
\begin{equation}\label{nouk0}
K_{q+1}^{0}(\B,\alpha,s;\ep)=K_{q}^{0}(\B,\alpha,s;\ep)+\ep ^{q+1}\bar K_{q}^{1}(\B,\alpha,s;0),
\end{equation}
where  $\bar K_{q}^{1}(\B,\alpha,s;0)=\bar K(\B,\alpha,s)$,
 given in Lemma~\ref{lem:homological},
is a  $\C ^{n-2q}$ function whose Fourier coefficients solve
equations~\eqref{averagingstepFourier}.

Moreover, the Hamiltonian $K_{q+1}^{0}(\B,\alpha,s;\ep)$
 verifies properties~[1.0], [1.1], [1.2]  up to order $q+1$
with $\tilde L$ replacing $L$.

Furthermore, $\ep^{q+2}K^{1}_{q+1}(B,\alpha,s;\ep)$ is a $\C ^{n-2q-2}$
function whose Taylor series coefficients with respect to $\ep$ are
trigonometric polynomials in $(\alpha,s)$.
\end{lemma}

Once we know how to solve any homological
equation~\eqref{averagingstepFourier}, we can proceed to obtain a suitable
global normal form of our reduced Hamiltonian by applying repeatedly the
procedure. The precise  result is formulated  in the following
Theorem~\ref{averaging1}, which is a straightforward generalization of
Theorem 8.9 in \cite{DelshamsLS06a}.

\begin{theorem}\label{averaging1}
Let  $\tilde K(\tilde I,\vp,s;\ep)$ be  a $\C^{n}$ Hamiltonian, $n > 1$, for
$I\in \JJ _{\delta}\subset \II ^*\subset \real ^{d}$ and consider any $1\le
m<n$, independent of $\ep$. Assume that
\begin{equation}\label{eq:hamilt}
\tilde K (\tilde I,\vp,s;\ep) = I_{0}+ h(I)+  \ep K(I,\vp,s;\ep).
\end{equation}
Let $K_{i}(I,\vp,s)$ $i = 1, , \dots ,m$ be the coefficients in the Taylor
expansion with respect to $\ep$ of $K(I,\vp,s;\ep)$, and assume that the
$K_{i}(I,\vp,s)$, $i = 1,\dots , m$ are trigonometric polynomials in $\vp,
s$.

Consider the  finite collection of sets $\R ^{i} \subset \II^*$, called
resonances activated at  order $i$, $i=1 \dots m$, following definition
\ref{def:activated}, as well as the resonances open up to order $m$:
$\R^{[\le m]}= \cup  _{1=1,\dots , m}\R ^{i}$.

Consider a number $0< L<1$ small enough such that:
\begin{itemize}
\item[L1] for any $\R_{k,l}, \R_{\tilde k, \tilde l} \in \R^{[\le m]}$,
    either $\R_{k,l} = \R_{\tilde k, \tilde l}$ or the tubular
    neighborhood of $\R_{k,l}$ of radius $L$ does not contain $\R_{\tilde
    k, \tilde l}$.
\end{itemize}

Then,  there exists a symplectic change of variables, depending on time, $
(I, \vp,s) \mapsto (\B, \alpha,s)$, periodic in $\vp$ and $s$, and  of class
$\C ^{n-2m}$, which is $\ep$-close to the identity in the $\C ^{n-2m-1}$
sense, such that transforms the Hamiltonian system associated to $\tilde
K(\tilde I,\vp,s;\ep)$ into a Hamiltonian system of Hamiltonian
\[
B_{0}+ \bar K(\B,\alpha,s ;\ep) = B_{0}+\bar K^0(\B,\alpha,s;\ep)
+\ep^{m+1} \bar K^{1}(\B,\alpha,s;\ep)
\]
where the function $\bar K^{0}$ is of class $\C ^{n-2m+2}$, and $\ep
^{m+1}\bar K^{1}$ is of class  $\C ^{n-2m}$, and they verify:

\begin{enumerate}
\item If $d(\B,\R ^{[\le m]})\ge 2L$, then
\[
\bar K^0(\B,\alpha,s;\ep) = h(\B)+ \ep \bar K^{0,0}(\B;\ep)
\]
where $\bar K^{0,0}(\B;\ep)$ is a polynomial of degree $m-1$ in $\ep$.
\item
If $d(\B,\R  ^{[ \le m]})\le L$, then we can find at least one (may be
more) $0\le j \le m$ such that $d(\B,\R  ^{j})\le L$, and therefore at
least one $k_i^j,l_i^j$ such that
$$
\R_{k_i^j,l_i^j} \in \R^j, \quad i=1,\dots n_j,
$$
and that $d(\B,\R _{k_i^j,l_i^j} ) \le L$.

Then,
\begin{equation} \label{avhamres3}
\begin{array}{rcl}
\bar K^0(\B,\alpha,s;\ep) &=& h(\B)+ \ep \bar K^{0,0}(\B;\ep) \\
&+& \sum _{0\le j \le m}\ep ^j U^{j,m}(\B, k_{1}^j\alpha + l_{1}^j s,
\dots , k_{n_j}^j\alpha +l_{n_j}^j s;\ep)
\end{array}
\end{equation}
where the functions $U^{j,m}(\B, \theta_1^j, \dots ,\theta _{n_j}^j;\ep)$
are polynomials in $\ep$ and  trigonometric polynomials in  the angle
variables $\theta_i^j$, $i=1,\dots, n_j$, with support of the Fourier
transform with respect to the $(\varphi,s)$ contained in $\N_1\cup \dots
\cup \N_m$.

Moreover, if $j=1$, the function $U^{1,m} (I, \theta_1^1, \dots ,\theta
_{n_1}^1;\ep)   $ is given by:
\[
\begin{split}
U^{1,m}=  \sum _{i=1}^{n_1} (\sum _{p=-N_i}^{N_i}  &
 K^{1}_{pk_{i}^1,pl_{i}^1} (\Gamma_{k_{i}^1,l_{i}^1}(\B))e^{2\pi i p(k_{i}^1 \alpha+l_{i}^1s)})\\
&  + \Or (\ep)
\end{split}
\]
where  $K^{1}_{k,l}(\B )$ are the Fourier coefficients of the
$K(\B,\alpha,s;0 )$ with respect to the angle variables $(\alpha,s)$.
\end{enumerate}
\end{theorem}

Note that in this Theorem  \ref{averaging1} we have not claimed anything in
the regions at a distance between $L$ and $2L$ of the resonance set $\R$.
This is not a problem because, by remembering that $L$ is arbitrary, we can
obtain the same results using $L/2$ in place of $L$.

Hence, the analysis that we will carry out in each of the different pieces
applies to the whole space.

\subsection{Averaging close to simple resonances}

Theorem \ref{averaging1} reduces the original system to a non-integrable
system of a partially simple form by a change of variables. The number of
angles that enter in the averaged Hamiltonian for a given $I$ depends on the
resonances which are close to $I$. The key of the following reasoning is to
understand the geometry of the set of points $I$ for which the averaged
system involves  only one angle. The set of points that we have found useful
to omit are the points in secular resonances (i.e. in the resonances which
appear in averaging to order $1$ or $2$) which are also part of another
resonance activated when averaging up to order $m$.

So, next step is to define a region
 $\II_{\delta}\subset \JJ_\delta$ where we take out the intersection of the secular resonances with
any other resonances which appear in the process of averaging up to order
$m$.

To this end we consider
$$
\mathbf{B}= \R^{[\le 2]} \cap \{ I \in \JJ_{\delta}, \ \mathbf{m}(\omega(I),m) \ge 2\}
$$
which is a finite union of surfaces of codimension two or higher in
$\JJ_{\delta}$, and we consider $\mathbf{B_{\delta}}$ a $\delta$-neighborhood
of these surfaces. Reducing $L=L(\delta)$ if necessary,  the set
$\II_{\delta}= \JJ_{\delta} \setminus \mathbf{B_{\delta}}$ verifies the
following property:

\begin{itemize}
\item[L2] If $I\in \II_{\delta}$ there is at most one  $\R_{k,l} \in
    \R^{[\le
2]}$ such that \\
$d(I,\R_{k,l})\le L$.
\end{itemize}
Theorem \ref{averaging1} in the domain $\II_{\delta}$ reads:

\begin{theorem}\label{averaging2}
Let  $\tilde K(\tilde I,\vp,s;\ep)$ be  the  $\C^{n}$ Hamiltonian of Theorem
\ref{averaging1}, $n> 1$, and consider any $1\le m<n$, independent of $\ep$.

Consider the finite collection of   sets $\R ^{i} \subset \II^*$, called
resonances activated at  order $i$, $i=1, \ldots m$, given in definition
\ref{def:activated}.

Let $0<\delta <1$ be any number and consider $0<L<1$  verifying $L1$ and $L2$
in $\II_{\delta}$.

Then,  the  symplectic change of variables given in Theorem \ref{averaging1},
$$
(I, \vp,s)  \mapsto (\B, \alpha,s),
$$
transforms the Hamiltonian system associated to $\tilde K(I,\vp,s;\ep)$ in
$\II_{\delta}$ into a Hamiltonian system of Hamiltonian
\begin{equation}\label{generalterm}
B_{0}+ \bar K(\B,\alpha,s ;\ep) = B_{0}+\bar K^0(\B,\alpha,s;\ep)
 +\ep^{m+1} \bar K^{1}(\B,\alpha,s;\ep)
\end{equation}
where the function $\bar K^{0}$ is of class $\C ^{n-2m+2}$, and
 $\ep ^{m+1}\bar K^{1}$ is of class  $\C ^{n-2m}$ and they verify:

\begin{enumerate}
\item \label{item1thm17} If $\B \in \II_{\delta}$, satisfies $d(\B,\R
    ^{[\le 2]})\ge 2L$, then
\[
\bar K^0(\B,\alpha,s;\ep) = h(\B)+ \ep \bar K^{0,0}(\B;\ep)+ O(\ep ^3)
\]
where $\bar K^{0,0}(\B;\ep)$ is a polynomial of degree $1$ in $\ep$.

\item
If $\B \in \II_{\delta}$, satisfies $d(\B,\R ^{[ \le 2]})\le L$, there
exists a unique resonance activated at order one or two
$$
\R_{k_0,l_0} \in \R^j, \quad j=1,2
$$
such that $d(\B,\R _{k_0,l_0} ) \le L$. Then
\begin{equation} \label{avhamres4}
\begin{array}{rcl}
\bar K^0(\B,\vp,s;\ep) &=& h(\B)+ \ep \bar K^{0,0}(\B;\ep) \\
&+& \ep^j U^{k_0,l_0}(\Gamma _{k_0,l_0}(\B), k_{0}\cdot \alpha + l_{0} s;\ep)
\end{array}
\end{equation}
where the functions $U^{k_0,l_0}(\Gamma _{k_0,l_0}(\B), \theta;\ep)$ are
polynomial in $\ep$ and trigonometric polynomial in the angle variable
$\theta=k_0 \cdot \alpha + l_0 s$.

Moreover, if $\R _{k_0,l_0} \subset \R^{1}$ the function $U^{k_0,l _0}$
is given by:
\[
U^{k_0,l_0}=  \sum _{p=-N_1}^{N_1}
K^{1}_{pk_{0},pl_{0}} (\Gamma_{k_{0},l_{0}}(\B) )e^{2\pi i p(k
_{0} \cdot \alpha+l_{0}s)}  + \Or (\ep)
\]
where  $ K^{1}_{k,l}(\B)$ are the Fourier coefficients of
$K(\B,\alpha,s;0 )$ with respect to the angle variables $(\alpha,s)$.
\end{enumerate}
\end{theorem}

The next goal is to study in more detail the behavior of the system predicted
by the averaged Hamiltonian.

The main remark is that, near simple resonances, the averaged system contains
only one angle and, therefore, it is integrable. This allows us to analyze
explicitly its  dynamics. Its turns out that, for the problem at hand, we
only need to study the resonances of order $1$ or $2$, which are called
``secular resonances" by astronomers.

\subsection{Geometric properties of the orbits of
the averaged Hamiltonian}\label{averaged}

In this section, we study the invariant tori of the averaged system obtained
in Theorem~\ref{averaging2}, that is, the system given by Hamiltonian $\bar
K^0$ in \eqref{generalterm}. Later, in section \ref{sec:non-resonantKAM},
\ref{sec:toriKAM}, we show that, under some non-degeneracy conditions, some
of these tori are also present in the original system. This is, basically,
the KAM theorem.

In Sections \ref{sec:non-resonant}, \ref{sec:resonant}, we will see that the
phase space $\II_\delta \times \torus^{d+1}$ is foliated by (quasi)-periodic
solutions of the averaged system. Nevertheless, the topology of the solutions
is very different in the non-resonant regions and in the resonant regions. We
define the {\emph {non-resonant region}} as the set:
\begin{equation}\label{non-resonantregion}
\S ^{L} =\{ (I,\vp,s) \in \II_{\delta}\times \torus ^{d+1}, \ d(I, \R^{[\le 2]}) \ge 2L\}
\end{equation}
In particular, $\S ^{L}$ includes the intersection of $\II_{\delta}$ with all
the resonances activated at order higher than 2. This region $\S ^{L}$ will
be  covered, up to very small gaps of order $O(\varepsilon ^{3/2})$, by KAM
tori.

In the resonant regions of $\II_{\delta}\times \torus ^{d+1}\setminus \S
^{L}$, we will obtain tori which are contractible to tori of lower dimension
and, therefore, are not homotopic to a torus present in the unperturbed
system. We call \emph{secondary} KAM tori the invariant tori which have
different topological type from the tori of the unperturbed system. We use
the name \emph{primary} tori for the invariant tori which are homotopic to
those of the unperturbed system. Primary tori are those usually considered in
the perturbative versions of KAM theorem for quasi-integrable systems.

The importance of the secondary tori is that they dovetail precisely into the
gaps between the set of KAM primary tori created by the resonances, so that
it will be possible to construct a web  of KAM tori, primary and secondary,
which are $\ep^{3/2}$-close.
 For systems with 2+1/2 degrees of freedom this was
introduced in \cite{DelshamsLS03a, DelshamsLS06a}. See also \cite{GideaL06a}.

In hypothesis {\bf H6} given in \ref{nondegeneracy7}, we formulate precisely
one non-degeneracy assumption on the averaged system which allows us to apply
the KAM theorem and conclude that  some of the solutions found in the
averaged system $\bar K^0$ (including secondary tori) are indeed present in
the full Hamiltonian \eqref{generalterm}, and therefore in the original
system \eqref{eq:hamiltonianrestricted}. Of course, since the averaged system
is computable from the original model, the non-degeneracy conditions on the
averaged system amount to some non-degeneracy conditions on the original
system.

\subsubsection{The  invariant tori  of
the averaged system in the non-resonant region of $\II_{\delta}$}
\label{sec:non-resonant}

By item (\ref{item1thm17}) in Theorem \ref{averaging2}, in the non-resonant
region $\S^{L}$ defined in \eqref{non-resonantregion}, the full averaged
hamiltonian \eqref{generalterm} reads
\[
\B_0+  h(\B)+ \ep \bar K^{0,0}(\B;\ep)+ O(\ep ^3).
\]
For the truncated Hamiltonian
\begin{equation}\label{averagednon-resonant}
\B_0+ h(\B)+ \ep \bar K^{0,0}(\B;\ep)
\end{equation}
the tori are given as the level sets of the averaged action variables

\[
\B_0 = c_{0}, \B_{1}= c_{1}, \ \dots , \B_{d}= c_{d}
\]
where the equation $\B_0=c_{0}$ is a reflection of the fact that the
Hamiltonian \eqref{averagednon-resonant} is autonomous.

When written in the original variables of the time dependent Hamiltonian
\eqref{eq:hamiltonianrestricted}, these tori in $\S^{L} \subset \II _{\delta}
\times \torus ^{d+1}$, as shown in Theorem~\ref{averaging1} and
Theorem~\ref{averaging2}, are given by the equations:
$$
F_{1}(I,\varphi,s;\ep) = c_{1}, \dots , F_{d}(I,\varphi,s;\ep) = c_{d}
$$
where
\begin{equation}
F_{1}(I,\varphi,s;\ep) = I_{1}+O(\ep),  \dots ,F_{d}(I,\varphi,s;\ep) = I_{d}+O(\ep) .\label{eq:primarytorus1}
\end{equation}

\subsubsection{The  invariant tori  in the non-resonant region of $\II_{\delta}$: KAM Theorem}
\label{sec:non-resonantKAM}

We note that in the non-resonant region $\S ^L$, we have managed to transform
the system into an integrable system up to an error which is $\ep^3$ when
measured in the $C^{r - 3}$ norm.

Furthermore, we point out that the averaged part has a frequency map which is
a diffeomorphism (it is an $O(\ep)$ perturbation of the diffeomorphism $I
\mapsto \frac{\partial h}{\partial I}$ in a smooth norm).

If $r$ is sufficiently large (so that $r -3$ is larger than $2 d + 3$) we can
apply a KAM theorem \cite{Poschel82} and conclude that there are invariant
tori which cover the non-resonant region $\S^L$ except for a set of measure
smaller than $O(\ep^{3/2})$.

\begin{theorem}\label{KAMnoresonant}
Under the conditions of theorem \ref{main}, there exists $\varepsilon_0$ such
that, for $0<|\varepsilon|<\varepsilon_0$, the region $\S^L$ can be covered
by $O(\ep^{3/2})$ neighborhoods of invariant objects under the Hamiltonian
flow of the Hamiltonian $K_{\ep}(I,\vp,s)$ in
\eqref{eq:hamiltonianrestricted}. Moreover:
\begin{itemize}
\item These invariant objects are given by the level sets $F=E$, for
    $|E-E'| \le \ep ^{ \frac{3}{2}}$.
\item The $\C^{2}$ function $F: \real^d\times \torus^d \times \torus \to
    \real$ is given
by \eqref{eq:primarytorus1}.
\item These invariant objects are  regular primary KAM
    $d+1$-tori\end{itemize}
\end{theorem}

Therefore, in the non-resonant region, each torus has several tori which  are
much closer than $O(\ep)$ to it. This is what in \cite{ChierchiaG98} was
called the \emph{``gap bridging mechanism''}.

\begin{remark} \label{rem:moreaveraging}
For experts, we note that there are different KAM theorems in the literature,
which differ in some subtle features; a systematic comparison can be found in
\cite{Llave01a}. The main difference in the literature is whether one step of
averaging requires to solve one cohomology equation or two. The methods which
use only one cohomology equation (e.g. the method in \cite{Arnold63a,
Poschel82, Salamon04}), called \emph{first order methods}, establish that the
gaps between tori are bounded by the  error to the power $1/2$.  Those that
use two cohomology  equations (e.g. the methods in \cite{Kolmogorov79,
Moser66a, Zehnder75,Zehnder76}), called \emph{second order methods}, lead to
gaps which are bounded by the error to power $1/4$. These quantitative
estimates for the Newton method are found in \cite{Zehnder76b} Very explicit
verifications of the quantitative estimates of the method of \cite{Arnold63a}
appear in \cite{Neishtadt81}. Simple examples show that the exponent $1/2$
cannot be increased.

In our case, either method could be applied. If we wanted to just refer to
the second order methods to obtain gaps of order $O(\varepsilon ^{1/4})$, it
would have been enough, to obtain gaps of order $O(\varepsilon^{3/2})$, to
define the non-resonant region as the region where one can average to order
$m=5$ instead of $m=2$ (in fact, it would be enough $m=4$ if we allow  gaps
of order $O(\varepsilon ^{5/4})$).

Another technical point is that some results in the literature lose more
derivatives. This is totally irrelevant for us since it only affects the
number of derivatives that we need to assume in the original Hamiltonian.
\end{remark}

\subsubsection{Primary and secondary invariant tori  of
the averaged system in the resonant region} \label{sec:resonant}

We define the \emph{secular resonant region} as:
\begin{equation}\label{secularregion}
\S ^{[\le 2]} =\{ (I,\vp,s) \in \II_{\delta}\times \torus ^{d+1}, \ d(I, \R^{[\le 2]}) \le L\}
\end{equation}

In this region we will  perform an elementary change of variables that makes
it clear that, close to a resonance activated at order one or two, the
Hamiltonian is a function of $d-1$ actions and one resonant angle.

The first observation is that the region $\S^{[\le 2]}$ is the union of the
regions $\R^{L}_{k,l}$, which consist on tubular neighborhoods of size $L$ of
the resonances $\R_{k,l}$ defined in Definition \ref{resonantset}.

For points  $(\B,\alpha,s) \in \R^{L}_{k_0,l_0}$, where
$$
\R_{k_{0},l_{0}} \subset \R^{j}\subset \R^{[\le 2]},\quad j=1,2
$$
the averaged Hamiltonian is
\begin{equation}\label{nonautonomousaveraged}
\B_{0}+ \bar K^0(\B,\alpha,s;\ep)
\end{equation}
 where, by Theorem \ref{averaging2}:
\begin{equation}\label{eq:averagedresonant}
\bar K^0(\B,\alpha,s;\ep) = h(\B)+ \ep \bar
 K^{0,0}(\B;\ep) + \ep ^{j}U^{k_0,l_0}(\Gamma_{k_0,l_0} (\B), k_{0} \cdot \alpha +l_{0}s;\ep)
\end{equation}
is given in \eqref{avhamres4}.

Assuming $k_{0}^{\j}\ne 0$ for some $1\le \j \le d$, to understand the
geometry of the averaged Hamiltonian we first perform the following change of
angles
\begin{equation}\label{eq:changesimplectic}
\tilde \theta = M\tilde \alpha ,
\end{equation}
where $\tilde \theta= (s,\theta)$, and $\tilde \alpha = (s,\alpha)$, and $M$
is the $d+1 \times d+1$ matrix:
\[
M=\left(\begin{array}{ccc}1&\dots&  \\ \cdot &\cdot & \\ l_{0} & k_{0}^{\top}\\\cdot &\cdot & \\ & & 1 \\
\end{array} \right)
\]

Let us observe that this change is just to take as a new angle the resonant
angle:
$$
\theta_i=\alpha_i, \quad i\ne \j , \quad   \theta_\j= k_{0}\cdot\alpha +l_{0}s
$$

To make the change symplectic we perform the change in actions:
\begin{equation}\label{eq:changesimplecticactions}
\tilde J = M^{-\top}\tilde \B .
\end{equation}

The change  $\tilde \B= M^{\top}\tilde J$ is equivalent to $
 \B= N ^{\top} J
$ where
\[
N=\left(\begin{array}{ccc}1&\dots&  \\ \cdot &\cdot & \\  & k_{0}^{\top}\\\cdot &\cdot & \\ & & 1 \\
\end{array} \right)
\]
and $\B_{0}= J_{0}+l_{0}J_{\j}$, which is just the change in action which
reflects the fact that we are doing a time dependent change of angles.

In components, this change is simply:
$$
\B_\j= k_0^\j J_\j, \quad \B_i= J_i + k_0^i J_\j ,\quad i\ne \j
$$
or, analogously
\begin{equation}\label{changeJ}
J_\j= \frac{ \B_\j}{k_0^\j}, \quad J_i= \B_i - \frac{\B_\j}{k_0^\j}k_0^i  ,\quad i\ne \j
\end{equation}

With this change, the  averaged Hamiltonian \eqref{nonautonomousaveraged}  is
given by:
\begin{equation}\label{averagedJ}
 J_{0}+ l_{0} J_{\j}+h(N^{\top} J)+ \ep \bar
 K^{0,0}(N^{\top} J;\ep) + \ep ^{j}U^{k_0,l_0}(\Gamma _{k_0,l_0}(N^{\top} J), \theta_{\j};\ep)
\end{equation}
which, in the region $\R ^{L}_{k_0,l_0} \subset \II_\delta \times \torus
^{d+1}$, corresponds to the autonomous Hamiltonian:
\begin{equation}\label{averagedJextended}
l_{0} J_{\j}+h(N^{\top} J)+ \ep \bar
 K^{0,0}(N^{\top} J;\ep) + \ep ^{j}U^{k_0,l_0}(\Gamma _{k_0,l_0}(N^{\top} J), \theta_{\j};\ep).
\end{equation}

Working in the variables $(J,\theta,s)$ makes easier to identify the
invariant tori. The invariant tori of Hamiltonian \eqref{averagedJextended}
will be given by prescribing the values of the $d-1$ action variables $J_i$,
for $i\ne \j$, and the value of the Hamiltonian \eqref{averagedJextended},
which is  a constant of motion.

Abusing slightly the notation, let us write $J=(\hat J, J_\j)$,  $\theta
=(\hat \theta, \theta_{\j})$, with
$$
\hat J= (J_{1}, \dots , J_{\j-1}, J_{\j+1}, \dots , J_{d}), \ \hat \theta=
(\theta_{1}, \dots ,\theta_{\j-1}, \theta_{\j+1}, \dots , \theta_{d}).
$$

Given a value of $\B = N^{\top}J = (\hat \B, \B_\j)$, we want to compute its
projection $\Gamma_{k_0,l_0}(\B)=\B ^*_{k_0,l_0}= N^{\top}J^*$.

By the $k_0$-characterization~\eqref{projection}
 of the projection $\Gamma_{k_0,l_0}$   we have that
 $$
 \B-\B ^*_{k_0,l_0} \in\langle k_0\rangle \mbox{ and }
 Dh(\B ^*_{k_0,l_0})\cdot k_{0}+ l_{0}=0.
 $$
Note that using the  form of the change \eqref{changeJ}, and since
$\B=N^{\top} J$ and $\B ^*_{k_0,l_0}=N^{\top}J^*$, we have that $J_i=J^*_i$,
if $i\ne \j$ and that $\B-\B ^*_{k_0,l_0} =  (J_\j -J^*_\j) k_0$.

Then, given $\B=N^{\top}J$, with $J=(\hat J, J_{\j})$, one can characterize
the projection $\Gamma_{k_0,l_0}(\B)$ as:

We  compute:
\begin{equation}\label{eq:projection}
\hat J= \hat B-\frac{\B_\j}{k_0^\j}\hat k_0 \ , \quad J_{\j}= \frac{\B_\j}{k_0^\j}
\end{equation}
and therefore
$$
\B= (\hat J,0)+ \frac{\B_\j}{k_0^\j} k_0,
$$
one can obtain the projection $\B^*_{k_0,l_0}=N^{\top}J^*$, with $J^*=(\hat
J, J_{\j}^*)$ in terms of $\hat J$:
\begin{equation}\label{eq:bestrella}
\B ^*_{k_0,l_0}= \B^*(\hat J)=(\hat J,0)+ \frac{\B_\j^*}{k_0^\j} k_0
\end{equation}
where $\B_\j^*= k_0^\j J_\j^*$ and $J_\j ^*= J_\j ^* (\hat J)$ is obtained
solving
\begin{equation}\label{eq:resonancia}
Dh(N^{\top} (\hat J, J^*_{\j}))\cdot k_{0}+ l_{0}=0.
\end{equation}

Let us emphasize that $J_\j^*$ and therefore $J^*$ and $\B ^*_{k_0,l _0}$ are
uniquely determined through \eqref{eq:resonancia}
 for a given $\B=N^{\top}J$,
and therefore for a given $\hat J$, assuming $\dist (\B, \R _{k_0,l_0}) <L$
and $L$  small enough.

Let us observe that for values $\B$ such that $(\B,\alpha,s)\in \II_\delta$,
one has that the corresponding values of $\hat J$ vary in a compact set that
we will denote by $\hat \J$.

For $\hat J \in \hat \J$, we will denote:
\begin{equation}\label{Uestrella}
U^{k_0,l_0, *}(\theta_{\j};\hat J,\ep)   =
U^{k_0,l_0}(\B^*_{k_0,l_0}(\hat J), \theta_{\j};\ep)
\end{equation}

Using this notation, system \eqref{averagedJextended} can be written as:
\begin{equation}\label{averagedJextendednew}
l_{0} J_{\j}+h(N^{\top} (\hat J, J_{\j}))+ \ep \bar
 K^{0,0}(N^{\top} (\hat J, J_{\j});\ep)+ \ep ^{j}U^{k_0,l_0,*}( \theta_{\j};\hat J,\ep).
\end{equation}

The next step is  to use Hypotheses {\bf H5} and {\bf H6} to obtain a change
of variables
$$
(J_\j, \theta _\j) \to (y,x)
$$
such that Hamiltonian \eqref{averagedJextendednew} in these new variables
becomes:

\begin{equation}\label{averagedxy}
\K_0(y,x;\hat J,\ep) = a(\hat J, \ep) \frac{y^2}{2} (1+O(y))+ \ep ^j U(x;\hat J,\ep).
\end{equation}

We first proceed to formulate hypothesis   {\bf H5} which states precisely
that the leading part of the kinetic energy is quadratic.

Taylor expanding the function $h(N^{\top}J)$ around the resonant point
$N^{\top}J^*$ and using \eqref{eq:resonancia}, we obtain:
\begin{equation}\label{taylorexpansion}
\begin{split}
l_{0} J_{\j}+ h(N^{\top} J)
&=  l_{0} J_{\j}+ h(N^{\top} J^*+ k_0(J_\j-J_\j^*)) \\
&= l_{0} J_{\j}^{*}+h\left(N^{\top}J^{*}\right)\\
&+\frac{1}{2} (J_{\j}-J_{\j}^{*})^{2} k_{0}^{\top}D^{2}h\left(N^{\top}J^{*}\right)k_{0} \\
&+ O((J_{\j}-J_{\j}^{*})^{3})
\end{split}
\end{equation}
where we have used that
$$
Dh(\B^*_{k_0,l_0}(\hat J))k_0+l_0=0, \quad \forall \hat J \in \hat \J.
$$
by the definition of $\B^*_{k_0,l_0}(\hat J) =N^{\top} J^*$ in
\eqref{eq:bestrella}. Therefore, introducing
\begin{equation}\label{coefquasiconvexity}
\a=\a_{k_0,l_0}:= k_{0}^{\top}D^{2}h\left(N^{\top}J^{*}\right)k
_{0}= k_{0}^{\top}D^{2}h\left(\B^*_{k_0,l_0}(\hat J)\right)k_{0},
\end{equation}
the hypothesis {\bf H5} is:
\begin{itemize}
\item[{\bf H5}] \label{nondegeneracy6} $\a \ne 0$, for any $(k_0,l_0)\in
    \R^{[\le 2]}$ and $\hat J \in \hat \J$.
\end{itemize}
With the notation \eqref{coefquasiconvexity},
equation~\eqref{taylorexpansion} becomes
\begin{equation}\label{taylorexpansion1}
l_{0} J_{\j}+ h(N^{\top} J)
= l_{0} J_{\j}^{*}+h(N^{\top}(\hat J, J_{\j}^{*}))
+\frac{ \a}{2}(J_{\j}-J_{\j}^{*})^{2}
+ O((J_{\j}-J_{\j}^{*})^{3}).
\end{equation}

Since the actions $\hat J$  are $d-1$ first integrals of the averaged
Hamiltonian \eqref{averagedJextended}, we have that the dynamics in the
$(J_\j, \theta_\j)$ variables is that of a nonlinear oscillator with
potential $U^{k_0,l_0,*}(\theta_\j;\hat J,\ep)$. We can think of the
variables $\hat J$ as parameters in the non-linear oscillator.


We now  introduce Hypothesis {\bf H6}, which is  the other  non-degeneracy
assumption which will make precise the heuristic notion  that  \emph{``the
averaged system near secular resonances looks like a pendulum''}.

Assumption {\bf H6} formulates precisely that the potential of the truncated
averaged Hamiltonian \eqref{eq:averagedresonant} at the resonance (see also
\eqref{averagedJextendednew}), if $a(\hat J)>0$, has a unique non-degenerate
maximum or, if $a(\hat J)<0$, has a unique non-degenerate minimum.
\begin{itemize}
\item[{\bf H6}] \label{nondegeneracy7} For any $ \B = N^{\top}J \in
    \S^{[\le 2]} \subset \II_\delta$, consider the value $(k_0,l_0)$ such
    that $d(\B, \R _{k_0,l_0})\le L$ and its $k_0$-projection
    $\Gamma_{k_0,l_0}(\B)=\B^{*}_{k_0,l_0}(\hat J)$. By hypothesis {\bf
    H5} we know that $\a=\a_{k_0,l_0}\ne 0$.

If $a(\hat J)>0$, we assume that there is a unique non-degenerate maximum
of the potential of Hamiltonian \eqref{eq:averagedresonant}
\[
U^{k_0,l_0}(\Gamma_{k_0,l_0}(\B), \theta_\j;0)=
U^{k_0,l_0}(\B^*_{k_0,l_0} (\hat J), \theta_\j;0)
\]
with respect to $\theta_\j$, which is uniformly non-degenerate with
respect to $\B \in \S^{[\le 2]}$. If $a(\hat J)<0$, we assume instead
that there is a unique non-degenerate minimum with the same uniformity
conditions.

That is,  there is a unique $\theta_\j^*$ such that
\begin{equation}\label{H7}
\begin{split}
& a(\hat J)U^{k_0,l_0}(\Gamma_{k_0,l_0}(\B), \theta_\j^* ;0)
=
\max_{\theta_\j}
a(\hat J)U^{k_0,l_0}(\Gamma_{k_0,l_0}(\B), \theta_\j;0) , \\
& a(\hat J)\frac{\partial^2}{\partial \theta_\j^2}
 U^{k_0,l_0}(\Gamma_{k_0,l_0}(\B), \theta_\j^*;0)
\le \beta < 0 .
\end{split}
\end{equation}
\end{itemize}

\begin{remark}
We call attention to the fact that, as $\Gamma_{k_0,l_0}(\B)\in \R_{k_0,l_0}
\subset \R^{[\le 2]}$ the conditions {\bf{H5}}, {\bf{H6}} need to be verified
only on the codimension one set  $\R ^{[\le 2]} \cap \II_\delta \times \torus
^{d+1}$ formed by the secular resonances in $\II_\delta \times \torus
^{d+1}$.

Note that the assumptions {\bf{H5}}, {\bf{H6}} are  $\C^r$ open conditions in
the space of Hamiltonians. If assumptions {\bf{H5}}, {\bf{H6}} are  verified
for a family, they will also be verified for all families close to it in a
$\C^r$ topology with $r\ge 2$ sufficiently large so that we can carry out the
averaging procedure. Therefore, the conditions also hold in a $\C^r$ set of
original Hamiltonians.
\end{remark}

Using the notation introduced in  \eqref{Uestrella} hypothesis {\bf H6} can
be written as:
\begin{equation}
\begin{split}
&a(\hat J) U^{k_0,l_0, *}( \theta _{\j}^*;\hat J,0)   =
\max_{\theta _{\j}} a(\hat J) U^{k_0,l_0,*}( \theta_{\j};\hat J,0), \quad  \\
& a(\hat J) \frac{\partial ^2}{\partial \theta _{\j}^{2}} U^{k_0,l_0, *}( \theta _{\j}^*;\hat J,0)\le \beta <0 .
\end{split}
\end{equation}

{From} now on we will assume that $a(\hat J)>0$. Moreover, it is uniform
respect to $\hat J \in \hat \J$. Therefore hypothesis {\bf H6} can be written
as:
\begin{equation}\label{H7ambhatJ}
\begin{split}
&U^{k_0,l_0, *}(\theta _{\j}^*;\hat J, 0)   = \max_{\theta _{\j}}
U^{k_0,l_0,*}(\theta_{\j};\hat J, 0), \\
& \frac{\partial ^2}{\partial \theta
_{\j}^{2}} U^{k_0,l_0, *}(\theta _{\j}^*;\hat J, 0)\le \beta <0 .
\end{split}
\end{equation}

The case $a(\hat J)<0$ can be done analogously.

Assumptions {\bf{H5}}, {\bf{H6}} imply that, as a function of
$(J_\j,\theta_\j)$, for  any value of $\hat J$, the Hamiltonian
\[
l_{0} J_{\j}+h( N^{\top} (\hat J, J_\j))+ \ep ^{j} U^{k_0,l_0,*}(\theta _\j;\hat J, 0)
\]
has a saddle   point  at $(J_\j ^*(\hat J), \theta_\j ^*(\hat J))$, which
gives rise to a saddle equilibrium point for the associated Hamiltonian
system.

We note that, because of uniformity of the second derivative of the potential
in \eqref{H7ambhatJ}   and the hypothesis {\bf{H5}}, we obtain that the point
$(J_\j ^*(\hat J), \theta_\j ^*(\hat J))$ is uniformly hyperbolic. For a
given $\ep >0$, the  Liapunov exponents are bounded away from zero for any
$\hat J\in \hat \J$  uniformly, and the angle between its  stable and
unstable directions is also bounded away from zero.

Therefore, by the implicit function Theorem, for $|\ep |< \ep _0$, the
Hamiltonian system associated to Hamiltonian \eqref{averagedJextendednew} in
the phase space of $(J_\j,\theta_\j)$, has a saddle equilibrium  point

\begin{equation}\label{eq:tildeJ}
(\tilde J_\j (\hat J, \ep),\tilde \theta_\j (\hat J, \ep))= (J_\j ^*(\hat J) ,\theta_\j ^*(\hat J))+ O(\ep).
\end{equation}
for any $\hat J \in \hat \J$.

To make the pendulum-like structure of the  system given by
Hamiltonian~\eqref{averagedJextendednew} more apparent and to analyze the
behavior, we will find it convenient to make the translation
\begin{equation}\label{gap:canvi2}
y=J_\j-\tilde J_\j(\hat J,\ep), \qquad x= \theta_\j-\tilde \theta
_{\j}(\hat J,\ep),\quad s=s,
\end{equation}
and we obtain the $\C^{r-2m-2}$ Hamiltonian
\begin{equation}\label{gap:hamiltonian}
\K_0(y,x; \hat J,\ep)=  h_{0}(y;\hat J,\ep)+\ep^{j} U(x;\hat J,\ep)
\end{equation}
where
\begin{equation}\label{h0US}
\begin{split}
h_{0}(y;\ep\hat J,)&= l_0 y+h\left(N^{\top}(\hat J,\tilde J_{\j}+y)\right)\\
&-
h\left(N^{\top}(\hat J,\tilde J_{\j})\right)\\
&+ \ep \bar K^{0,0}\left(N^{\top}(\hat J,\tilde J_{\j}+y);\ep\right)\\
&- \ep \bar K^{0,0}\left(N^{\top}(\hat J,\tilde J_{\j});\ep\right) \\
U(x; \hat J,\ep)  &=
U^{k_0,l_0,*}( \tilde \theta_{\j}+x;\hat J,\ep)-
U^{k_0,l_0,*}(\tilde \theta_{\j};\hat J, \ep)
\end{split}
\end{equation}
where we have subtracted a constant term to the averaged
Hamiltonian~\eqref{averagedJextendednew}, the energy of the saddle $(\tilde
J_\j ,\tilde \theta_\j)=(\tilde J_\j (\hat J, \ep),\tilde \theta_\j (\hat J,
\ep))$, to normalize:
\begin{equation} \label{sella}
\begin{split}
&h_{0}(0; \hat J,\ep) = \frac{\partial h_{0}}{\partial
y}(0;\hat J,\ep)=0, \\
&\frac{\partial ^2 h_{0}}{\partial
y^2}(0;\hat J,\ep)= a(\hat J, \ep)=\a  + O(\ep)\ne 0, \\
&U(0;\hat J,\ep)= \frac{\partial U}{\partial
x}(0;\hat J,\ep)=0, \ \frac{\partial^{2} U}{\partial
x^{2}}(0;\hat J,\ep) \le \beta <0,
\end{split}\end{equation}
therefore the averaged Hamiltonian \eqref{gap:hamiltonian} can be written as:
\begin{equation}\label{gap:hamiltonian1}
\K_0(y,x;\hat J,\ep) = a(\hat J, \ep) \frac{y^2}{2} (1+O(y))+ \ep ^j U(x;\hat J,\ep)
\end{equation}
and  $(0,0)$  is a saddle point of Hamiltonian \eqref{gap:hamiltonian}, with
energy level $\K_0(0,0;\hat J,\ep)=0$.

Once we have the averaged Hamiltonian written in the form \eqref{averagedxy},
we discuss the geometry and the dynamics on the sets obtained by fixing the
energy level.

The main observation is that if we fix $\hat J=\hat c$, there is a critical
value $c_\j^*(\hat c, \ep)=0=\K_0(0,0;\hat J,\ep)$   for $c_\j$ at which the
topology and the dynamics of the level sets of the Hamiltonian $\K_0(y,x;\hat
J,\ep)$, and therefore of the Hamiltonian \eqref{averagedJextendednew},
change.

Now we describe the invariant sets of Hamiltonian $\K_{0}(y,x,\hat J;\ep)$
given in~\eqref{gap:hamiltonian} in the region:
 \begin{equation}\label{domainD}
D =\{ (\hat J, \hat \theta, y,x, s) \in \hat \J \times \torus  ^{d-1}\times \real
\times \torus ^2 , \quad \abs{y}\le
\bar L\}.
\end{equation}
for some $0<\bar L<L$, where $\hat \theta =(\theta _{1}, \dots ,\theta
_{\j-1},\theta _{\j +1},\dots ,\theta _{d})$.

Given any value $\hat J=\hat c$ we consider in the $(y,x)$ space, the level
set $\K_0(y,x;\hat J,\ep)=c_\j$:
\begin{itemize}
\item When $c_\j >  0$ but close enough to zero, the level set in the
    $(y,x)$ annulus is composed by two  non-contractible circles.
\item When $c_\j <   0$ but close enough to zero, the level  set in the
    $(y,x)$ annulus is a circle which, however, is contractible to a
    point.
\item When $c_\j = 0$, the level set is the union of two separatrices and
    the hyperbolic critical point $(0,0)$.
\end{itemize}

Therefore, the region $D$ is filled by the the level sets of the constants of
motion, that is, the energy surfaces of the Hamiltonian $\K_0$, and the
corresponding $\hat J$:
\begin{equation}\label{torusxy}
\T_{c}^{0} = \{ (\hat J, \hat \theta, y,x,s) \in \real ^{d-1} \times \torus  ^{d-1} \times \real
\times \torus ^2 : \, \K_{0}( y,x;\hat J, \ep)= c_{\j}, \hat J=\hat c \}.
\end{equation}
$\T_{c}^{0}$ will, of course, be invariant by the Hamiltonian flow of $\K
_0$.

The sets $\T_c^0$ consist on:
\begin{itemize}
\item When $c_\j >  0$ but close enough to zero, the level set
    $\T_{c}^{0}$ is composed by two primary tori (non-contractible tori
    of dimension $d+1$).
\item When $c_\j <   0$ but close enough to zero, the level  set
    $\T_{c}^{0}$ is a secondary torus (torus of dimension $d+1$
    contractible  to a $d$- dimensional torus).
\item The level set corresponding to $c_\j =  0$ consists of  one
    whiskered torus and its coincident whiskers: the hyperbolic torus
    $\torus^{d }\times \{(0,0)\}$ and the homoclinic  orbits to it. We
    will refer to $\T^0_{0}$ as the separatrix loop.
\end{itemize}

Formula  \eqref{torusxy} gives an implicit equation for the tori
$\T_{c}^{0}$. To compute the images of these tori under the scattering map in
section  \ref{sec:scattering}, it will be convenient to have the explicit
equation of these tori.

\begin{lemma}\label{toriexplicit}
There exists $\rho>0$, such that the two primary tori (components of the
secondary tori) $\T_{c}^{0}$ of Hamiltonian $\K_0$ can be written as graphs
of the variables $(\hat J, y)$ over the angle variables $(\hat \theta, x)$,
for $\rho\le x \le 2\pi -\rho$:
\begin{equation}\label{torusxyexplicit}
\begin{aligned}
\T_{c}^{0} =& \{ (\hat J, \hat \theta, y,x,s) \in \real ^{d-1} \times \torus  ^{d-1} \times \real
\times [\rho,2\pi-\rho] \times \torus  :\\
&\quad \hat J=\hat c, \  y= \pm \tilde \Y(x;c,\ep)\}.
\end{aligned}
\end{equation}
where the function $\tilde \Y(x;c,\ep)$ has different expressions depending
of the value $c_\j$:
\begin{enumerate}
\item If $0<c_\j\le\ep ^j$:
\begin{equation}
\tilde \Y(x;c,\ep)= \tilde \ell(x;c,\ep) (1+O(\ep ^\frac{j}{2}))
\end{equation}
\item If $\ep ^j\le c_\j<1$, for $c_\j = d_\j \ep ^\alpha$, with
    $0<\alpha <j$:
\begin{equation}
\tilde \Y(x;c,\ep)= \tilde \ell(x;c,\ep)(1 +O(\ep ^\frac{\alpha}{2}))
\end{equation}
\item If $c_\j=O(1)$:
\begin{equation}
\tilde \Y(x;c,\ep)= h_0^{(-1)}(c_\j)(1+O(\ep ^j))
\end{equation}
\end{enumerate}
where the function where $h_0$ is given in \eqref{sella} and $\tilde \ell$ is
given by:
$$
\tilde \ell(x;c,\ep)= \sqrt{\frac{2}{a(\hat c, \ep)}(c_\j-\ep ^j U(x;\hat c,\ep))}
$$
\end{lemma}
Once we know the structure of the level sets of the averaged Hamiltonian
\eqref{averagedJextendednew} in terms of the variables $(\hat J,\hat
\theta,y,x)$, we can write the equations of these sets in the original
variables of the problem.

First, in terms of the variables
 $(\B, \alpha ,s)\in \real ^{d} \times \torus  ^{d+1} $, using \eqref{changeJ}, \eqref{eq:projection},
 \eqref{gap:canvi2} and \eqref{gap:hamiltonian} and \eqref{h0US}, the tori $\T_{E}^{0} $, with $E=(\hat E, \tilde E_\j)$,
 are given by:
\begin{equation}\label{eq:toriaveraged}
\begin{split}
\hat \B - \frac {\B_{\j}} {k_0^\j}\hat { k_0 } &= \hat E \\
\frac{\B_\j}{k_0^\j}l_0+h(\B)  + \ep K^{0,0}(\B;\ep)+ \ep ^{j} U^{k_0,l_0,*}( k_0 \alpha +l_0 s;\hat E,\ep)  &= \tilde E_\j .
\end{split}
\end{equation}
$\T_{E}^{0}$ will, of course, be invariant by the Hamiltonian flow of the
averaged system \eqref{eq:averagedresonant}.

\begin{remark}
The equations \eqref{eq:toriaveraged} are a natural consequence of the fact
that the non-autonomous Hamiltonian \eqref{eq:averagedresonant} has as first
integrals the functions that are at the left hand side of
\eqref{eq:toriaveraged}.
\end{remark}
Let us observe that the relation between the constants $c$ and $E$ is given
by
\begin{eqnarray*}
\hat E&=&\hat c\\
\tilde E_\j &=& c_\j+\tilde E_\j ^*(\hat E)=c_\j +  l_0 \tilde J_\j  +h(N^{\top}( \hat E, \tilde J_\j ) )\\
&+& \ep  \bar K^{0,0} (N^{\top}(\hat E, \tilde J_\j ))
+ \ep ^{j} U^{k_0,l_0,*}(\tilde \theta _\j;\hat E, \ep)
\end{eqnarray*}
and the critical value where the topology of the invariant tori change is now
$\tilde E_\j ^*=\tilde E_\j ^*(\hat E)$,
which is the energy level of the saddle $(\tilde J_\j,\tilde \theta
_\j)=(\tilde J_\j (\hat E, \ep),\tilde \theta _\j(\hat E, \ep))$ and
corresponds to taking the critical value $c_\j=0$, $\hat c= \hat E$, in
\eqref{torusxy}, of Hamiltonian $\K_0$ in \eqref{gap:hamiltonian} with
variables \eqref{gap:canvi2}.

It is important to note that  equations \eqref{eq:toriaveraged} can also be
written, using \eqref{changeJ},  \eqref{gap:canvi2}  and
\eqref{eq:projection}, as:
\begin{equation}\label{eq:toriaveragedy}
\begin{aligned}
\hat \B - \frac {\B_{\j}} {k_0^\j}\hat { k_0 }&= \hat E \\
a(\hat E,\ep)\frac{y^2}{2}(1+O(y))  + \ep ^{j} U^{k_0,l_0,*}( k_0 \alpha +l_0 s;\hat E,\ep) &= E_\j
\end{aligned}
\end{equation}
where
$$
y=\frac{\B_\j}{k_0^{\j}}- \tilde J_\j (\hat E,\ep)=\frac{\B_\j- \B_\j^*(\hat E)}{k_0^{\j}}+O(\ep).
$$
The value of the critical value $E_\j^*$ where the topology of the tori
changes is given by
$$
E^*_\j= \ep ^j U^{k_0,l_0,*}(\tilde \theta _m(\hat E,\ep);\hat E,\ep),
$$
which is just the value of the potential at the saddle point. Again, using
the changes \eqref{changeJ}, \eqref{gap:canvi2} and formula
\eqref{torusxyexplicit} of Lemma \ref{toriexplicit}, we can obtain explicit
formulae for these tori:
\begin{eqnarray*}
\hat \B&=& \hat E +\frac{\B_\j}{k_0^{\j}}\hat k_0 \\
\B_\j &=& k_0^\j \tilde J_\j (\hat E,\ep) \pm  k_0^\j \Y(k_0 \alpha +l_0 s; E,\ep)\\
 &=& k_0^\j J_\j ^*(\hat E) \pm  k_0^\j \Y(k_0 \alpha +l_0 s; E,\ep) + O(\ep) \\
 &=& \B_\j ^*(\hat E) \pm  k_0^\j \Y(k_0 \alpha +l_0 s; E,\ep) + O(\ep)
\end{eqnarray*}
where
$$
\Y(\theta_\j ; E,\ep)= \tilde \Y(\theta_\j -\tilde \theta_\j (\hat E,\ep); E;\ep).
$$

Going back to the original variables $ (I, \varphi, s) \in \real ^{d} \times
\torus  ^{d+1} $ we can write the implicit equations for these tori
$\T_{E}^{0}$ as:
\begin{eqnarray*}
\hat I - \frac {I_{\j}} {k_0^\j}\hat { k_0 } + O(\ep)&=&\hat E,\\
\frac{I_\j}{k_0^\j}l_0 +h(I)  +  \ep K^{0,0}(I;\ep)+ \ep ^{j}U^{k_0,l_0,*}( k_0 \varphi +l_0 s;\hat E,\ep)  + O(\ep ) &=& \tilde E_\j
\end{eqnarray*}

That can also be  also written as:
\begin{eqnarray}
\hat I - \frac {I_{\j}} {k_0^\j}\hat { k_0 } + O(\ep)&=&\hat E, \label{eq:toriaveragedI} \\
a(\hat E,\ep)\frac{y^2}{2}(1+O(y))  + \ep ^{j} U^{k_0,l_0,*}( k_0 \varphi +l_0 s;\hat E,\ep) +O(\ep^{j+1}) &=&E_\j  \nonumber
\end{eqnarray}
where:
\begin{equation}\label{eq:toriaveraged2}
y=\frac{I_\j- \B_\j^*(\hat E)}{k_0^{\j}}+O(\ep), \quad E^*_\j= \ep ^jU^{k_0,l_0,*}(\tilde \theta _m (\hat E,\ep);\hat E,\ep).
\end{equation}

Finally, these tori can be also written explicitly as:
\begin{eqnarray}
\hat I &=&\hat E +\frac{I_\j}{k_0^\j}\hat k_0 +O(\ep)  \label{eq:toriexplicitI} \\
I_\j
&=& \B_\j ^*(\hat E) \pm  k_0^\j \Y(k_0 \varphi +l_0 s; E,\ep) + O(\ep)\nonumber
\end{eqnarray}

Let us observe that the function $\Y$ verifies, according to Lemma
\ref{toriexplicit}:
\begin{itemize}
\item If $0<| E_\j-  E_\j^*|\le \ep ^j$:
\begin{equation}
\Y(\theta_\j ; E,\ep)= \ell(\theta_\j ; E,\ep) (1+O(\ep ^\frac{j}{2}))
\end{equation}
\item If $\ep ^j<| E_\j-  E_\j^*|<1$, writing $| E_\j-  E_\j^*| = d_\j
    \ep ^\gamma$, with $0<\gamma <j$:
\begin{equation}
\Y(\theta_\j ; E,\ep)= \ell(\theta_\j ; E,\ep)(1 +O(\ep ^\frac{\gamma}{2}))
\end{equation}
\item If $| E_\j-  E_\j^*|=O(1)$:
\begin{equation}
\Y(\theta_\j ; E,\ep)= h_0^{(-1)}(E_\j)(1+O(\ep ^j))
\end{equation}
\end{itemize}
where the function $\ell$ is given by (see \eqref{h0US}:
\begin{equation}\label{eq:ell}
\ell(\theta_\j ; E,\ep)= \sqrt{\frac{2}{a(\hat E, \ep)}( E_\j-\ep ^j U^{k_0,l_0,*}(\theta_\j ;\hat E,\ep))}.
\end{equation}

\subsubsection{Primary and secondary tori near the
secular resonances: KAM Theorem}\label{sec:toriKAM}

If we apply the changes of variables \eqref{eq:changesimplectic},
\eqref{changeJ} and \eqref{gap:canvi2} to Hamiltonian \eqref{generalterm} we
obtain:
\begin{equation} \label{hamiltonianaveragedresonance}
\K(\hat J,\hat \theta, y,x,s;\ep) =\K_{0}(\hat J, y,x;\ep) + \ep ^{m+1} S(\hat J,\hat \theta, y,x,s;\ep).
\end{equation}

First, we change to action--angle variables of the integrable part $\K_{0}$.
The two only difficulties are that the action angle variables become singular
near the level set $\K_{0}=0$, which is the separatrix  of the torus $\{
(x,y)=(0,0)\}$, and also the twist condition becomes singular. In fact, the
twist goes to $\infty$ when one approaches the separatrix and this is
favorable to application of the KAM theorem \cite{Herman85, DelshamsLS06a})
because this theorem only requires lower bounds on the twist and a larger
twist improves the quantitative assumptions of the theorem.

If the number of averaging steps  $m$ is large enough we can ensure that
there exist KAM tori (both primary and secondary) that cover the whole
resonant region up to distances $O(\ep^{3/2})$ and which are $\ep^3$ close to
the level sets $\K_0=c_\j$ of the averaged Hamiltonian $\K_0$. Of course, we
could get even higher powers in the density by averaging more times.

\begin{itemize}
\item We select a region $| c_{\j} - c^{*}_{\j}| \le  \ep ^{\alpha}$
    surrounding the separatrix $\K_0=c_\j ^*=0$. In this region (the
    chaotic zone) we will not perform any further analysis. We will just
    remark that it is small. In particular, if $\alpha
    =\frac{3}{2}+\frac{j}{2}$, the level sets of energy $ c_\j$ and $
    c^{*}_{\j}$ are at a distance $\ep^{3/2}$ if $( c_{\j} - c^{*}_{\j})
    = \pm   \ep ^{\alpha}$ .
\item In the complementary region: $| c_{\j} - c^{*}_{\j}| \ge  \ep
    ^{\alpha}$ we  change to action angle variables adapted to the level
    sets, $F=c=(\hat c,c_\j)$, of the function $F=(\hat J,\K _0)$ defined
    in \eqref{torusxy}.  We note that one of the components of $F$ is
    precisely $\K_0$, the integrable part of the averaged Hamiltonian
    \eqref{hamiltonianaveragedresonance}.  The action variables can be
    obtained geometrically integrating the canonical form over the loops
    in a torus \cite{Arnold78, AbrahamM78}.

It is well known  that the singularities of the action variable are only
 a power of $(c_{\j} - c^{*}_{\j})$. Therefore, since the size of the
 remainder in Hamiltonian \eqref{hamiltonianaveragedresonance}
is $O(\ep^{m+1})$, when expressed in the action angle variables in the
region $|c_{\j} - c^{*}_{\j}| \ge \ep ^{\alpha}$ the smallness of the
remainder  will be $O(\ep ^{m+1 - A \alpha})$, for some value $A>0$.
\item The KAM theorem in action-angle variables \cite{Poschel01} gives
    tori which are at a distance $\ep ^{m+1 - A \alpha}$ of the level
 sets of the action variables. The gaps between these tori are
$O(\ep^{(m+1 - A \alpha)/2})$.
\item Coming back to variables $(y,x,\hat J, \hat \theta)$ we obtain tori
    at a distance between them of order $O(\ep^{(m+1 - A \alpha)/2
    -\alpha A })$, and that are at a distance $O(\ep^{m+1 - 2 A\alpha}) $
    from the level sets of $(\hat J,\K_0)$.

We note that, if we fix $\alpha = \frac{3}{2}+\frac{j}{2}$, taking into
account that $A$ is a fixed number, we obtain that, taking $m$ large
enough we can ensure that  we have tori for $|c-c'|< \ep
^{\frac{3}{2}+\frac{j}{2}}$, that is, the gaps are smaller than
$\ep^{3/2}$ as claimed.
\item Going back to the original variables through changes
    \eqref{eq:projection}, \eqref{gap:canvi2} (which are close to the
    identity),  we obtain the result in next theorem \ref{innermap0}.
\end{itemize}

\begin{theorem}\label{innermap0}
Under the conditions of theorem \ref{main}, there exists $\varepsilon_0$ such
that, for $0<|\varepsilon|<\varepsilon_0$, the secular resonant region
$\S^{[\le 2]}$ can be covered by $O(\ep^{3/2})$ neighborhoods of invariant
objects under the Hamiltonian flow of the Hamiltonian $K_{\ep}(I,\vp,s)$ in
\eqref{eq:hamiltonianrestricted}. Moreover:
\begin{itemize}
\item These invariant objects are given by the level sets $F=E=(\hat
    E,E_\j)$, for $|E-E'| \le \ep ^{ \frac{3}{2}+\frac{j}{2}}$ and where
    $\ep ^{   \frac{3}{2}+\frac{j}{2} } \le | E_\j- E_\j ^*|\le 1$ with
    $E _\j^*$ given in \eqref{eq:toriaveraged2}.
\item The $\C^{2}$ function $F: \real^d\times \torus^d \times \torus \to
    \real$ is given
 by \eqref{eq:toriaveragedI}.
\item These invariant objects are either regular primary KAM $d+1$-tori,
    secondary $d+1$-KAM tori of class $1$ (i.e. $d+1$-dimensional
    invariant tori which are contractible to a $\torus^{d}$) or invariant
    manifolds of $d$-dimensional whiskered invariant tori.
\end{itemize}
\end{theorem}

\subsection{Second step: The generation of
a homoclinic manifold and computation of the scattering map
}\label{sec:scattering0}

Let us observe that by hypothesis {\bf H2}, for $\ep=0$, the manifold
$\tilde\Lambda _{0}$ has stable and unstable manifolds which coincide along a
homoclinic manifold
$$
\tilde \Gamma_0=W^s(\tilde\Lambda _{0})=W^u (\tilde\Lambda
_{0})
$$
with
$$
\tilde \Gamma_0=\{ (p^*(\tau), q^*(\tau), I, \varphi,s) , \ (I,\vp, s,\tau) \in
\II^* \times \torus ^{d +1}\times \real ^n\}
$$
 The first result in this section is that, if system
\eqref{eq:hamiltonianequations} satisfies the non-degeneracy assumption {\bf
H7}, then for all $0 < |\ep |< \ep_0$, $W^s(\tilde\Lambda _{\ep})$, $W^u
(\tilde\Lambda _{\ep})$, the stable and unstable manifolds of the normally
hyperbolic invariant manifold $\tilde \Lambda _\ep$ introduced in Section
\ref{normalhyp}, have a transversal intersection along a homoclinic manifold
$\tilde \Gamma_{\ep}$. Then, following \cite{DelshamsLS08} we use this
intersection to define the scattering map in $H_-\subset
\tilde\Lambda_{\ep}$.

We will use a notation very similar to that of \cite{DelshamsLS06a} and,
indeed refer to this paper for a series of detailed calculations. The proof
of the next proposition is identical to proposition 9.2 in
\cite{DelshamsLS06a}.

\begin{prop}\label{prop:transverseintersection}
Assume that hypothesis \textbf{H7} is fulfilled. Then, given
$(I,\varphi,s)\in H_- \subset \II^*  \times \torus ^{d+1}$, for $\ep$ small
enough, there exists a locally  unique point $\tilde z^*$ of the form
\begin{equation}
\label{eq:intersectionform}
 \tilde z^*(I,\vp, s; \ep) =
(p^{*}(\tau^{*}(I,\vp,s))+\Or(\ep),q^{*}(\tau^{*}(I,\vp,s))+\Or(\ep),I,\vp ,s)
\end{equation}
such that
 $W^{\st}(\tilde \Lambda _{\ep})\itr
 W^{\un}(\tilde \Lambda _{\ep})$ at $\tilde z^*$,
that is,
 \[\tilde z ^*\in W^{\st}(\tilde \Lambda _{\ep})\cap W^{\un}(\tilde
\Lambda _{\ep}) \mbox{ and }
 T_{\tilde z^*} W^{\st}(\tilde  \Lambda _{\ep})+
T_{\tilde z^*} W^{\un}(\tilde
 \Lambda _{\ep})= T_{\tilde z^*} \tilde \M ,
 \]
where $\tilde \M = \real ^{n}\times \real ^{n}\times  \II \times\torus ^{d}
\times \torus. $

In particular, there exist unique points
\[
 \tilde x_{\pm}=\tilde x_{\pm}(I,\vp,s;\ep)=(0,0,I,\vp,s)
 +\Or _{\C ^{1}}(\ep)\in \tilde \Lambda _{\ep}
\]
such that
\begin{equation} \label{eq:asymptotic}
\abs{\tilde\Phi_{\ep,t}(\tilde z^*)-\tilde\Phi_{\ep,t}(\tilde
x_{\pm})} \le \cte e^{-\alpha \abs{t}/2} \ \mbox{ for }\  t \to \pm
\infty.
\end{equation}
Moreover, expressing the points $\tilde x_{\pm}=k_{\ep} (I_{\pm},\vp
_{\pm},s_{\pm};\ep)$ in terms of the symplectic parametrization of $\tilde
\Lambda_{\ep}$ introduced in section \ref{normalhyp}, the following formulas
hold:
\[
I(\tilde x_{\pm}) = I +\Or _{\C ^{1}}(\ep), \qquad \vp (\tilde
x_{\pm}) = \vp +\Or _{\C ^{1}}(\ep), \qquad s(\tilde x_{\pm}) = s,
\]
and
\begin{equation}\label{eq:scattering}
 I(\tilde x_{+}) -I(\tilde x_{-}) = \ep  \frac{\partial  \LL^{*}}{\partial \vp}
 (I,\vp-\omega(I)s)+\Or _{\C ^{1}}(\ep^{1+\varrho}),
\end{equation}
where $\LL^{*}(I,\theta)$ is given in \eqref{generator}, and $\varrho>0$.
\end{prop}

\subsection{The scattering map }
{From} now on we will take the homoclinic manifold

$$
\tilde\Gamma_{\ep}=\{  \tilde z^*(I,\vp, s; \ep), \ (I,\vp, s) \in
H_- \subset \II^* \times \torus ^{d+1}\}
$$
given by  proposition \ref{prop:transverseintersection}. Following
\cite{DelshamsLS08} we will call $\tilde\Gamma_{\varepsilon}$  an homoclinic
channel  and we can define the scattering map $s_{\ep}$ on
$\tilde\Lambda_{\varepsilon}$, also called outer map, associated to
$\tilde\Gamma_{\ep}$.

Following \cite{DelshamsLS06a,DelshamsLS08} the scattering map is defined as
follows: for any two points $\tilde x_{\pm} \in \tilde \Lambda_\ep$, we say
that $\tilde x_{+}=s_\ep (\tilde x_{-})$, if there exists a point $\tilde z
\in \tilde\Gamma_{\varepsilon}$ such that
\[
\dist({\tilde \Phi _{\ep, t}(z),\tilde \Phi_{\ep, t}(x_{\pm})})\to 0, \quad
\mbox{for}\ t\to \pm \infty .
\]
Since the  unperturbed system is a product system, it is clear that,
independently of what is the homoclinic manifold, the stable manifold of one
point in $\tilde \Lambda_0$ is the same as its unstable manifold. Therefore,
$s_0=\mathrm{\Id}$.

As shown in \cite{DelshamsLS08} the scattering map is an exact symplectic map
and depends smoothly on parameters because the homoclinic manifold depends
smoothly on parameters even through $\ep =0$.

It is well known (see \cite{LlaveMM86}) that a family of exact symplectic
mappings $s_\ep$ is conveniently described using a generator
$\mathcal{S}_\ep$ and the associated Hamiltonian $S_\ep$:
\[
\frac{d}{d \ep} s_\ep = \mathcal{S}_\ep \circ s_\ep;
\quad \quad  \imath_{\mathcal{S}_\ep} \omega = d S_\ep .
\]
Indeed in \cite{DelshamsLS08} it is shown that the Hamiltonian $S_\ep$ is
given,  up to first order in $\ep$ by the function
$-\L^*(I,\varphi-\omega(I)s)$, where the reduced Melnikov potential
$\L^*(I,\theta)$ is given in \eqref{generator}:
$$
S_\ep =  -\L^* + \ep S_1 + O(\ep ^2)
$$
Therefore, the scattering map can be written, using the coordinates
$(I,\varphi,s)$ as:
\begin{equation}\label{scatteringformula}
\begin{array}{rcl}
s_\varepsilon(I,\varphi,s)&=&
(I+\ep \partial _{\theta}\L^{*}(I,\varphi-\omega(I)s)+O(\ep^2), \\
&&\varphi-\ep \partial _I\L^{*}(I,\varphi-\omega(I)s) + O(\ep^2), s)
\end{array}
\end{equation}
and, for any fixed $s\in \torus$, up the first order in $\ep$ it is given, in
the coordinates $(I,\varphi)$, as the time $-\ep$ map of the Hamiltonian flow
of Hamiltonian $\L^*(I,\theta)$ evaluated at $(I,\varphi-\omega(I)s)$.

The fundamental property to have instability will be to check, for any fixed
$s$,  that the tori invariant for the inner flow in $\tilde \Lambda _\ep$ are
not invariant by the perturbed scattering map $s_\ep$. Therefore, we will pay
attention at how the scattering map moves the tori $\T_{E}$ given in
\eqref{eq:primarytorus1}, \eqref{eq:toriaveragedI},
 \eqref{eq:toriaveraged2}
using the results in \cite{DelshamsH09}.

\subsection{Interaction between the inner flow and the
scattering map and hypothesis {\bf H6}} \label{sec:scattering}

We have already shown  in theorems \ref{KAMnoresonant} and \ref{innermap0},
that the KAM tori $\T _E$ (both primary and secondary) are the level sets of
an $\real^d$-valued function $F_\ep$. Indeed we have approximate expressions
for it in \eqref{eq:primarytorus1}, \eqref{eq:toriaveragedI},
 \eqref{eq:toriaveraged2}  (see also \eqref{eq:toriexplicitI}).

The scattering map transports the level sets of $F_\ep $ into other
manifolds, which are the level sets of $F_\ep \circ s_\ep^{-1}$.

The key observation relies on Lemma 10.4 in \cite{DelshamsLS06a} (see also
\cite{DelshamsLS00}), which states that, given two invariant manifolds for
the inner flow  $\Sigma_i \subset \Lambda$, $i=1,2$, if $\Sigma _1$
intersects transversally $s_\ep(\Sigma_2)$ in $\tilde \Lambda _\ep$, then
$W^u_{\Sigma_2} \pitchfork W^s_{\Sigma_1} $.

Our next goal will be to make explicit the conditions to ensure that the
scattering map creates heteroclinic intersections between the KAM tori,
primary or secondary, created in sections \ref{sec:non-resonantKAM} and
\ref{sec:toriKAM}.

Fix $s\in \torus$. The tori $\T_{E'}$ and $s_\ep (\T_E)$ intersect if there
exists a point $\tilde x=(I,\varphi,s)$ such that:
\begin{equation}\label{intersectiontori}
\begin{array}{rcl}
F_{\ep}(I,\varphi,s;\ep)=E \\
F_{\ep}\circ s_{\ep}(I,\varphi,s;\ep)=E'.
\end{array}
\end{equation}
Let us observe that the first equation is the implicit equation for the torus
$\T_E$. Instead, we can use its explicit equation $I=\lambda
_E(\varphi,s;\ep)$ to eliminate the first $d$ equations. The manifolds
$\T_{E'}$ and $s_\ep (\T_E)$ intersect if there exists $(\varphi,s)$ such
that:
\begin{equation}\label{eq:intersection}
F_{\ep}\circ s_{\ep}(\lambda _E(\varphi,s;\ep),\varphi,s;\ep)=E'
\end{equation}
and the intersection will be transversal if
\begin{equation}\label{transversal1}
\det D(F_\ep \circ s_\ep (\lambda _E(\varphi,s;\ep),\varphi,s;\ep))\ne 0 ,
\end{equation}
where $D=D_{\varphi}$.

Using  formula (159) in \cite{DelshamsH09}, we know that, given a function
$F$:
\begin{equation}\label{scateringDH}
F\circ s_{\ep}= F - \ep \{F,\L^{*}\}+ \frac{\ep^2}{2}(\{\{F,\L^*\},\L^*\}+\{F,S_1\}) +O(\ep ^3).
\end{equation}
Therefore equation \eqref{eq:intersection} reads:
\begin{equation}\label{eq:heterocliniques}
 - \{F_\ep,\L^{*}\}(\lambda _E(\varphi,s;\ep),\varphi,s;\ep)+ O(\ep)=\frac{E'-E}{\ep}
\end{equation}
and we will have intersection as long as  $\frac{E'-E}{\ep} $ is small
enough, close to the non-degenerate zeros of
$$
\{F_\ep,\L^{*}\}(\lambda _E(\varphi,s;\ep),\varphi,s;\ep)=0
$$
therefore the transversality condition \eqref{transversal1}  is equivalent to
\begin{equation}\label{transversal11}
D(\{ F_\ep, \L^{*}\})(\lambda _E(\varphi,s;\ep),\varphi,s;\ep)\ne 0
\end{equation}
for a point in each of the level sets of $F_\ep$.

\subsubsection{The non-degeneracy condition {\bf H8} in the non-resonant region}

The non-resonant region $\S^L$ (see \eqref{non-resonantregion}) is of $O(1)$
and is covered by $\ep^{\frac{3}{2}}$ neighborhoods of tori which are given
by the level sets of the function:
$$
F_\ep (I, \varphi ,s)= I+ O(\ep) =E.
$$
Therefore $\lambda _E(\varphi,s;\ep)=E+O(\ep )$, and equation
\eqref{eq:heterocliniques} reads:
\begin{equation}\label{heteroclinicsclosetoh}
\frac{\partial \L^*}{\partial \theta}(E, \varphi-\omega (E)s) +O(\ep) =\frac{E'-E}{\ep}
\end{equation}
Moreover, by the KAM theorem \ref{KAMnoresonant} given in section
\ref{sec:non-resonantKAM} we have tori for $|E-E'|\le c \ep ^{3/2}$, and
therefore  equation \eqref{heteroclinicsclosetoh} has solutions for
$\varphi$, which are non-degenerate if condition \eqref{transversal11} is
verified, which in our case, becomes:
$$
{\rm det} \left | \frac{\partial ^{2}  \L^{*} }{\partial \theta ^{2}} (E, \varphi-\omega (E)s)\right |\ne 0,
$$
and  is guaranteed if
\begin{equation}\label{eq:nondegeneratenores}
{\rm det} \left |\frac{\partial ^{2}  \L^{*} }{\partial \theta ^{2}} (I, \varphi-\omega (I)s)\right | \ne 0
\end{equation}
is satisfied for $(I,\varphi,s)  \in H_- \subset \II^*\times \torus ^{d+1}$,
and is one of the non-degeneracy conditions included in  Hypothesis {\bf H8}.

\subsubsection{Heteroclinic orbits close to homoclinic ones  in the non-resonant region}

If condition \eqref{eq:nondegeneratenores} is verified in the region  $\S^L$
we can guarantee the existence of heteroclinic connections between
neighboring KAM tori in this region. If we look for  heteroclinic connections
close to homoclinic ones, one can obtain  a more explicit sufficient
condition for equations \eqref{heteroclinicsclosetoh} to have a solution. The
main idea is to solve equations  \eqref{heteroclinicsclosetoh}, using the
implicit function theorem. The small parameter will be
$$
\delta=\frac{E'-E}{\ep}+\varepsilon
$$
and then equation \eqref{heteroclinicsclosetoh} read:
$$
\frac{\partial \L^*}{\partial \theta}(E, \varphi-\omega (E)s)  =O(\delta)
$$
Therefore, a non-degeneracy condition which guarantees that equation
\eqref{heteroclinicsclosetoh} have solutions close to the solutions of:
$$
\frac{\partial \L^*}{\partial \theta}(E, \varphi-\omega (E)s)  =0
$$
is that the function $\L^*$ has non-degenerate critical points, that is:

\begin{equation}\label{eq:h8non-resonant}
\frac{\partial \L^*}{\partial \theta}(E, \varphi-\omega (E)s)  =0, \implies \quad
{\rm det}\left |\frac{\partial ^{2}  \L^{*} }{\partial \theta ^{2}} (E, \varphi-\omega (E)s)\right |\ne 0,
\end{equation}
in the region $\S^L\cap H_-$. Equation \eqref{eq:h8non-resonant} is part of
Hypothesis {\bf H8}.

\subsubsection{The non-degeneracy condition {\bf H8} in the resonant region}

Now, we study  the intersection equation \eqref{intersectiontori} in the
secular resonant region \eqref{secularregion} $\S^{[\le 2]}\cap H_-$, to
ensure that the image under the scattering map of a primary or secondary
torus intersects other  nearby tori. We  denote  by $F_\ep$ again the
function whose level sets gives the tori. We recall that the secular resonant
region $\S^{[\le 2]}$ is the union of the tubular neighborhoods  $\R^L_{k,l}$
of the secular resonances $\R_{k,l}$, for $(k,l) \in \N^{[\le 2]}$.

If $\R_{k_0,l_0}$ is a resonance of order $j$, $j=1,2$, in the  region
$\R^L_{k_0,l_0}$, according to \eqref{eq:toriaveragedI},
\eqref{eq:toriaveraged2} and the KAM theorem \ref{innermap0}, the invariant
tori are given by the level sets of a function:
$$
F_\ep = (\hat F, F_\j) =E =(\hat E, E_\j),
$$
for $\ep ^{\frac{3}{2}+\frac{j}{2}} \le |E_\j -E^*_\j| \le 1$, where $E^*_\j
= \ep ^j U^{k_0,l_0,*}(\tilde \theta _\j(\hat E,\ep);\hat E,\ep)$ (see
\cite{DelshamsLS06a}), with $U^{k_0,l_0,*}$  given in \eqref{Uestrella},
\begin{equation}\label{above1}
\hat F(I, \varphi ,s;\ep)= \hat I-\frac{I_\j}{k_0^\j}\hat k_0 +O(\ep),
\end{equation}
and
\begin{equation}\label{above2}
F_\j (I, \varphi ,s;\ep)
= a(\hat E,\ep)\frac{y^2}{2}(1+O(y))  +
\ep ^{j} U^{k_0,l_0,*}( k_0 \varphi +l_0 s;\hat E,\ep) +O(\ep^{j+1}) \\
\end{equation}
with $a(\hat E,\ep)$ given in \eqref{sella}, \eqref{coefquasiconvexity},  and
\[
y=\frac{I_\j- \B_\j^*(\hat E)}{k_0^{\j}}+O(\ep)
\]
where $\B_\j ^*(\hat E) =\Gamma _{k_0,l_0}(I)$ is the $k_0$-projection in the
resonance $\R_{k_0,l_0}$.

Moreover, by the KAM theorem  \ref{innermap0}, we know that there exist tori
$F_\ep=E$, $F_\ep=E'$ for $| E- E'| \le   \ep ^{\frac{3}{2}+\frac{j}{2}}$.

The way to solve equation \eqref{intersectiontori} is slightly different for
a resonance or order one or for a resonance of order  two. We give all the
details in the case of a first order resonance. The case of a resonance of
order two can be done with minor modifications, as it is explained in Remark
\ref{rem:secondorder}.

For  $F_\ep =(\hat F, F_\j)$, the set of equations \eqref{intersectiontori}
reads
\begin{equation*}
\begin{split}
\hat F(I, \varphi ,s;\ep)&= \hat E \\
F_\j (I, \varphi ,s;\ep)&=E_\j  \\
\hat F (I, \varphi ,s;\ep)- \ep\{ \hat F,\L^*\}(I, \varphi ,s;\ep)+ O(\ep ^2)&=\hat E' \\
F_\j (I, \varphi ,s;\ep)-  \ep \{ F_\j,\L^*\}(I, \varphi ,s;\ep)+ O(\ep ^2)&=E'_\j
\end{split}
\end{equation*}
and, using  \eqref{above1}, \eqref{above2}, these equations  are equivalent
to
\begin{equation}\label{intersectiona}
\begin{array}{@{}rcl@{}}
\hat F(I, \varphi ,s;\ep)&=& \hat E \\
F_\j (I, \varphi ,s;\ep)&=&E_\j \\
\multicolumn{3}{@{}l@{}}{\hskip -1cm\displaystyle\partial _{\hat \theta}\L ^* (I, \varphi-\omega(I)s) -
\frac{1}{k_0^\j}\partial _{\theta _\j}\L^* (I, \varphi-\omega(I)s)  \hat k_0}\\
+O(\ep )&=&\dfrac{\hat E'- \hat E }{\ep}\\
\multicolumn{3}{@{}l@{}}{\hskip -1cm\displaystyle a(\hat E,\ep)\frac{I_\j- \B_\j^*(\hat
E)}{(k_0^{\j})^2}\partial _{\theta
_\j} \L^* (I, \varphi-\omega(I)s}\\
+O(\ep)&=&\dfrac{E'_\j- E_\j}{\ep }.
\end{array}
\end{equation}
{From} the  first two equations we obtain, using  \eqref{eq:toriexplicitI}
for $E_\j -E_\j ^*=O(\ep ^\gamma)$, $0< \gamma \le
\frac{3}{2}+\frac{j}{2}=2$, we obtain:
\begin{eqnarray*}
I &=&(\hat E,0) +
\frac { I_\j } {k_0^\j} k_0+  O(\ep )\\
I_\j &=& \B^*_\j (\hat E)\pm k_0^\j \Y(k_0 \varphi+l_0 s;E;\ep) + O(\ep)
\end{eqnarray*}
where $E=(\hat E, E_\j)$ and
\begin{equation}\label{yele}
\Y(k_0 \varphi+l_0 s;E;\ep)= \ell (k_0 \varphi + l_0 s; E, \ep) (1+O(\ep ^{\frac{\gamma}{2}}))
\end{equation}
and the function $\ell$ is given in \eqref{eq:ell}:
$$
\ell(k_0 \varphi +l_0 s; E;\ep)= \sqrt{\frac{2}{a(\hat E, \ep)}( E_\j-\ep  U^{k_0,l_0,*}(k_0 \varphi +l_0 s;\hat E,\ep))}
$$
which gives, using that $E_\j -E_\j ^*=O(\ep ^\gamma)$:
\begin{eqnarray*}
I &=&(\hat E,0) + \frac{\B^*_\j(\hat E)} {k_0^\j} k_0= I^*(\hat E)
+O(\ep, \ep ^{\frac{\gamma}{2}})
\end{eqnarray*}
and therefore, the two last equations of \eqref{intersectiona} read:
\begin{equation}\label{hetero}
\begin{aligned}
\partial _{\hat \theta} \L ^* - \frac{1}{k_0^\j} \partial _{\theta _\j} \L^* \hat k_0
+ O(\ep , \ep ^{ \frac{\gamma}{2} })
&= -\frac{\hat E- \hat E'}{\ep}\\
\pm \frac{a(\hat E; \ep)} {k_{0}^{\j}} \left[ \Y(k_0 \varphi+l_0 s;E;\ep)+O(\ep)\right]\times&\\
\times \left(\partial _{\theta _\j}\L^*
+O(\ep , \ep ^{ \frac{\gamma}{2} })\right)&= -\frac{ \tilde E_\j-  E_\j'}{\ep}
\end{aligned}
\end{equation}
where, to shorten the notation, we have just written
$$
\L^*= \L^*(I^*(\hat E),\varphi-\omega (I^*(\hat E))s).
$$
Using equation~\eqref{yele} and that, by \eqref{sella},  $a(\hat
E,\ep)=a(\hat E) +O(\ep)$, equations \eqref{hetero} are equivalent to:
\begin{equation}\label{intersectionb}
\begin{aligned}
\partial _{\hat \theta} \L ^* - \frac{1}{k_0^\j} \partial _{\theta
_\j}
\L^*\hat k_0   +O(\ep , \ep ^{ \frac{\gamma}{2} } )
&= \frac{\hat E'- \hat E}{\ep} \\
\pm \frac{\left(a(\hat E)+O(\ep)\right)} {k_{0}^{\j}}
\left[ \ell(k_0 \varphi+l_0 s;E;\ep)   (1+O(\ep^{\frac{\gamma}{2}}))+O(\ep)\right]&\\
\times\left(\partial _{\theta _\j}\L^* +O(\ep, \ep ^{\frac{\gamma}{2}}))\right)
&= \frac{ \tilde E'_\j-  E_\j}{\ep}.
\end{aligned}
\end{equation}

We will see that we will have a solution of  equations \eqref{intersectionb}
for $\varphi$ if $| E'- E| \le O(\ep ^{ \frac{3}{2}+\frac{j}{2} })= O(\ep
^2)$.

It will be useful to work in the variables $\theta=(\hat \theta ,\theta
_{\j})=( \hat \varphi, k_0 \varphi +l_0 s)$. Observe that conversely:
\begin{equation}\label{eq:canvipoincare}
\hat \varphi = \hat \theta, \quad \varphi_\j = \frac{\theta _\j -\hat k_0 \hat \varphi -l_0 s}{k_0^\j}
 \end{equation}
Now, if we define the auxiliary function:
\begin{equation}\label{eq:poincareauxiliary}
\begin{split}
\L^*_{k_0,l_0}(\hat \theta, \theta _\j,s;\hat E)=& \L^* (I^*(\hat E),\varphi -\omega (I^*(\hat E))s) \\
=& \L^* (I^*(\hat E),\hat \varphi -\hat \omega (I^*(\hat E))s, \varphi
_\j-\omega _\j(I^*(\hat E))s)
\end{split}
\end{equation}
Using that $\omega (I^*(\hat E))\cdot k_0 + l_0=0$, we obtain:
\begin{multline*}
\L^*_{k_0,l_0}(\hat \theta, \theta _\j,s;\hat E)=\\
\L^* (I^*(\hat E),\hat \theta-\hat \omega(I^*(\hat E)) s, \frac{\theta_\j
-(\hat \theta -\hat \omega (I^*(\hat E)) s) \hat k_0}{k_0^\j}),
\end{multline*}
and then, taking derivatives with respect to $\hat \theta$ and $\theta_\j$:
\begin{eqnarray*}
\frac{\partial}{\partial \hat \theta}\L^*_{k_0,l_0} ( \hat \varphi, k
_0 \varphi +l_0 s,s;\hat E)&=&
\frac{\partial}{\partial \hat \theta}\L^* (I^*(\hat E),\varphi -\omega (I^*(\hat E))s)\\
&-& \frac{1}{k_0^\j}\frac{\partial}{\partial \theta _\j}\L^*(I^*(\hat E),\varphi -\omega (I^*(\hat E))s) \ \hat k_0\\
\frac{\partial}{\partial  \theta_\j}\L^*_{k_0,l_0}( \hat \varphi, k
_0 \varphi +l_0 s,s;\hat E) &=&  \frac{1}{k_0^\j}\frac{\partial}{\partial \theta _\j}\L^*(I^*(\hat E),\varphi -\omega (I^*(\hat E))s)
\end{eqnarray*}

Therefore equations \eqref{intersectionb} become:
\begin{equation}\label{hetero1}
\begin{aligned}
\frac{\partial}{\partial \hat \theta}\L^*_{k_0,l_0} (\hat \theta, \theta_\j,s;\hat E)
+ O(\ep ,\ep ^{\frac{\gamma}{2}}) &=\frac{\hat E'- \hat E}{\ep}\\
\pm \frac{(a(\hat E)+O(\ep))} {k_{0}^{\j}}
\left[\ell(\theta_\j;E;\ep) (1+O(\ep^{\frac{\gamma}{2}}))+O(\ep)\right]\times&\\
\times\left( \frac{\partial}{\partial  \theta_\j}\L^*_{k_0,l_0} (\hat \theta, \theta_\j,s;\hat E)
+O(\ep, \ep ^{\frac{\gamma}{2}})\right)&= \frac{ \tilde E'_\j-  E_\j}{\ep}
\end{aligned}
\end{equation}
which are the generalization to higher dimensions of the function $\M$ in
\cite[page 108]{DelshamsLS06a}.

Before looking for the solutions of these equations, me make a further
simplification. First  observe that there exist primary and secondary tori
close to the separatrix of the averaged Hamiltonian for energies
$$| E_\j - E_\j^* |> \ep ^{\frac{3}{2}+\frac{j}{2}}=O(\ep ^2), \quad
E_\j^*=\ep U^{k_0,l_0,*}(\tilde \theta (\hat E,\ep);\hat E,\ep)),
$$
therefore it makes sense to scale $E_\j = \ep  e_\j$ in the function $\ell$
of  equations \eqref{hetero1} obtaining:
\begin{equation}\label{heteroscaled}
\begin{aligned}
\frac{\partial}{\partial \hat \theta}\L^*_{k_0,l_0} (\hat \theta, \theta_\j,s;\hat E)
+ O(\ep ,\ep ^{\frac{\gamma}{2}}) &= \frac{\hat E'- \hat E}{\ep} \\
\pm \frac{(a(\hat E)+O(\ep))} {k_{0}^{\j}}
\left[\bar \ell(\theta_\j;E;\ep) (1+O(\ep^{\frac{\gamma}{2}})+O(\ep^{\frac{1}{2}}))\right]\times&\\
\times\left( \frac{\partial}{\partial  \theta_\j}\L^*_{k_0,l_0} (\hat \theta, \theta_\j,s;\hat E)
+O(\ep, \ep ^{\frac{\gamma}{2}})\right)&= \frac{ \tilde E'_\j-  E_\j}{ \ep^{1+\frac{1}{2} } }
\end{aligned}
\end{equation}
and the function $\bar \ell$ is a scaled version of the one given in
\eqref{eq:ell}:
$$
\bar \ell(\theta_\j; \hat E, e_\j;\ep)=
\sqrt{\frac{2}{a(\hat E, \ep)}(e_\j- U^{k_0,l_0,*}(\theta_\j;\hat E,\ep))}.
$$
The function $\bar \ell =O(1)$, but the important observation is that if we
take $\rho>0$ and we exclude a small region around the critical point $\tilde
\theta_\j$, (that is, for $\rho<\theta_\j-\tilde\theta_\j \le 2\pi-\rho$)
$\bar \ell$ never vanishes. In fact one has:
\begin{equation}\label{eq:boundlbar}
\bar \ell(\theta_\j; \hat E, e_\j;\ep)\ge d >0, \quad \mbox{for}\quad
\rho<\theta_\j-\tilde\theta_\j \le 2\pi-\rho .
\end{equation}

To have  non-degenerate solutions of equations \eqref{heteroscaled} it
suffices to assume:
\begin{equation}\label{determinant}
\det D_{\theta}\left( \begin{array}{c}
\frac{\partial}{\partial \hat \theta}\L^*_{k_0,l_0}  (\hat \theta, \theta_\j,s;E)\\
\pm \frac{a(\hat E)} {k_{0}^{\j}}  \bar \ell (\theta_\j;E;\ep)\frac{\partial}{\partial  \theta_\j}\L^*_{k_0,l_0} (\hat \theta, \theta_\j,s;E) \end{array}\right)
\neq 0 .
\end{equation}
Making explicit the derivatives in \eqref{determinant}  and separating in
blocks corresponding to $\hat \theta$ and $\theta _\j$, one obtains:
\begin{equation}
\pm \frac{a(\hat E)} {k_{0}^{\j}}  \left| \begin{array}{cc}
\frac{\partial^2}{\partial \hat \theta ^2}\L^*_{k_0,l_0} & \frac{\partial^2}{\partial \hat \theta \partial  \theta_\j}\L^*_{k_0,l_0}\\
\frac{\partial}{\partial \hat \theta}( \bar \ell\frac{\partial}{\partial  \theta_\j}\L^*_{k_0,l_0}) & \frac{\partial}{\partial \theta_\j}( \bar
\ell\frac{\partial}{\partial  \theta_\j}\L^*_{k_0,l_0})  \end{array}\right| (\hat \theta, \theta_\j,s; E, \ep)\neq 0
\end{equation}
which gives:
\begin{equation}\label{det}
\left| \begin{array}{cc}
\frac{\partial^2}{\partial \hat \theta ^2}\L^*_{k_0,l_0} & \frac{\partial^2}{\partial \hat \theta \partial\theta_\j}\L^*_{k_0,l_0}\\
 \bar \ell\frac{\partial ^2}{\partial \hat \theta \partial  \theta_\j}\L^*_{k_0,l_0} &  \bar \ell\frac{\partial ^2}{\partial \theta_\j ^2}\L^*_{k_0,l_0} - \frac
{(U^{k_0,l_0,*})'}{ a(\hat E)\  \bar \ell} \frac{\partial}{\partial  \theta_\j}\L^*_{k_0,l_0} \end{array}\right|(\hat \theta, \theta_\j,s; E,\ep)\neq 0 ,
\end{equation}
which, using that neither $\bar \ell$ nor $a(\hat E)$ vanish, is equivalent
to:
\begin{eqnarray}
&&\left(2(e_m-U^{k_0,l_0,*}) [ \frac{\partial^2}{\partial \hat \theta ^2}  \L^*_{k_0,l_0} \ \frac{\partial ^2}{\partial  \theta_\j^2}\L^*_{k_0,l_0}-
(\frac{\partial ^2}{\partial \hat \theta \partial  \theta_\j}\L^*_{k_0,l_0})^2]\right.-\nonumber \\
&& \left.\frac{\partial^2}{\partial \hat \theta ^2}\L^*_{k_0,l_0}  \frac{\partial}{\partial  \theta_\j}\L^*_{k_0,l_0} (U^{k_0,l_0,*})' \right)(\hat \theta,
\theta_\j,s; E,\ep)\neq 0 \label{eq:nondegenerateresonant}
\end{eqnarray}
This inequality (or \eqref{det}) constitutes  part of hypothesis {\bf H8},
and  is the generalization of the  non-degeneracy conditions  {\sl H5' } and
{\sl H5''} in \cite{DelshamsLS06a}. We call attention to the fact that
\eqref{det} takes a value for $\varepsilon =0$.

An equivalent formulation for this non-degeneracy conditions can be written
using the symplectic structure of the system.

Introducing the poisson brackets:
\begin{eqnarray*}
\{\hat F, \cdot\} =\frac{\partial }{\partial \hat \theta}\\
\{F_\j ,\cdot\}=  \bar \ell \frac{\partial }{\partial \theta_\j}
\end{eqnarray*}
we see that equations \eqref{hetero1} read:
\begin{eqnarray*}
\{\hat F, \L^*_{k_0,l_0}\} + O(\ep ,\ep ^{\frac{\gamma}{2}})=\frac{\hat E- \hat E'}{\ep}\\
\{F_\j ,\L^*_{k_0,l_0}\} + O(\ep ^{\frac{1}{2}} ,\ep ^{\frac{\gamma}{2}})=  \frac{\hat E_\j- \hat E_\j'}{\ep}
\end{eqnarray*}
and the non-degeneracy condition \eqref{det} becomes:
$$
\left| \begin{array}{cc}
\{\hat F, \{\hat F, \L^*_{k_0,l_0}\}\} & \{F_\j ,\{\hat F, \L^*_{k_0,l_0}\} \} \\
\{\hat F, \{\{F_\j ,\L^*_{k_0,l_0}\}\} & \{\{F_\j , \{\{F_\j , \L^*_{k_0,l_0}\}\}
\end{array}\right| (\hat \theta, \theta_\j,s; E,\ep)\ne 0 .
$$

\subsubsection{Heteroclinic connections between primary
tori and secondary tori close to homoclinic connections}

If condition \eqref{eq:nondegenerateresonant} is verified in the region
$\R_{k_0,l_0}^{[\le L]}$ we can guarantee the existence of heteroclinic
connections between the primary and secondary tori in this region. If we look
for this heteroclinic connections close to homoclinic ones, one can obtain a
more explicit sufficient condition  to have a solution of equations
\eqref{heteroscaled}.

The main idea is to solve equations  \eqref{heteroscaled} using the implicit
function theorem. The small parameters will be
$$\hat \delta=\frac{\hat E- \hat E'}{\ep} + O(\ep, \ep ^{\gamma/2})
, \quad \delta_\j =
\frac{ E_\j-  E_\j'}{\ep^{\frac{1}{2}+1}} + O(\ep ^{\frac{1}{2}}, \ep ^{\gamma/2}),
$$
and then equations \eqref{heteroscaled} read:

\begin{eqnarray*}
\frac{\partial}{\partial \hat \theta}\L^*_{k_0,l_0} (\hat \theta, \theta_\j,s;\hat E) &=& O(\hat \delta )\\
\frac{a(\hat E)} {k_{0}^{\j}}
\bar \ell(\theta_\j;E;\ep)   \frac{\partial}{\partial  \theta_\j}\L^*_{k_0,l_0} (\hat \theta, \theta_\j,s;\hat E)
&=& O(\delta_\j)
\end{eqnarray*}
and, using that, by \eqref{eq:boundlbar}, the function $\bar \ell$ never
vanishes neither does $a(\hat E)$, they are equivalent to:
\begin{eqnarray}
\frac{\partial}{\partial \hat \theta}\L^*_{k_0,l_0} (\hat \theta, \theta_\j,s;\hat E) &=& O(\hat \delta )\\
\frac{\partial}{\partial  \theta_\j}\L^*_{k_0,l_0} (\hat \theta, \theta_\j,s;\hat E)
&=& O(\delta_\j)
\end{eqnarray}

Therefore, the non-degeneracy condition which guarantees that these equations
have solutions close to the solutions of:
\begin{eqnarray}\label{homoclinic}
\frac{\partial}{\partial \hat \theta}\L^*_{k_0,l_0} (\hat \theta, \theta_\j,s;\hat E) &=& 0 \nonumber\\
\frac{\partial}{\partial  \theta_\j}\L^*_{k_0,l_0} (\hat \theta, \theta_\j,s;\hat E)
&=& 0
\end{eqnarray}
is simply:
\begin{equation}\label{nondegenhom}
{\rm det} \left |\frac{\partial ^{2}  \L^{*} _{k_0,l_0}}{\partial \theta ^{2}} (\theta,s;\hat E)\right |\ne 0 .
\end{equation}
holding in the region $\R_{k_0,l_0}^{[\le L]}$.
This is part of  the non-degeneracy conditions which constitute Hypothesis
{\bf H8}.

In summarizing, the Hypothesis {\bf H8}  consists in assuming inequalities
\eqref{eq:nondegeneratenores}, \eqref{eq:h8non-resonant}, \eqref{det},
\eqref{nondegenhom}.

It is important to note that the function $\L^*_{k_0,l_0}(\hat
\theta,\theta_\j)$ is the Poincar\'{e} function $\L^* (I^*(\hat
E),\vp-\omega(I^*(\hat E))s)$ after the linear change of variables
\eqref{eq:canvipoincare}. Therefore condition \eqref{nondegenhom} is
equivalent to condition \eqref{eq:h8non-resonant} that ensures that the
Poincar\'{e} function has non-degenerate critical points.

Any of the non-degeneracy simplified conditions \eqref{homoclinic},
\eqref{nondegenhom}, or equivalently \eqref{eq:h8non-resonant}, constitute
Hypothesis {\bf H8'} stated after Theorem \ref{main}, since  they are
sufficient conditions to ensure that the surface $\T_E'$ intersects
transversally $s_\ep (\T_E)$ for $|E-E'| =O(\ep ^{\frac{3}{2}+\frac{1}{2}})$.

\begin{remark}\label{rem:secondorder}
In the case of a second order resonance, one needs to take into account the
terms of order $\ep ^2$ in  equation \eqref{scateringDH}. Nevertheless, if
one looks for heteroclinic solutions close to homoclinic ones some  easy
computations show that  these heteroclinic connections exist  if equations
\eqref{homoclinic} have non-degenerate zeros, and this is also guaranteed by
condition \eqref{nondegenhom}.

\end{remark}


\subsection{Constructing chains of invariant tori. Contouring the
resonances of higher multiplicity. Formulation of the symbolic dynamics}

In this section, we will see how to put together the information we have
gathered on the scattering map and the KAM tori, and show that we can
construct largely arbitrary motions in action space. In particular, we can go
around double resonances and other effects of codimension 2.

We will prove the following result which clearly implies Theorem~\ref{main}
since it has the same hypothesis and clearly stronger conclusions.

\begin{theorem} \label{thm:path}
Let $H_\ep$ be a family of the form \eqref{modelsconsidered}. Assume that
$H_\ep$ satisfies all the hypothesis {\bf H1}--{\bf H8}. In particular, it is
$C^r$ for $r \ge r_0$.

Let $m_0$ be a sufficiently large number. Fix $\delta > 0$ sufficiently small
and consider the set $\I_\delta \subset \II^* \subset \I$ defined before
theorem \ref{averaging2} to verify condition L2. Then, there exists $\ep_0 >
0 $ such that for all $|\ep| \le \ep_0$, given any $C^1$ path $\gamma: [0,1]
\rightarrow \I_\delta$  in $\I_\delta$ there exists $x_\ep(t)$ a trajectory
of the flow generated by $H$ and a time reparameterization $\Psi_\ep$ (i.e. a
diffeomorphism  $\Psi_\ep: \real^+ \to [0,1]$) in such a way that
\begin{equation} \label{thereispath}
| I(x_\ep(t)) - \gamma( \Psi_\ep(t)) | \le C \ep^{1/2}
\end{equation}

\end{theorem}

Of course, Theorem~\ref{thm:path}  immediately implies Theorem~\ref{main}.
Clearly, the hypothesis of both theorems are the same and, for $\delta$
sufficiently small so that given any two points $I_-, I_+ \in \I^*$, we can
get a path contained in $\I_\delta$ which starts at a distance less than
$\delta$ from the $I_-$ and ends at a distance less than $\delta$ from $I
_+$. Applying Theorem~\ref{thm:path} to this path we obtain the statement of
Theorem~\ref{main} for $\delta + C \ep_0^{1/2}$.

As a corollary of the proof of Theorem~\ref{thm:path}, we obtain that it is
possible to construct orbits that are $\delta$ dense on invariant manifold
$\Lambda_\ep$ for $\ep$ small enough. These orbits also include excursions on
the stable and unstable manifolds. So that they are dense in a larger domain.
Some constructions of models with orbits dense on submanifolds appear also in
\cite{FontichM03}.

\begin{remark}
Note that we do not prescribe first the path and then state conditions on the
perturbations. We have identified conditions on the Hamiltonian that give the
simultaneous existence of trajectories that follow any path in $\II_\delta$.
\end{remark}

\begin{remark}
As we will see, the estimate in \eqref{thereispath}, is rather pessimistic
for most of the paths.  Indeed, except when the path is close to the resonant
region we can have a bound $C\ep$ in \eqref{thereispath}.
\end{remark}

\subsubsection{Proof of Theorem~\ref{thm:path} }

The proof of Theorem~\ref{thm:path} will consist in recalling all the
information that we have been gathering to construct a transition chain of
whiskered tori that follows the indicated path. Then, it will suffice to
invoke an obstruction argument that establishes that given

a transition chain of whiskered tori (i.e. a sequence of whiskered tori
$\Tau_i$ such that $W^u_{\Tau_i} \pitchfork W^s_{\Tau_{i+1}}$), there is an
orbit that follows the path.

Recall that we have shown that there is a normally hyperbolic invariant
manifold $\tilde \Lambda_\ep$.

We have  shown that under the nondegeneracy assumptions. {\bf H5}, {\bf H6}
we can define a scattering map in the region $\I^*$, which is of a size
independent of $\ep$. In this region, we could define the scattering map and
give explicit formulas for its leading behavior.

Independently of the scattering map, we have developed averaging theory and
obtained information about a geography of the resonances that appear when

averaging. It is important to note that the geography of resonances depends
only on the integrable flow. The perturbations activate some of them at the
order that we consider.

We recall that the set $\II_\delta$ was obtained by removing from the set
$\I^*$ (defined through hypotheses {\bf H3}---{\bf H8}) all the points at a
distance less than $\delta$ from either of
\begin{itemize}
\item[ a)] The set of  double resonances activated at order smaller than
    $m _0$ one of whose resonances is a secular resonance (i.e. a
    resonance of order $1$,$2$, see \eqref{multiplicity} and
    Definition~\ref{def:activated}).
\item[b)] The set of points in \eqref{degenerateresonances} for which the
    secular resonance is degenerate.
\end{itemize}

Note that by assumption {\bf H3} and {\bf H8}, the sets  involved in a), b)
above are the union of a finite number of codimension $2$ manifolds---b) will
be empty for quasi-convex Hamiltonians---. Hence, for sufficiently small
$\delta$, the set $\I_\delta$ will be connected. Note also that $\I_\delta$
is independent  of $\ep$ and that has a size of order $1$.

Recall that in Section~\ref{sec:analyzing} we have shown that the region
$\I_\delta$ can be covered by a collection  of KAM tori which are $\ep^{3/2}$
close to each other (as mentioned in Remark~\ref{rem:moreaveraging}, we could
have obtained a larger power of $\ep$ simply by averaging more times, which
requires to remove some more double resonances and assume more derivatives in
the model). We will refer to this collection as the \emph{scaffolding} since
the motions we construct consist on jumping from one element of the
scaffolding  to the next by the scattering map and moving along the element
for a while.

We have shown that, under the hypothesis {\bf H8}, we have that the image
under the scattering map of any of the tori constructed in
Section~\ref{averaged} intersects transversally all the other tori which are
at a distance smaller that  a quantity $O(\ep)$.

That is, if $\Tau$ is an invariant torus in $\tilde \Lambda_\ep$ -- hence a
whiskered torus in the whole phase space -- we have $\Tau \pitchfork_{\tilde
\Lambda_\ep} \Tau' $ for all other tori $\Tau'$ at a distance smaller than $C
\ep$.

Given a  $C^1$ path as in the conclusion of Theorem~\ref{thm:path}, we can
find a sequence $\{\Tau_i\}_{i = 0}^\infty$ of tori at a distance $O(\ep)$
from each other and from  the path $\gamma$ (recall that we have shown that
these tori are at a distance not more than $O(\ep^{3/2})$.) These tori
satisfy
\[
S_\ep(\Tau_i) \pitchfork_{\tilde \Lambda_\ep} \Tau_{i+1}.
\]
By Lemma 10.4 in \cite{DelshamsLS00}, we obtain that these invariant tori in
$\tilde \Lambda_\ep$ -- hence whiskered tori in the full phase space --
satisfy
\[
W^u_{\Tau_i} \pitchfork_{\tilde \Lambda_\ep} W^s_{\Tau_{i+1}}.
\]
That is, they constitute a transition chain.

In these circumstances, there are theorems that show that there are orbits
that follow the transition whole transition chain. One theorem particularly
well suited for our purposes is that of \cite{FontichM00} (See also some
extensions \cite{FontichM00, DelshamsLS00, DelshamsLS06b}.) Of course, there
are many versions of these results in the literature, but some of them
include the extra assumption that the Birkhoff normal form of the tori does
not contain some terms, or that the system is $C^\infty$ or that the
transition chain is finite. The paper \cite{FontichM00} does not have any of
these limitations and also does not need any assumptions on the topology of
the embedding of the torus. It applies just as well to chains in which some
of the tori are primary and other that are secondary.

\section{An example}\label{sec:example}

In this section, we present an explicit example where one can check it
verifies the conditions {\bf H1} to {\bf H8}. Consider the Hamiltonian:
\begin{multline}
 H(I_1,I_2,\varphi_1,\varphi_2,p,q,t,\ep) =\pm \left(\frac{p^2}{2}
 + \cos q -1\right)+ h(I_1,I_2)\\
 +\ep \cos q \ g(\varphi_1, \varphi_2,t)\label{eq:example}
\end{multline}
where
 $$
 h(I_1,I_2)= \Omega _1 \frac{I_1^2}{2}+ \Omega _2 \frac{I_2^2}{2},
$$
and
$$
g(\varphi_1, \varphi_2,t)= a_1\cos \varphi_1 +a_2 \cos \varphi
_2
+a_3 \cos (\varphi_1 + \varphi_2 - t).
$$
\begin{prop}
Assume that $a_0$, $a_1$, $a_2$, $\Omega _1$, $\Omega _2$, $\Omega _1 +\Omega
_2$, $4\Omega _1 +\Omega _2$ and $\Omega _1 +4\Omega _2$ are non zero. Then
Hamiltonian \eqref{eq:example} verifies hypotheses {\bf H1} to {\bf H8} of
Theorem \ref{main}.
\end{prop}

As we will see, the proof of this proposition is very explicit and we can
give a rather precise description of the geometric objects involved in the
construction. This proof also shows that there are other heteroclinic
connections which could be used to construct unstable orbits. These other
choices would lead, through similar calculations, to other regions of
parameters where Theorem \ref{main} applies.

\proof The first observation is that $g$ is a trigonometric polynomial in the
angles $\varphi _1$, $\varphi_2$, $t$, so it is clear that Hamiltonian
\eqref{eq:example} satisfies hypotheses {\bf H1} to {\bf H4}. The Hamiltonian
of one degree of freedom $P_{\pm}(p,q)=\pm\left(p^{2}/2+\cos q -1\right)$ is
the standard pendulum when we choose the $+$ sign, and its separatrix for
positive $p$ is given by:
\[
q_{0}(t)= 4\arctan e^{\pm t}, \quad p_{0}(t) = 2/{\cosh t}.
\]

An important feature of the Hamiltonian~\eqref{eq:example} is that the
5-dimen\-sional hyperbolic invariant manifold
\[
\tilde \Lambda=\{(0,0,I_1,I_2,\varphi_1, \varphi_2,s): (I_1,I_2,\varphi_1, \varphi_2,s) \in \real ^2\times
\torus^{3}\}
\]
is \emph{preserved} for $\ep\neq 0$: $p=q=0\Rightarrow\dot p=\dot q=0$.
However, in contrast with the example in \cite{Arnold64}, the perturbation
does not vanish on $\tilde \Lambda$. Indeed, the dynamics on $\tilde \Lambda$
is provided simply by the restriction of $H_{| \tilde \Lambda}$, which is a
$2$ and a half degrees of freedom Hamiltonian taking the form
$$h(I_1,I_2)+\ep  g(\varphi_1, \varphi_2,t).$$
However, for any  $I=(I_1,I_2)$, the $3$-dimensional whiskered tori
\[
\T^{0}_{I}=\{(0,0,I,\varphi_1, \varphi_2,s): (\varphi_1, \varphi
_2,s)\in
\torus^{3}\}
\]
are not preserved if $a_i\ne 0$, and the resonances activated at order one
are given by the equations $\omega_i=0$, $i=1,2,3$, where we introduce the
notation
$$\omega_1=\Omega_1 I_1, \quad \omega_2=\Omega_2, \quad
\omega_3=\Omega_1 I_1+\Omega_2 I_2 -1.$$ Therefore, \eqref{eq:example}
presents the large gap problem, because it has ``large gaps'' associated to
any of these resonances activated at order one (and also to the resonances
activated at order two that will be introduced later on).

The Melnikov potential \eqref{Melnikovpot} associated to the Hamiltonian
\eqref{eq:example} is given by
\[
L(\tau, I,\vp_1,\vp_2,s)=\frac{1}{2}\int_{-\infty}^{\infty}
p_{0}^{2}(\tau+\sigma)g(\vp_1+\Omega _1 I_1 \sigma, \vp_2+\Omega
_2 I_2 \sigma,s+\sigma) d\sigma,
\]
and computing the integrals by the residue theorem, we obtain
\[
L(\tau, I,\vp_1,\vp_2,s)= \sum_{i=1}^{3}A_{i}\cos(\vp_{i}-\omega_{i}\tau)
\]
where we introduce $\vp_3:=\vp_1+\vp_2-s$, and
\[
A_{i}=A_{i}(\omega_{i})=\frac{2\pi\omega_{i}}{\sinh(\pi\omega_{i}/2)}\,a_{i},
\qquad i=1,2,3.
\]
Since $\tau \in \real$, it can be written as
$$
L(\tau, I,\vp_1,\vp_2,s)=
\LL (I,\vp_1-\omega_1\tau,\vp_2-\omega _2\tau,\vp_3-\omega _3\tau)
$$
with
\[
\LL( I,\vp)=\sum_{i=1}^{3}A_{i}\cos \vp_{i}.
\]
Therefore
\[
\frac{\partial L}{\partial \tau}(\tau, I,\vp_1,\vp_2,s)
= \sum_{i=1}^{3}\omega_{i}A_{i}\cos(\vp_{i}-\omega_{i}\tau)
\]
so that, given $(I,\vp_1,\vp_2,s)$, the condition
$$
\frac{\partial L}{\partial \tau}(\tau,I,\vp_1,\vp_2,s) =0
$$
is equivalent to the search of critical points $\tau ^*$ of the map
\begin{equation}\label{eq:critical}
\tau \in \real \mapsto L(\tau, I,\vp_1,\vp_2,s)
=\sum_{i=1}^{3}A_{i}\cos(\vp_{i}-\omega_{i}\tau)
\end{equation}
that is, of the function $\LL$ restricted to the straight line in $\torus^3$:
\begin{equation}\label{eq:straighline}
{\mathbf R}= {\mathbf R}(I,\vp)= \{ \vp-\omega\tau, \ \tau \in \real \}.
\end{equation}
Fixing $I\in \real ^2$, the $8$ critical points of
$$
\vp\in\torus ^3\mapsto \LL(I,\vp)
$$
satisfy $\tau ^*=0$, as well as the points $(I,\vp)$ in
$$
\C (I)=\{ \vp\in\torus ^3, \sum_{i=1}^{3}\omega_{i}A_{i}\cos\vp_{i}=0,\
 \sum_{i=1}^{3}\omega_{i}^2A_{i}\sin\vp_{i}\neq 0\}.
$$
As a consequence, the search for  critical points of the map
\eqref{eq:critical} is equivalent to the search for intersections between the
straight line ${\mathbf R}(I,\vp)$ and the set $\C (I)$.

For Hamiltonian \eqref{eq:example} the equation of $\C(I)$ is simply
\begin{equation}\label{eq:crest}
\omega_1 A_1 \sin \vp_1+\omega_2 A_2 \sin \vp_2+\omega_3 A_3 \sin \vp_3 =0
\end{equation}
which is just $D \LL(\vp)\omega=0$ and defines, locally, the equation of a
regular surface in the angles $\vp=(\vp _1,\vp_2,\vp_3=\vp_1+\vp_2-s)$ as
long as $\omega\top D^{2}\LL(\vp)\omega\neq 0$ holds. We notice that for any
$I \in \real^2$, the points $\vp^*_M$ and $\vp^*_m$ where the Melnikov
potential $\LL$ reaches its maximum and minimum ($\vp^*_i=0$ or $\pi$) belong
to the set $\C (I)$, so there exist at least two zones contained in $\C (I)$
where this set behaves as a local regular surface $\C_M (I)$, $\C_m (I)$,
respectively, which will be called \emph{crests} in analogy with the case
when $\vp$ is two-dimensional (see \cite{DelshamsH09}).

Once the set $\C(I)$ is known to be formed at least by the two crests
$\C_m(I)$ and $\C_M(I)$, it is clear that, for any $\vp$, there exist several
possible intersections of the straight line ${\mathbf R}(I,\vp)$ given in
\eqref{eq:straighline} with the crests $\C_m(I)$ and $\C_M(I)$, parameterized
by several values $\tau ^*$ of the parameters $\tau$ which give rise to
several scattering maps.

{From} now on, we will choose \emph{only one} of these intersections, the
``first one'' with the crest $\C_M(I)$. Given $(I,\vp_1,\vp_2,s)$, we define
$\tau ^*(I,\vp_1,\vp_2,s)\allowbreak =\tau ^*_M(I,\vp)$ as the real number
$\tau$ with minimum absolute value $|\tau|$ among all $\tau$ satisfying:
$$
\vp-\omega\tau \in \C_M(I).
$$
To determine a domain of definition of $\tau^*$ in the variables $(I,\vp)$,
it suffices to check that the straight line ${\mathbf R }(I,\vp)$ intersects
transversally $\C_M(I)$, that is, that $\omega\top D^{2}\LL(\vp)\omega\neq 0$
which is exactly the inequality satisfied by $\C(I)$ and a fortiori by
$\C_M(I)$.

We can now choose the domain of definition $H_-= H_M$, where $\tau ^*$ is
continuous simply by taking $H_-$ as an appropriate neighborhood of
$\vp^*_M$, so that hypothesis {\bf H7} is fulfilled.

Recall that the reduced Poincar\'{e} function defined in \eqref{generator} is
$$
\L^* (I,\theta)= L(\tau ^*(I,\theta,0),I,\theta, 0)
 = \LL (I,\theta -\omega\tau ^*(I,\theta)).
$$
Given $(I,\theta)$, $\L^*(I,\theta)$ is the value of $\L$ on ${\mathbf
R}(I,\theta)\cap \C_M(I)$ and it is constant along ${\textbf R}(I,\theta)$ so
$\L^*(I,\theta)$ is well defined on $\C_M = \cup _{I\in \II^*} \C_M(I)$.

Recall that the scattering map written in coordinates $(I,\vp,s)$ takes the
form \eqref{scatteringformula}, which, in coordinates $(I, \theta =
\vp-\omega(I)s)$ becomes
$$
s_\varepsilon(I,\theta) =
(I+\ep \partial _{\theta}\L^{*}(I,\theta)+O(\ep^2), \
\theta-\ep \partial _I\L^{*}(I,\theta) + O(\ep^2), s).
$$

We will check the hypotheses in the non-resonant region and in the resonances
activated up to order two. There are three resonances activated at order one
in this model
\begin{eqnarray*}
\mathbf{R}_1 &=& \R _{1,0,0}=\{ (I_1, I_2), \ I_1=0\} \\
\mathbf{R}_2 &=& \R _{0,1,0}=\{ (I_1, I_2), \ I_2=0\}\\
\mathbf{R}_3 &=& \R _{1,1,-1}=\{ (I_1, I_2), \ \Omega_1 I_1+\Omega_2 I_2=1\}
\end{eqnarray*}
and four more activated at order two:
\begin{eqnarray*}
\mathbf{R}_4 &=& \R _{1,0,-1}=\{ (I_1, I_2), \ \Omega_1 I_1
=1\} \\
\mathbf{R}_5 &=& \R _{0,1,-1}=\{ (I_1, I_2), \ \Omega_2 I_2
=1\}\\
\mathbf{R}_6 &=& \R _{2,1,-1}=\{ (I_1, I_2), \ 2\Omega_1 I_1+\Omega_2 I_2=1\}\\
\mathbf{R}_7 &=& \R _{1,2,-1}=\{ (I_1, I_2), \ \Omega_1 I_1+2\Omega_2 I_2=1\}
\end{eqnarray*}

For $(I_1,I_2) $ in the non-resonant region, the condition to have
heteroclinic orbits between the KAM tori are given by
\eqref{eq:h8non-resonant}. In the resonant regions, one has to check
\eqref{homoclinic} and \eqref{nondegenhom}. In our example, one can easily
check that both conditions are implied by the conditions
\[
D\LL(\theta)\omega=0,\quad \omega^\top D^2 \LL(\theta)\omega\neq 0
\]
defining $\C(I)$, which are a fortiori satisfied by the crest $\C_M(I)$.

To check conditions {\bf H5}, {\bf H6} and {\bf H8} in $\mathbf{R}_i$ we
simply need to impose that $\Omega _1$, $\Omega _2$, $\Omega _1 +\Omega _2$,
$4\Omega _1 +\Omega _2$ and $\Omega _1 +4\Omega _2$ are non zero. Moreover,
the potential at the resonance $\mathbf{R}_2$ is given by
$$
U^1 (0,I_2, \vp_1)= a_1 \cos \vp_1
$$
and therefore hypothesis {\bf H6} is also verified. The study of the
potential in the other resonances $\mathbf{R}_i$ is analogous.

\section*{Acknowledgements}
R.L. has been supported by NSF grant DMS1162544 and he acknowledges the
hospitality of UPC for many visits and he was affiliated with Univ. of Texas
during part of the work.

A. Delshams and Tere M-Seara have  been partially supported by the Spanish
MINECO-FEDER Grants MTM2009-06973, MTM2012-31714 and the Catalan Grant
2009SGR859.

We thank CRM for hospitality during the program ``Stability and Instability
in Dynamical Systems'' in 2008. We also thank. M. Gidea, G. Huguet, V.
Kaloshin, P. Rold\'an, C. Sim\'o, D. Treschev for several discussions and
encouragements.


\begin{thebibliography}{FGKR11}

\bibitem[AKN88]{ArnoldKN88}
V.~Arnold, V.~Kozlov and A.~Neishtadt.
\newblock \emph{Dynamical Systems III}, volume~3 of \emph{Encyclopaedia Math.
  Sci.}
\newblock Springer, Berlin, 1988.

\bibitem[AM78]{AbrahamM78}
R.~Abraham and J.~E. Marsden.
\newblock \emph{Foundations of mechanics}.
\newblock Benjamin/Cummings Publishing Co. Inc. Advanced Book Program, Reading,
  Mass., 1978.
\newblock ISBN 0-8053-0102-X.
\newblock Second edition, revised and enlarged, With the assistance of Tudor
  Ra\c tiu and Richard Cushman.

\bibitem[{A}rn63]{Arnold63a}
V.~I. {A}rnol'd.
\newblock Proof of a theorem of {A}. {N}. {K}olmogorov on the invariance of
  quasi-periodic motions under small perturbations.
\newblock \emph{Russian Math. Surveys}, 18(5):9--36, 1963.

\bibitem[Arn64]{Arnold64}
V.~Arnold.
\newblock Instability of dynamical systems with several degrees of freedom.
\newblock \emph{Sov. Math. Doklady}, 5:581--585, 1964.

\bibitem[Arn78]{Arnold78}
V.~I. Arnold.
\newblock \emph{Mathematical methods of classical mechanics}.
\newblock Springer-Verlag, New York, 1978.
\newblock ISBN 0-387-90314-3.
\newblock Translated from the Russian by K. Vogtmann and A. Weinstein, Graduate
  Texts in Mathematics, 60.

\bibitem[BB02]{BertiB02}
M.~Berti and P.~Bolle.
\newblock A functional analysis approach to {A}rnold diffusion.
\newblock \emph{Ann. Inst. H. Poincar{\'e} Anal. Non Lin{\'e}aire},
  19(4):395--450, 2002.

\bibitem[BBB03]{BertiBB03}
M.~Berti, L.~Biasco and P.~Bolle.
\newblock Drift in phase space: a new variational mechanism with optimal
  diffusion time.
\newblock \emph{J. Math. Pures Appl. (9)}, 82(6):613--664, 2003.

\bibitem[BCV01]{BessiCV01}
U.~Bessi, L.~Chierchia and E.~Valdinoci.
\newblock Upper bounds on {A}rnold diffusion times via {M}ather theory.
\newblock \emph{J. Math. Pures Appl. (9)}, 80(1):105--129, 2001.

\bibitem[Ber10]{Bernard10}
P.~Bernard.
\newblock Arnold's diffusion: from the {\it a priori} unstable to the {\it a
  priori} stable case.
\newblock In \emph{Proceedings of the {I}nternational {C}ongress of
  {M}athematicians. {V}olume {III}}, pages 1680--1700. Hindustan Book Agency,
  New Delhi, 2010.

\bibitem[Car81]{Cary81}
J.~R. Cary.
\newblock Lie transform perturbation theory for {H}amiltonian systems.
\newblock \emph{Phys. Rep.}, 79(2):129--159, 1981.

\bibitem[CG94]{ChierchiaG94}
L.~Chierchia and G.~Gallavotti.
\newblock Drift and diffusion in phase space.
\newblock \emph{Ann. Inst. H. Poincar{\'e} Phys. Th{\'e}or.}, 60(1):1--144,
  1994.

\bibitem[CG98]{ChierchiaG98}
L.~Chierchia and G.~Gallavotti.
\newblock Erratum: Drift and diffusion in phase space.
\newblock \emph{Ann. Inst. H. Poincar{\'e} Phys. Th{\'e}or.}, 68:135, 1998.

\bibitem[Che08]{Cheng08}
C.-Q. Cheng.
\newblock Variational methods for the problem of {A}rnold diffusion.
\newblock In \emph{Hamiltonian dynamical systems and applications}, NATO Sci.
  Peace Secur. Ser. B Phys. Biophys., pages 337--365. Springer, Dordrecht,
  2008.

\bibitem[Che10]{Cheng10}
C.-Q. Cheng.
\newblock Variational construction of diffusion orbits for positive definite
  {L}agrangians.
\newblock In \emph{Proceedings of the {I}nternational {C}ongress of
  {M}athematicians. {V}olume {III}}, pages 1714--1728. Hindustan Book Agency,
  New Delhi, 2010.

\bibitem[{Che}12]{Cheng12}
C.-Q. {Cheng}.
\newblock {Arnold diffusion in nearly integrable Hamiltonian systems}.
\newblock \emph{ArXiv e-prints}, July 2012.

\bibitem[Chi79]{Chirikov79}
B.~Chirikov.
\newblock A universal instability of many-dimensional oscillator systems.
\newblock \emph{Phys. Rep.}, 52(5):264--379, 1979.

\bibitem[DGLS08]{DelshamsGLS08}
A.~Delshams, M.~Gidea, R.~de~la Llave and T.~M. Seara.
\newblock Geometric approaches to the problem of instability in {H}amiltonian
  systems. {A}n informal presentation.
\newblock In \emph{Hamiltonian dynamical systems and applications}, NATO Sci.
  Peace Secur. Ser. B Phys. Biophys., pages 285--336. Springer, Dordrecht,
  2008.

\bibitem[DH09]{DelshamsH09}
A.~Delshams and G.~Huguet.
\newblock Geography of resonances and {A}rnold diffusion in a priori unstable
  {H}amiltonian systems.
\newblock \emph{Nonlinearity}, 22(8):1997--2077, 2009.

\bibitem[DLS00]{DelshamsLS00}
A.~Delshams, R.~de~la Llave and T.~M. Seara.
\newblock A geometric approach to the existence of orbits with unbounded energy
  in generic periodic perturbations by a potential of generic geodesic flows of
  $\bf{T}^{\rm 2}$.
\newblock \emph{Comm. Math. Phys.}, 209(2):353--392, 2000.

\bibitem[DLS03]{DelshamsLS03a}
A.~Delshams, R.~de~la Llave and T.~M. Seara.
\newblock A geometric mechanism for diffusion in {H}amiltonian systems
  overcoming the large gap problem: announcement of results.
\newblock \emph{Electron. Res. Announc. Amer. Math. Soc.}, 9:125--134
  (electronic), 2003.

\bibitem[DLS06a]{DelshamsLS06a}
A.~Delshams, R.~de~la Llave and T.~M. Seara.
\newblock A geometric mechanism for diffusion in {H}amiltonian systems
  overcoming the large gap problem: heuristics and rigorous verification on a
  model.
\newblock \emph{Mem. Amer. Math. Soc.}, 179(844):viii+141, 2006.

\bibitem[DLS06b]{DelshamsLS06b}
A.~Delshams, R.~de~la Llave and T.~M. Seara.
\newblock Orbits of unbounded energy in quasi-periodic perturbations of
  geodesic flows.
\newblock \emph{Adv. Math.}, 202(1):64--188, 2006.

\bibitem[DLS08]{DelshamsLS08}
A.~Delshams, R.~de~la Llave and T.~M. Seara.
\newblock Geometric properties of the scattering map of a normally hyperbolic
  invariant manifold.
\newblock \emph{Adv. Math.}, 217(3):1096--1153, 2008.

\bibitem[DMM09]{DuToitMM09}
P.~{Du Toit}, I.~Mezi{\'c} and J.~Marsden.
\newblock Coupled oscillator models with no scale separation.
\newblock \emph{Phys. D}, 238(5):490--501, 2009.

\bibitem[Fen72]{Fenichel71}
N.~Fenichel.
\newblock Persistence and smoothness of invariant manifolds for flows.
\newblock \emph{Indiana Univ. Math. J.}, 21:193--226, 1971/1972.

\bibitem[Fen74]{Fenichel74}
N.~Fenichel.
\newblock Asymptotic stability with rate conditions.
\newblock \emph{Indiana Univ. Math. J.}, 23:1109--1137, 1973/74.

\bibitem[FGKR11]{FejozGKR11}
J.~Fejoz, M.~Guardia, V.~Kaloshin and P.~Rold{\'a}n.
\newblock Diffusion along mean motion resonance in the restricted planar
  three-body problem.
\newblock \emph{ArXiv e-prints}, 2011.

\bibitem[FM00]{FontichM00}
E.~Fontich and P.~Mart{\'i}n.
\newblock Differentiable invariant manifolds for partially hyperbolic tori and
  a lambda lemma.
\newblock \emph{Nonlinearity}, 13(5):1561--1593, 2000.

\bibitem[FM03]{FontichM03}
E.~Fontich and P.~Mart{\'i}n.
\newblock Hamiltonian systems with orbits covering densely submanifolds of
  small codimension.
\newblock \emph{Nonlinear Anal.}, 52(1):315--327, 2003.

\bibitem[GL06a]{GideaL06a}
M.~Gidea and R.~de~la Llave.
\newblock Arnold diffusion with optimal time in the large gap problem, 2006.
\newblock Preprint.

\bibitem[GL06b]{GideaL06b}
M.~Gidea and R.~de~la Llave.
\newblock Topological methods in the instability problem of {H}amiltonian
  systems.
\newblock \emph{Discrete Contin. Dyn. Syst.}, 14(2):295--328, 2006.

\bibitem[Hal97]{Haller97}
G.~Haller.
\newblock Universal homoclinic bifurcations and chaos near double resonances.
\newblock \emph{J. Statist. Phys.}, 86(5-6):1011--1051, 1997.

\bibitem[Hal99]{Haller99}
G.~Haller.
\newblock \emph{Chaos near resonance}.
\newblock Springer-Verlag, New York, 1999.
\newblock ISBN 0-387-98697-9.

\bibitem[Her83]{Herman85}
M.-R. Herman.
\newblock \emph{Sur les courbes invariantes par les diff{\'e}o\-mor\-phismes de
  l'anneau. {V}ol. 1}, volume 103 of \emph{Ast{\'e}risque}.
\newblock Soci{\'e}t{\'e} Math{\'e}matique de France, Paris, 1983.

\bibitem[HM82]{HolmesM82}
P.~Holmes and J.~Marsden.
\newblock Melnikov's method and {A}rnold diffusion for perturbations of
  integrable {H}amiltonian systems.
\newblock \emph{J. Math. Phys.}, 23(4):669--675, 1982.

\bibitem[HPS77]{HirschPS77}
M.~Hirsch, C.~Pugh and M.~Shub.
\newblock \emph{Invariant manifolds}, volume 583 of \emph{Lecture Notes in
  Math.}
\newblock Springer-Verlag, Berlin, 1977.

\bibitem[Kol54]{Kolmogorov79}
A.~N. Kolmogorov.
\newblock On conservation of conditionally periodic motions for a small change
  in {H}amilton's function.
\newblock \emph{Dokl. Akad. Nauk SSSR (N.S.)}, 98:527--530, 1954.
\newblock English translation in {\it Stochastic Behavior in Classical and
  Quantum Hamiltonian Systems (Volta Memorial Conf., Como, 1977)}, Lecture
  Notes in Phys., 93, pages 51--56. Springer, Berlin, 1979.

\bibitem[KZ12]{KaloshinZ12}
V.~{Kaloshin} and K.~{Zhang}.
\newblock {A strong form of Arnold diffusion for two and a half degrees of
  freedom}.
\newblock \emph{ArXiv e-prints}, December 2012.

\bibitem[Lla01]{Llave01a}
R.~de~la Llave.
\newblock A tutorial on {K}{A}{M} theory.
\newblock In \emph{Smooth ergodic theory and its applications (Seattle, WA,
  1999)}, pages 175--292. Amer. Math. Soc., Providence, RI, 2001.

\bibitem[LM88]{LochakM88}
P.~Lochak and C.~Meunier.
\newblock \emph{Multiphase Averaging for Classical Systems}, volume~72 of
  \emph{Appl. Math. Sci.}
\newblock Springer, New York, 1988.

\bibitem[LMM86]{LlaveMM86}
R.~de~la Llave, J.~M. Marco and R.~Moriy{\'o}n.
\newblock Canonical perturbation theory of {A}nosov systems and regularity
  results for the {L}iv\v sic cohomology equation.
\newblock \emph{Ann. of Math. (2)}, 123(3):537--611, 1986.

\bibitem[LT83]{LiebermanT83}
M.~A. Lieberman and J.~L. Tennyson.
\newblock Chaotic motion along resonance layers in near-integrable
  {H}amiltonian systems with three or more degrees of freedom.
\newblock In J.~C.~W. Horton and L.~E. Reichl, editors, \emph{Long-time
  prediction in dynamics (Lakeway, Tex., 1981)}, pages 179--211. Wiley, New
  York, 1983.

\bibitem[Mar13]{Marco13}
J.-P. Marco.
\newblock Generic hyperbolic properties of classical systems on the torus,
  2013.
\newblock Preprint.

\bibitem[Mat04]{Mather04}
J.~N. Mather.
\newblock Arnol\cprime d diffusion. {I}. {A}nnouncement of results.
\newblock \emph{J. Math. Sci. (N. Y.)}, 124(5):5275--5289, 2004.

\bibitem[Mey91]{Meyer91}
K.~R. Meyer.
\newblock Lie transform tutorial. {II}.
\newblock In K.~R. Meyer and D.~S. Schmidt, editors, \emph{Computer aided
  proofs in analysis (Cincinnati, OH, 1989)}, volume~28 of \emph{IMA Vol. Math.
  Appl.}, pages 190--210. Springer, New York, 1991.

\bibitem[Mos66]{Moser66a}
J.~Moser.
\newblock A rapidly convergent iteration method and non-linear partial
  differential equations. {I}.
\newblock \emph{Ann. Scuola Norm. Sup. Pisa (3)}, 20:265--315, 1966.

\bibitem[Ne{\u\i}81]{Neishtadt81}
A.~I. Ne{\u\i}shtadt.
\newblock Estimates in the {K}olmogorov theorem on conservation of
  conditionally periodic motions.
\newblock \emph{Prikl. Mat. Mekh.}, 45(6):1016--1025, 1981.

\bibitem[Pes04]{Pesin04}
Y.~B. Pesin.
\newblock \emph{Lectures on partial hyperbolicity and stable ergodicity}.
\newblock Zurich Lectures in Advanced Mathematics. European Mathematical
  Society (EMS), Z{\"u}rich, 2004.
\newblock ISBN 3-03719-003-5.

\bibitem[P{\"o}s82]{Poschel82}
J.~P{\"o}schel.
\newblock Integrability of {H}amiltonian systems on {C}antor sets.
\newblock \emph{Comm. Pure Appl. Math.}, 35(5):653--696, 1982.

\bibitem[P{\"o}s01]{Poschel01}
J.~P{\"o}schel.
\newblock A lecture on the classical {KAM} theorem.
\newblock In \emph{Smooth ergodic theory and its applications (Seattle, WA,
  1999)}, volume~69 of \emph{Proc. Sympos. Pure Math.}, pages 707--732. Amer.
  Math. Soc., Providence, RI, 2001.

\bibitem[PT07]{PiftankinT07}
G.~N. Piftankin and D.~V. Treshch{\"e}v.
\newblock Separatrix maps in {H}amiltonian systems.
\newblock \emph{Russian Math. Surveys}, 62(2):219--322, 2007.

\bibitem[Sal04]{Salamon04}
D.~A. Salamon.
\newblock The {K}olmogorov-{A}rnold-{M}oser theorem.
\newblock \emph{Math. Phys. Electron. J.}, 10:Paper 3, 37 pp. (electronic),
  2004.

\bibitem[Ten82]{Tennyson82}
J.~Tennyson.
\newblock Resonance transport in near-integrable systems with many degrees of
  freedom.
\newblock \emph{Phys. D}, 5(1):123--135, 1982.

\bibitem[Thi97]{Thirring97}
W.~Thirring.
\newblock \emph{Classical mathematical physics}.
\newblock Springer-Verlag, New York, third edition, 1997.
\newblock ISBN 0-387-94843-0.
\newblock Dynamical systems and field theories, Translated from the German by
  Evans M. Harrell, II.

\bibitem[Tre12]{Treschev12}
D.~Treschev.
\newblock Arnold diffusion far from strong resonances in multidimensional {\it
  a priori} unstable {H}amiltonian systems.
\newblock \emph{Nonlinearity}, 25(9):2717--2757, 2012.

\bibitem[Zeh75]{Zehnder75}
E.~Zehnder.
\newblock Generalized implicit function theorems with applications to some
  small divisor problems/{I}.
\newblock \emph{Comm. Pure Appl. Math.}, 28:91--140, 1975.

\bibitem[Zeh76a]{Zehnder76}
E.~Zehnder.
\newblock Generalized implicit function theorems with applications to some
  small divisor problems/{II}.
\newblock \emph{Comm. Pure Appl. Math.}, 29:49--111, 1976.

\bibitem[Zeh76b]{Zehnder76b}
E.~Zehnder.
\newblock Moser's implicit function theorem in the framework of analytic
  smoothing.
\newblock \emph{Math. Ann.}, 219(2):105--121, 1976.

\bibitem[Zha11]{Zhang11}
K.~Zhang.
\newblock Speed of {A}rnold diffusion for analytic {H}amiltonian systems.
\newblock \emph{Invent. Math.}, 186(2):255--290, 2011.

\end{thebibliography}

\def\cprime{$'$} \def\cprime{$'$}

\end{document}